\documentclass[11pt,a4paper]{article}

\usepackage[a4paper,margin=30mm]{geometry}
\usepackage[T1]{fontenc}
\usepackage[utf8]{inputenc}
\usepackage{mathpazo}
\usepackage{mathrsfs}
\usepackage{microtype}
\usepackage{amsmath,amssymb,amsthm,mathtools}

\usepackage{enumitem}
\usepackage{xcolor}
\usepackage{hyperref}

\hypersetup{
  colorlinks=true,
  linkcolor=blue!55!black,
  citecolor=blue!55!black,
  urlcolor=blue!55!black,
  pdftitle={Liouville systems II: Critical points at infinity in the resonant case},
  pdfauthor={Mohameden Ahmedou, Jonas Fey, Zhengni Hu, Lei Zhang},
  pdfsubject={Critical points at infinity and noncompact ends of a pseudogradient flow for resonant Liouville systems},
  pdfkeywords={Liouville systems, critical points at infinity, resonant parameters, pseudogradient flows, blow-up analysis}
}

\usepackage{cite}
\allowdisplaybreaks
\numberwithin{equation}{section}
\newtheorem{theorem}{Theorem}[section]
\newtheorem{lemma}[theorem]{Lemma}
\newtheorem{proposition}[theorem]{Proposition}
\newtheorem{corollary}[theorem]{Corollary}
\theoremstyle{definition}
\newtheorem{definition}[theorem]{Definition}
\theoremstyle{remark}
\newtheorem{remark}[theorem]{Remark}
\newcommand{\oH}{\mathring{\mathrm H}^{1}}
\newcommand{\ff}{\mathcal F^N_{\rho^*,h}}
\newcommand{\rg}{\mathfrak g}
\newcommand{\rf}{\mathrm f}
\newcommand{\fJ}{\mathfrak c}
\title{Liouville systems II:\\Critical points at infinity in the resonant case}
\author{Mohameden Ahmedou\thanks{Mathematisches Institut, Justus-Liebig-Universit\"at Giessen, Germany.}
\and Jonas Fey\thanks{Mathematisches Institut, Justus-Liebig-Universit\"at Giessen, Germany.}
\and Zhengni Hu\thanks{Department of Mathematics, Shanghai Jiao Tong University, Shanghai, China.}
\and Lei Zhang\thanks{Department of Mathematics, University of Florida, Gainesville, FL, USA.}}
\date{November 22, 2026}
\begin{document}
\maketitle
\begin{abstract}
For a Liouville system at a resonant parameter, we  construct a
bounded, locally Lipschitz pseudogradient, which classifies the non-compact ends of its flow at finite energy levels.  If the reduced functional is
Morse and the leading coefficient is nonzero at every balanced critical
pair, each unbounded forward trajectory with energy bounded below
converges in the bubbling parameters to a pair with negative
coefficient. Conversely, each such pair, up to a  permutation, is realized
by an open set of escaping trajectories.
\end{abstract}



\tableofcontents

\section{Introduction}

Under the hypotheses \ref{II:assump:H1}--\ref{II:assump:H2} stated below, we classify the bubbling configurations that a descending trajectory can approach at finite energy when a Liouville system is at resonance. 

Let $(M,g)$ be a closed Riemann surface of unit volume. We consider the generalized Liouville system
\begin{equation}\label{II:eq:L-System}
 -\Delta_g u_i=
 \sum_{j=1}^n a_{ij}\rho_j
 \left(\frac{h_je^{u_j}}{\int_Mh_je^{u_j}\,dV_g}-1\right),
 \qquad i\in I:=\{1,\ldots,n\},
\end{equation}
where $h_j\in C^4(M)$ is positive, $\rho_j>0$, and $A=(a_{ij})$ satisfies \ref{II:assump:H1}--\ref{II:assump:H2} stated below. Set  $A^{-1}=(a^{ij})$. Fixing the component means gives the space
\[
 \mathcal H:=(\oH)^n,
 \qquad
 \oH:=\left\{v\in H^1(M):\int_Mv\,dV_g=0\right\},
\]
with Dirichlet norm
\[
 \|u\|_{\mathcal H}^2:=\sum_{i=1}^n\int_M|\nabla_gu_i|_g^2\,dV_g.
\]
 The corresponding Euler--Lagrange functional of \eqref{II:eq:L-System} is
\begin{equation}\label{II:def:E-functional}
 J_\rho(u)=\frac12\sum_{i,j=1}^na^{ij}
 \int_M\nabla_gu_i\cdot\nabla_gu_j\,dV_g
 -\sum_{i=1}^n\rho_i\log\int_Mh_ie^{u_i}\,dV_g.
\end{equation}

We assume the Lin--Zhang interaction hypotheses
\begin{enumerate}[label=(H\arabic*), ref=(H\arabic*)]
 \item \label{II:assump:H1} $A$ is symmetric, entrywise nonnegative, irreducible, and
 invertible;
 \item \label{II:assump:H2} $a^{ii}\le0$, $a^{ij}\ge0$ for $i\ne j$, and
 $\sum_{j=1}^na^{ij}\ge0$ for every $i$.
\end{enumerate}

For each integer $N\ge1$ define the critical  hypersurface by
\begin{equation*}
 \Gamma_N:=\{\rho\in(\mathbb R_+)^n:\Lambda_{I,N}(\rho)=0\}, 
\end{equation*}
where $ \Lambda_{I,N}(\rho):=
 4\sum_{i=1}^n\frac{\rho_i}{2\pi N}
 -\sum_{i,j=1}^na_{ij}
 \frac{\rho_i}{2\pi N}\frac{\rho_j}{2\pi N}.$
 If   $\rho\in \Gamma_N$ for some $N\in\mathbb N $, we call $\rho$  a resonant parameter.
 
For parameters in compact subsets of  $(0, +\infty)^n  \setminus \cup_{N\in \mathbb N}\Gamma_N$, Lin and Zhang \cite{LinZhang2011} established uniform bounds for solutions of \eqref{II:eq:L-System} and computed the Leray--Schauder degree in each chamber  
\[
\mathcal O_N
:=
\left\{
\rho\in(\mathbb R_+)^n:
8\pi(N-1)\sum_{i=1}^n\rho_i
<
\sum_{i,j=1}^n a_{ij}\rho_i\rho_j
<
8\pi N\sum_{i=1}^n\rho_i
\right\}\qquad \text{ for } N\in \mathbb N .
\]
Moreover,
in the one-bubble case, Lin and Zhang obtained sharp estimates for
the blow-up location and for the rate at which the parameters of a
bubbling sequence approach the critical hypersurface in 
\cite{LZone}.
For multiple bubbles, Huang and Zhang \cite{HuangZhang} proved comparability of the concentration scales, obtained sharp parameter estimates, and derived relations among the coefficients at distinct blow-up points. Cheng, Li and Zhang constructed bubbling solutions on a flat torus \cite{CLZconstruction} and, under additional hypotheses, proved local uniqueness and nondegeneracy for such solutions \cite{CLZ}.

We now fix $\rho^*\in\Gamma_N$ and study the flow of $J_{\rho^*}$. Our previous paper \cite{AFHZI} gives an orthogonal decomposition of blow-up solutions into an $N$-bubble approximate solution and a remainder, and proves compactness of the solution set at $\rho^*$ under the hypotheses used in this paper. Resonance affects the variational problem differently. Along a balanced bubbling ray, $\Lambda_{I, N}(\rho^*)=0$ cancels the terms in the quadratic energy and the logarithmic integrals that grow with the common logarithmic scale. The Dirichlet norm diverges, the energy has a finite limit, and the gradient tends to zero (see Theorems~\ref{II:thm:fixed-balanced-ray} and~\ref{II:thm:mixed-ray}, and the proof of Theorem~\ref{II:thm:ambient-rays}). The direction of descent along the common scale determines whether the ray defines an end of the flow.

For a balanced critical  pair $(p^*,r^*)$ (see  Definition \ref{def:S-balance}), let
$\mathcal L_\infty(p^*,r^*)$ be the leading coefficient in the parameter
expansion of \cite{AFHZI}, defined in \eqref{II:def:L-infty}. For any sequence of  bubbling solutions, a nonzero coefficient $\mathcal L_\infty(p^*, r^*)$ has the opposite sign to
$\Lambda_{I,N}(\rho^k)$ for large $k$. Thus, an approach from below
requires $\mathcal L_\infty(p^*,r^*)<0$. At the resonance parameter, the same
coefficient determines the sign of the energy derivative with
respect to the common scale along the associated bubbling ray.

\begin{theorem}\label{II:thm:characterization-cpi}
Fix $N\ge1$ and $\rho^*=(\rho_1^*, \ldots, \rho_n^*)\in\Gamma_N$. Assume
\hyperref[II:assump:H1]{\textup{(H1)}}--
\hyperref[II:assump:H2]{\textup{(H2)}}, and $\ff$ is Morse and
$\mathcal L_\infty(p,r^*)\ne0$ at every balanced critical pair.

There exists a bounded, locally Lipschitz vector field
$W\colon\mathcal H\to\mathcal H$ such that
\begin{equation*}
DJ_{\rho^*}(u)[W(u)]<0
\qquad\text{when }DJ_{\rho^*}(u)\ne0.
\end{equation*}
Let $u(s)$ be a forward trajectory of $W$.

\begin{enumerate}[label=\textup{(\roman*)}]
\item
If $u(s)$ is bounded in $\mathcal H$, then its orbit is relatively compact
and
\begin{equation*}
\operatorname{dist}_{\mathcal H}
\bigl(u(s),\operatorname{Crit}(J_{\rho^*})\bigr)
\to0
\qquad\text{as }s\to\infty.
\end{equation*}

\item
Suppose that $u(s)$ is unbounded and
$
\inf_{s\ge0}J_{\rho^*}(u(s))>-\infty.
$
Then $u(s)$ eventually remains in a full $N$-bubble chart. After one
fixed relabeling of the bubbles, its modulation parameters satisfy
\begin{equation*}
L(s):=\frac1N\sum_{t=1}^N\log\lambda_t(s)
\to+\infty,
\qquad
(p(s),r(s))\to(p^*,r^*),
\end{equation*}
where $\lambda_1(s),\ldots,\lambda_N(s)$ are its bubble
scales, and  $(p^*,r^*)$ is a balanced critical pair such that
\begin{equation*}
\nabla\ff(p^*)=0,
\qquad
\mathcal L_\infty(p^*,r^*)<0.
\end{equation*}
Moreover,
\begin{equation*}
J_{\rho^*}(u(s))
\to
C(v,N)-\frac1N\ff(p^*)
\qquad\text{as }s\to\infty,
\end{equation*}
where $C(v,N)$ depends only on the normalized entire profile $v$, the
fixed resonant parameter $\rho^*$, and $N$, and is independent of
$(p^*,r^*)$.
\end{enumerate}

Conversely, every balanced critical pair $(p^*,r^*)$ satisfying
$\mathcal L_\infty(p^*,r^*)<0$ admits a forward-invariant open set of
initial data whose trajectories escape with limiting configuration
$(p^*,r^*)$. Therefore, in the sense of
Definition~\ref{II:def:critical-point-infinity}, the critical points at
infinity are in bijection with
\begin{equation*}
\mathcal K_\infty^-/\mathcal S_N,
\end{equation*}
where $\mathcal S_N$ acts by simultaneously permuting the centers and
their relative log-scales.
\end{theorem}

In this paper, we employ Bahri’s framework of critical points at infinity \cite{Bahri1989}. The application of this framework to critical elliptic problems \cite{BahriCoron1988,BahriLiRey1995} and to resonant scalar Liouville-type equations \cite{AhmedouBenAyed2017, MBL2017, ABAMorse, AHW, AHWZ} motivates the present construction. For the scalar Dirichlet problem, Ahmedou, Ben Ayed, and Lucia \cite{MBL2017} minimize the energy in the normal directions and utilize derivatives with respect to the bubble parameters to construct a pseudogradient. In contrast, for the Liouville system \eqref{II:eq:L-System}, the quadratic form of $A^{-1}$ is indefinite. Consequently, we solve the normal equation using a uniformly bounded inverse. Another challenge arises because varying the component amplitudes alters the mass balance, while the same relative scales must serve every component. We address this by solving the coupled equations for the amplitudes and relative scales at finite scale. Subsequently, we compute the leading energy derivative with respect to the common scale in both decay regimes. The sign of this derivative determines the direction of the flow near a balanced ray. Furthermore, we establish that if a trajectory enters arbitrarily small neighborhoods of a balanced pair with a negative coefficient at arbitrarily large times, it eventually remains in any fixed neighborhood of that pair.

The paper is organized as follows. Section \ref{sec:2} introduces the resonant data, balanced configurations, and bubbling charts. Section \ref{sec:3} develops the normal reduction and finite-scale centering. Section \ref{sec:4} analyzes the transverse equations and the common-scale derivative. Section \ref{sec:5} constructs the pseudogradient and establishes the classification of critical points at infinity.

\section{Preliminary}\label{sec:2}
\subsection{Resonance and balanced configurations}\label{II:sec:resonant-data}
Fix $\rho^*\in\Gamma_N$, set $\sigma_i=\rho_i^*/(2\pi N)$ and $D_\sigma=\operatorname{diag}(\sigma_i)$, and set 
\[
 m=A\sigma,\quad q=m-2\mathbf1,\quad \widehat m=\min_i m_i,\quad
 \widehat I=\{i:m_i=\widehat m\}, 
\]
where $\mathbf 1=(1,\ldots,1)\in \mathbb R^n$.

Under \ref{II:assump:H1}--\ref{II:assump:H2}, positivity of $\sigma$ and resonance give the strict
proper-subsystem mass inequalities \cite[Theorem~2.2]{LinZhang2011}.
The existence and symmetry theorems of Chipot--Shafrir--Wolansky \cite{CSW1997}, stated as Theorems~A and B in \cite{LinZhang2010}, together with radial uniqueness \cite[Theorem~1.1]{LinZhang2010}, yield the entire profile
$v$ centered at the origin and normalized by $\max_i v_i(0)=0$. It satisfies
\begin{equation*}
 -\Delta v_i=\sum_j a_{ij}e^{v_j},\qquad
 (2\pi)^{-1}\int_{\mathbb R^2}e^{v_i}=\sigma_i,\qquad
 v_i(y)=-m_i\log|y|+\fJ_i+O(|y|^{-(\widehat m-2)}) .
\end{equation*}
Here, $m_i>2$ by finite mass and the radial asymptotics of  \cite[Section~2]{AFHZI}.
The Pohozaev identity from \cite{CK1995} shows  $\sum_i\sigma_i(4-m_i)=0$ at resonance. Thus either $m_i=4$ for all $i$ or
$2<\widehat m<4$. We set
\begin{equation*}
 \Psi_{\widehat m}(\lambda)=
 \begin{cases}\lambda^{-(\widehat m-2)},&\widehat m<4,\\
 \lambda^{-2}\log\lambda,&\widehat m=4.\end{cases}
\end{equation*}

The Green function $G^g(\cdot,p)$ is normalized by
\[
-\Delta_g G^g(\cdot,p)=\delta_p-1,
\qquad
\int_M G^g(x,p)\,dV_g(x)=0.
\]
Here $\delta_p$ denotes the Dirac measure at $p$. 
Choose smoothly centered isothermal coordinates $(y_p, U(p))$ and a radial cutoff
$\chi_p$, equal to one on $|y_p|<r_0/2$ and supported on $|y_p|<r_0$.
The metric is $e^{\widehat\varphi_p(y)}|dy|^2$, with
$\widehat\varphi_p(0)=0$, $\nabla\widehat\varphi_p(0)=0$ and
$\Delta\widehat\varphi_p(0)=-2K_g(p)$.
Denote  $\varphi_p(x)=\widehat\varphi_p(y_p(x))$. Define the regular part of the Green function and the Robin function as follows: 
\[
 H^g(x,p)=G^g(x,p)+(2\pi)^{-1}\chi_p(x)\log|y_p(x)|,
 \qquad R^g(p)=H^g(p,p).
\]
Let
\[
\Delta
:=
\left\{
(p_1,\ldots,p_N)\in M^N:
p_s=p_t \text{ for some } s\ne t
\right\}.
\]
For
$
p=(p_1,\ldots,p_N)\in M^N\setminus\Delta,$
define
\begin{equation*}
\ff(p)
:=
\sum_{t=1}^N\sum_{i=1}^n \rho_i^*
\left[
\log h_i(p_t)
+\pi m_i
\left(
R^g(p_t)
+\sum_{\substack{s=1\\ s\ne t}}^N G^g(p_t,p_s)
\right)
\right].
\end{equation*}
Set
\begin{equation*}
 \widetilde G_t(x;p)=H^g(x,p_t)+\sum_{s\ne t}G^g(x,p_s),\qquad
 \rf_{p,t}^i(x)=h_i(x)e^{2\pi m_i\widetilde G_t(x;p)} .
\end{equation*}
It follows that 
\[
\nabla_{p_t}\ff=\sum_i\rho_i^*\nabla_x\log\rf_{p,t}^i(p_t).
\]

Set
\[
\ell_t:=\log\lambda_t,\qquad
L:=\frac1N\sum_{t=1}^N\ell_t,\qquad
\delta_t:=\ell_t-L.
\]
Then $\delta=(\delta_t)\in E_0:=\{z\in\mathbb R^N:\sum_tz_t=0\}$. Fix an
orthonormal basis $\{e^\alpha\}_{\alpha=1}^{N-1}$ of $E_0$ and denote
\[
\delta_t=\sum_{\alpha=1}^{N-1}e_t^\alpha r_\alpha.
\]
Also set
\[
\lambda:=e^L,\qquad
\varepsilon:=e^{-L},\qquad
\varepsilon_t:=\lambda_t^{-1}=e^{-L-\delta_t},
\qquad
\Psi_L:=\Psi_{\widehat m}(e^L).
\]
For $N=1$, $E_0=\{0\}$, so no relative scale parameter is needed.

\begin{definition}\label{def:S-balance}
Let $p\in M^N\setminus\Delta$ and $r^*\in\mathbb R^{N-1}$, and set
$
\delta_t^*:=\sum_{\alpha=1}^{N-1}e_t^\alpha r_\alpha^*.
$
The pair $(p,r^*)$ is \emph{balanced} if, for every $i$,
\[
\rf_{p,t}^i(p_t)e^{q_i\delta_t^*}
\]
is independent of $t$. Since $q_i=m_i-2>0$, for each $p$ there is at most one
$r^*$ satisfying this condition. If $N=1$, then $E_0=\{0\}$ and no
parameter $r^*$ is needed.
\end{definition}

For $\widehat m=4$, define
\begin{align}
T_{i,t}^{\rm raw}
&:=
\Delta_{g,x}\log\rf_{p,t}^i(p_t)
+\left|\nabla_x\log\rf_{p,t}^i(p_t)\right|^2
-2K_g(p_t),
\nonumber\\
\mathcal L_\infty(p,r^*)
&:=
\frac1N\sum_{i=1}^n\sum_{t=1}^N
e^{\fJ_i-2\delta_t^*}T_{i,t}^{\rm raw}.
\label{II:def:L-infty}
\end{align}
For $2<\widehat m<4$, we define $\mathcal L_\infty$ by taking the finite
part in \eqref{II:eq:mixed-coefficient}. Both definitions agree with
\cite{AFHZI}.

We assume:
\begin{enumerate}[label=(A\arabic*),ref=(A\arabic*)]
\item\label{II:item:a1}
$\ff$ is Morse on $M^N\setminus\Delta$.
\item\label{II:item:a2}
$\mathcal L_\infty(p,r^*)\ne0$ whenever
$\nabla\ff(p)=0$ and $(p,r^*)$ is balanced.
\end{enumerate}

Set
\[
\mathcal K_\infty^-
:=
\left\{
(p,r^*):
\nabla\ff(p)=0,\ (p,r^*)\text{ is balanced},\
\mathcal L_\infty(p,r^*)<0
\right\}.
\]
The group $\mathfrak S_N$ acts on this set by simultaneously permuting
the centers and the entries of $\delta^*$.

For $N\ge2$, Proposition~5.1 of \cite{AFHZI} gives a uniform positive
lower bound for the distances between the centers of every critical
point of $\ff$. Hence these critical points lie in a compact subset of
$M^N\setminus\Delta$. This is immediate for $N=1$. By
Assumption~\ref{II:item:a1}, there are only finitely many such points.
For each of them, at most one $r^*$ satisfies the balance condition.

For $u=(u_1,\ldots,u_n)$, set
\[
H_i(u):=\int_M h_i e^{u_i}\,dV_g,\qquad
D_i(u):=\frac{\rho_i^*h_i e^{u_i}}{H_i(u)},\qquad
d\mu_i(u):=\frac{D_i(u)}{\rho_i^*}\,dV_g.
\]
For $b=(b_{i,t})$ near zero, define
\begin{align}
\alpha^{i,j,t}
&:=\delta_{ij}+a^{ji}b_{i,t},
\qquad
\sigma_{i,t}(b)
:=
\frac1{2\pi}\int_{\mathbb R^2}
e^{\sum_j\alpha^{i,j,t}v_j},
\nonumber\\
\rg_{p,b}^i(x)
&:=
\exp\left(
2\pi\sum_{j=1}^n\sum_{s=1}^N
m_j(\alpha^{i,j,s}-\delta_{ij})H^g(x,p_s)
\right),
\nonumber\\
\mathfrak M_{i,t}
&:=
2\pi\sigma_{i,t}(b)
e^{-\sum_j\alpha^{i,j,t}\fJ_j}
(\rf_{p,t}^i\rg_{p,b}^i)(p_t)
\lambda_t^{q_i+\sigma_i b_{i,t}},
\qquad
\mathfrak S_i
:=
\sum_{t=1}^N\mathfrak M_{i,t}.
\label{II:def:full-local-mass}
\end{align}
Set 
\[
\Omega_{i,t}:=\frac{\mathfrak M_{i,t}}{\mathfrak S_i},\qquad
\theta_{i,t}:=1-N\Omega_{i,t},\qquad
\Delta_\Omega:=\max_{i,t}|\theta_{i,t}|.
\]
Thus, $\theta_{i,t}=0$ for all $t$ precisely when
$
\mathfrak M_{i,1}=\cdots=\mathfrak M_{i,N}.
$

For  $b_t=(b_{1,t},\ldots,b_{n,t})$, set
\[
b^0:=\frac1N\sum_{t=1}^N b_t,\qquad
b^\alpha:=\sum_{t=1}^N e_t^\alpha b_t,\qquad
c^\alpha:=Lb^\alpha,\qquad
\beta^0:=Lb^0,
\]
and
\[
|b|:=\max_{i,t}|b_{i,t}|,\qquad
K:=A^{-1},\qquad
\xi_t:=D_\sigma b_t.
\]

\subsection{Projected bubbles and neighborhood at infinity}\label{II:sec:bubbling-chart}
For $f\in L^1(M)$, set
$
\overline f:=\frac 1 {\operatorname{vol}_g(M)}\int_M f\,dV_g.
$
Let $p$ range over a compact subset of $M^N\setminus\Delta$ on which the
supports of the cutoffs are disjoint, and assume that $|\ell_t-L|$ is
bounded. Set $\chi_t:=\chi_{p_t}$. For
$Y=\lambda_t y_{p_t}(x)$, define
\begin{equation*}
-\Delta_gPW_{i,t}=V_{i,t}-\overline{V_{i,t}},
\qquad
\overline{PW_{i,t}}=0,
\end{equation*}
where  $V_{i,t}
:=
\chi_t e^{-\varphi_{p_t}}\lambda_t^2e^{v_i(Y)}$. Set
\begin{equation*}
PW_{i,t}^*:=\sum_j a_{ij}PW_{j,t},
\qquad
U_{z,i}:=\sum_{t=1}^N
\bigl(PW_{i,t}^*+b_{i,t}PW_{i,t}\bigr),
\end{equation*}
and
\begin{equation*}
z:=(p,\ell,b),
\qquad
E_z:=\bigl(\operatorname{range}DU_z\bigr)^\perp.
\end{equation*}

Choose $C>1$, $\eta>0$, and $\epsilon_0>0$ uniformly for all balanced
critical pairs. 
For a critical point $p^*\in M^N\setminus\Delta$, set
\[
M_{p^*,\tau}
:=
\left\{
p\in M^N\setminus\Delta:
\max_t d_g(p_t,p_t^*)\le\tau,\quad
U(p_s)\cap U(p_t)=\varnothing\ \text{for }s\ne t
\right\}.
\]
For $\epsilon_0>0$, define the neighborhood of the potential critical point at infinity:
\begin{align*}
V^\eta_{\rho^*}(N,\epsilon_0)
:=
\Bigl\{u\in\mathcal H:\;&
\|\nabla J_{\rho^*}(u)\|<\epsilon_0;\\
&\exists\,\lambda\in(\mathbb R_+)^N,\quad
\max_t\lambda_t^{-1}<\epsilon_0,\quad
\max_{s,t}\frac{\lambda_t}{\lambda_s}<C;\\
&\exists\,p\in M_{p^*,\tau/2},\quad
d_g(p_s,p_t)>2\eta\ \text{for }s\ne t;\\
&\max_i\|u_i-U_{i,p,\lambda,\mathbf 0}\|<\epsilon_0
\Bigr\}.
\end{align*}
Here $\epsilon_0$ is smaller than the radius in
\cite[Lemma~2.3]{AFHZI}.

The modulation of \cite[Lemma 2.3]{AFHZI} gives
the following decomposition for $u\in V^\eta_{\rho^*}(N,\epsilon_0)$
\[
u=U_z+w,\qquad w\in E_z,
\]
and hence
\[
\int_M V_{i,t}w_i\,dV_g=0
\qquad\text{for all }i,t.
\]

For bounded $L|b|$ and large $L$, the normalized tangent frame has a
uniformly positive definite Gram matrix and uniformly bounded first
derivatives \cite[Lemma~2.3 and Lemmas~A.2--A.4]{AFHZI}.

\subsection{Critical points at infinity}
Following Bahri's variational framework \cite{Bahri1989,BahriCoron1988,BahriLiRey1995}, we describe
the loss of compactness through noncompact trajectories of a descending
field.
\begin{definition}
\label{II:def:critical-point-infinity}
A positive orbit $s\mapsto u(s)$ is a \emph{chart-escaping orbit} if its
tail remains in the union of the full $N$-bubble charts and, in modulation
coordinates,
\begin{equation*}
L(s)\to+\infty,\qquad
J_{\rho^*}(u(s))\to c_\infty\in\mathbb R,\qquad
(p(s),r(s))\to(p^*,r^*).
\end{equation*}
The pair $(p^*,r^*)$ is called its \emph{limiting configuration}.

Two chart-escaping orbits are equivalent if their limiting configurations
agree after a simultaneous permutation of the centers and relative
log-scales. A \emph{critical point at infinity} is an equivalence class
of chart-escaping orbits.
\end{definition}

	\section{Normal reduction and finite-scale centering}\label{sec:3}
	

\subsection{Uniform estimates on the original normal fibers}
Let $Z_{0,i}=2+y\cdot\nabla v_i$ and $Z_{k,i}=\partial_{y_k}v_i$, for 
$k=1,2$. The inverse matrix linearization is
$$(\mathscr L_v\phi)_i=-\sum_j a^{ij}\Delta\phi_j-e^{v_i}\phi_i.$$

The orthogonality conditions give uniform weighted estimates on each
rescaled disk. We use these estimates to invert the Hessian restricted
to $E_z$.
\begin{lemma}
\label{II:lem:normal-model}
Fix the entire profile, a positive lower bound on the distances
between the centers, and upper bounds on the ratios of the scales
and on $L|b|$. Let
\[
B_1:=\{y\in\mathbb R^2:|y|<1\},
\qquad
A_j:=\{y\in\mathbb R^2:2^j<|y|<2^{j+1}\},\quad j\ge0.
\]
For $w\in E_z$ with $\|w\|_{\mathcal H}=1$, denote $w_i$ also for its
rescaling on the disk at $p_t$. Then
\begin{equation}\label{II:eq:m4-anchored-test-moments}
 \|w_i\|_{L^2(B_1)}\le C,\qquad
 \|w_i\|_{L^2(A_j)}\le C2^j\sqrt{1+j}.
\end{equation}
 For every $\nu>0$,
there are $c_\nu,C_\nu>0$ such that
\[
 \int_{\lambda_t y_{p_t}(U(p_t))}
 (1+|y|)^{-2-\nu}\exp(c_\nu w_i^2)\,dy\le C_\nu.
\]
The same moment bounds hold for any $w\in\mathcal H$ with
$\|w\|_{\mathcal H}=1$ and
\[
 \max_{i,t}\left|\int_M V_{i,t}w_i\,dV_g\right|\le C_0,
\]
with constants depending also on $C_0$.

Let $P_z$ be the Dirichlet orthogonal projection onto $E_z$. For
large $L$, the operator on $E_z$ associated with the form
\[
 \sum_{i,j}a^{ij}\int_M
 \langle\nabla_g w_i,\nabla_g w'_j\rangle_g\,dV_g
 -\sum_{i,t}\int_M V_{i,t}w_iw'_i\,dV_g
\]
has a uniformly bounded inverse. If $\Delta_\Omega$ is sufficiently
small, the same holds for $P_zD^2J_{\rho^*}(U_z)|_{E_z}$. On a fixed
small ball in $E_z$, the map
\[
 h\mapsto P_zD^2J_{\rho^*}(U_z+h)|_{E_z}
\]
is uniformly Lipschitz in operator norm.
\end{lemma}


\begin{proof}
On a rescaled disk, let $a=\frac 1{|B_1|}\int_{B_1} w_i$. Radial
Cauchy--Schwarz and Poincar\'e inequalities  give
\[
 \|w_i-a\|_{L^2(B_1)}\le C,\qquad
 \|w_i-a\|_{L^2(A_j)}\le C2^j\sqrt{1+j}.
\]
The same estimates apply to the truncated outer annulus. The cutoff is
one on the rescaled $B_1$, where $\int_{B_1}e^{v_i}\,dy>0$. Thus
\[
 |a|\le C\left(
  \left|\int_M V_{i,t}w_i\,dV_g\right|
    +1+\sum_{j\ge0}2^{(2-m_i)j}\sqrt{1+j}\right)\le C.
\]
This proves \eqref{II:eq:m4-anchored-test-moments}, including the case
of bounded integrals against $V_{i,t}$.

Apply Moser--Trudinger to the oscillation about the annular mean on a
fixed rescaled annulus \cite[Theorem~1.7]{Fontana}. The mean is
$O(\sqrt{1+j})$, so
\[
 \int_{A_j}(1+|y|)^{-2-\nu}e^{c_\nu w_i^2}\,dy
 \le C\,2^{-\nu j}e^{C c_\nu(1+j)}.
\]
Choose $c_\nu$ so that $Cc_\nu<\nu\log 2$. Use Moser--Trudinger inequality  on $B_1$
and on $M$ for the inner disk and the fixed exterior, respectively.
For each fixed $k\ge0$, the annular bounds also give
\[
 \int_{A_j}e^{v_l}(1+\log(2+|y|))^k |w_i|^2\,dy
 \le C\,2^{(2-m_l)j}(1+j)^{k+1}.
\]

Let $d\omega$ be the spherical area and let $d$ denote the distance to the
point at infinity. Then
\[
 k_i:=e^{v_i}\frac{dy}{d\omega}=O(d^{m_i-4}),
 \qquad k_i\in L^s(S^2)\quad\text{for some }s>1.
\]
A rescaled limit $f$ has finite Dirichlet energy and
$\int_{A_j}|f_i|^2\,d\omega\le C2^{-2j}(1+j)$. Thus $f$ belongs to
$H^1(S^2)^n$ by conformal invariance and the zero $2$-capacity of the
omitted point, applied first to bounded truncations.

Let $\mathcal B_\infty$ be the plane form of $\mathscr L_v$ and let
$\langle\cdot,\cdot\rangle_D$ denote the Dirichlet pairing on the plane.
Suppose $f$ satisfies
\[
 \int_{\mathbb R^2}e^{v_i}f_i\,dy=0,\qquad
 \langle f,Z_k\rangle_D=0\quad(0\le k\le2)
\]
and $\mathcal B_\infty(f,\cdot)$ vanishes on this constrained space.
For arbitrary tests $\phi$,
\[
 \mathcal B_\infty(f,\phi)
 =\sum_i\beta_i\int_{\mathbb R^2}e^{v_i}\phi_i\,dy
   +\sum_{k=0}^2\mu_k\langle Z_k,\phi\rangle_D.
\]
The tests $\phi=e_i$ give $\beta_i=0$. Since $\mathscr L_v Z_l=0$,
testing with $\phi=Z_l$ gives
\[
 \sum_{k=0}^2\mu_k\langle Z_k,Z_l\rangle_D=0
 \quad(0\le l\le2),\qquad \mu=0.
\]
The equation on the sphere has right-hand side $A(k_i f_i)_i\in L^p$
for some $p>1$, by H\"older and Sobolev embedding. Elliptic regularity
gives boundedness. By \cite[Theorem~2.1(2)]{LZone},
$f\in\operatorname{span}\{Z_0,Z_1,Z_2\}$, and the Dirichlet
orthogonality gives $f=0$.

Let $T_z$ be the Dirichlet representative of the model form. If the
uniform lower bound fails, then along a sequence with $L\to\infty$,
\[
 w\in E_z,\qquad \|w\|=1,\qquad
 P_zT_zw=o_{\mathcal H}(1),\qquad
 T_zw=o_{\mathcal H}(1)+\sum_a\mu_a T_a ,
\]
where $T_a$ are the normalized tangent vectors. The pairings with the
position and scale directions are bounded by the weighted moments.
For an amplitude direction, the projected Poisson equation and
\cite[Lemma~A.1]{AFHZI} give
\[
 \left|\left\langle T_zw,\frac{e_iPW_{i,t}}{\sqrt{\ell_t}}
          \right\rangle_{\mathcal H}\right|\le C L^{-1/2}.
\]
Here, in each potential pairing, we subtract the constant of $PW_{i,t}$
on the corresponding core using $\int_M V_{i,s}w_i\,dV_g=0$. The
remainder is $O(1+\log(2+|y|))$ on its own disk and bounded on the
others. The inverse Gram bound gives $|\mu|\le C$. Since the mixed Gram
block tends to zero, the amplitude multipliers tend to zero
\cite[Lemma~A.4]{AFHZI}.

For the dilation direction, use $Z_{0, i}+q_i=m_i+y\cdot\nabla v_i$.
Its value and gradient have orders $O(|y|^{-(\widehat m-2)})$ and
$O(|y|^{-1-(\widehat m-2)})$. The corresponding orders for the
translations are $O(|y|^{-1})$ and $O(|y|^{-2})$. Cutting off these
representatives at diverging radii therefore costs $o(1)$ in Dirichlet
norm. Local Rellich compactness, the weighted tails, and the tangent
expansions in \cite[Lemmas~A.1--A.4]{AFHZI} give, on each core, a limit
$f$ with
\[
 \int_{\mathbb R^2}e^{v_i}f_i\,dy=0,\qquad
 \langle f, Z_k\rangle_D=0.
\]
For a compactly supported rescaled test, subtract its surface mean
before testing the equation for $T_zw$. Its Dirichlet pairing with each
normalized amplitude is $O(L^{-1/2})$ by the projected Poisson
equation. Passing to the limit gives
\[
 \mathcal B_\infty(f,\phi)=
   \sum_{k=0}^2\widetilde\mu_k\langle Z_k,\phi\rangle_D .
\]
The identity extends to $H^1(S^2)^n$ by density across the omitted point
and $k_i\in L^s(S^2)$. The plane argument gives $f=0$. Testing the
identity with $Z_k$ gives $\widetilde\mu=0$, and hence the position and
scale multipliers on that core tend to zero. Local compactness and
the weighted tails give $\int_M V_{j,t}w_i^2\,dV_g=o(1)$ for all
$i,j,t$. Since $KA=I$,
\[
 1=\|w\|^2=
 \langle T_zw,Aw\rangle_{\mathcal H}
 +\sum_{i,t}\int_M V_{i,t}w_i(Aw)_i\,dV_g=o(1).
\]
The compression $P_zT_z|_{E_z}$ is self-adjoint. Its lower bound gives
closed range and trivial kernel, hence surjectivity.

The projected expansion \cite[Lemma~A.1]{AFHZI}, with the masses
in \eqref{II:def:full-local-mass}, gives the radial factor
$e^{v_i+b_{i,t}\psi_i}$, where $\psi_i=\sum_j a^{ij}v_j$. After division
by their values at the center, the smooth factors differ from one by
$O(e^{-\eta L}(1+|y|)^\eta)$ for some
$0<\eta<\min(1,\widehat m-2)$. The difference between the normalized
radial densities at $b$ and at $0$ is bounded by
\[
 C|b|(1+\log(2+|y|))(1+|y|)^{-2-\nu}
\]
for some $\nu>0$. On normal tests, the coefficient of the leading
density on disk $t$ is $N\Omega_{i,t}$. The weighted moments bound the
errors in the linear and quadratic pairings by
$C(|b|+e^{-\eta L})$, including the fixed exterior. The anchor gives the
same bound for the means in the covariance. Consequently
\[
 \bigl\|P_z(D^2J_{\rho^*}(U_z)-T_z)|_{E_z}\bigr\|
 \le C(\Delta_\Omega+|b|+e^{-\eta L}).
\]
For small $\Delta_\Omega$ and large $L$, the inverse follows by a
Neumann series.

For a tilt $h\in E_z$ in a fixed small ball, Jensen's inequality gives
$\int_M e^{h_i}d\mu_i(U_z)\ge e^{-C\|h\|}$. H\"older and the Gaussian
bound give uniform third absolute moments of unit normal tests for
\[
 d\mu_i(U_z+h)=\frac{e^{h_i}\,d\mu_i(U_z)}
                    {\int_M e^{h_i}\,d\mu_i(U_z)}.
\]
The third derivative of $\log H_i$ is the centered triple moment under
this measure. Thus $D^3J_{\rho^*}(U_z+h)$ is uniformly bounded on three
normal unit tests, which proves the Lipschitz estimate.
\end{proof}

We will also use the ambient moment estimate
\begin{equation}\label{II:eq:ambient-moments}
 \int_M |w_i|^2\,d\mu_i(u)\le C(1+\|u\|^2)\|w\|^2 .
\end{equation}
Assume $\|w\|=1$. Since $u_i$ has mean zero, Jensen's inequality gives
$
H_i(u)=\int_M h_i e^{u_i}\,dV_g
\ge \min_M h_i.
$
By H\"older's and Moser--Trudinger inequalies,
\[
\begin{aligned}
\mu_i(u)(\{|w_i|>s\})
&\le C\|e^{u_i}\|_{L^2(M)}
       \operatorname{vol}_g(\{|w_i|>s\})^{1/2}\le \min\{1,Ce^{C\|u\|^2-cs^2}\}.
\end{aligned}
\]
Integrating this bound and splitting at
$s=C_1(1+\|u\|^2)^{1/2}$ for a sufficiently large $C_1$, we obtain \eqref{II:eq:ambient-moments}. 

\subsection{The corrected family at the borderline exponent}

Assume $A\sigma=4\mathbf 1$ throughout this subsection and the next
subsection on normal reduction. Set
\[
a=(2\max_i\sigma_i)^{-1/2},
\qquad
u_a(y)=\log\frac{8a^2}{(1+a^2|y|^2)^2}.
\]
Let $P_t$ be the mean-zero solution of
\[
-\Delta_gP_t
=\chi_t e^{-\varphi_{p_t}}\lambda_t^2e^{u_a(\lambda_t y_{p_t})}
-\overline{\chi_t e^{-\varphi_{p_t}}\lambda_t^2e^{u_a(\lambda_t y_{p_t})}}.
\]
Radial uniqueness \cite[Theorem~1.1]{LinZhang2010} and linearity give
\begin{equation}\label{II:eq:m4-explicit-profile}
\begin{aligned}
v_i&=u_a+\log(\sigma_i/4),&
e^{\fJ_i}&=2\sigma_i/a^2,\\
PW_{i,t}&=\tfrac14\sigma_iP_t,&
U_{z,i}&=\sum_t s_{i,t}P_t,\qquad
s_{i,t}=1+\tfrac14\sigma_i b_{i,t}.
\end{aligned}
\end{equation}

Set $d_t=(a\lambda_t)^{-1}$. The scalar projection satisfies
\cite[Lemma~A.1]{AFHZI}
\begin{equation}\label{II:eq:m4-mean-free-scalar-model}
\begin{aligned}
P_t&=\mathcal T_t+\gamma_t+\mathfrak e_t,
&\gamma_t&=-\overline{\mathcal T_t},\\
\mathcal T_t&=-2\chi_t\log(|y_{p_t}|^2+d_t^2)
              +8\pi H^g(\cdot,p_t),
&\overline{\mathfrak e_t}&=0.
\end{aligned}
\end{equation}

For every fixed $k$,
\begin{equation}
    \label{eq:assertion1} \|\mathfrak e_t\|_{C^k(M)}\le C_k d_t^2,
\end{equation}
with the same bound after the derivatives in $p_t$ and $\ell_t$ used
below. {\it Indeed}, away from $p_t$,
\[
\mathcal T_t
=8\pi G^g(\cdot,p_t)
 -2\chi_t\log\!\left(1+\frac{d_t^2}{|y_{p_t}|^2}\right).
\]
The Dirac terms cancel, and
\[
-\Delta_g\mathcal T_t
=\chi_t e^{-\varphi_{p_t}}\lambda_t^2e^{u_a(\lambda_t y_{p_t})}
 -8\pi+R_t,
\qquad
\|R_t\|_{C^k(M)}\le C_kd_t^2.
\]
Here $R_t$ is supported where the cutoff varies. The same estimates hold
after differentiation in $p_t$ and $\ell_t$. Subtracting the mean gives
$-\Delta_g\mathfrak e_t=-R_t+\overline{R_t}$, then  the estimate \eqref{eq:assertion1} follows
from the mean-zero Poisson estimate.

\subsection{Projection of the first angular mode}

Assume $A\sigma=4\mathbf 1$. In centered orthonormal coordinates, set
\[
g_{i,t}:=\nabla_x\log\rf_{p,t}^i(p_t).
\]
We adapt the scalar correction in
\cite[Section~2 and Appendix~A]{AhmedouBenAyed2017} to the system.

\begin{lemma}\label{II:prop:m4-first-mode-sphere}
At a fixed center, denote $g_i=g_{i,t}$. There is a unique first angular
mode $\Phi$ whose stereographic transform belongs to $W^{2,p}(S^2)$
for $1<p<2$ and satisfies
\[
\begin{aligned}
(\mathscr L_v\Phi)_i
  &=e^{v_i}g_i\cdot y-\sum_{k=1}^2c_k(-\Delta Z_{k,i}),\qquad
\sum_i\int_{\mathbb R^2}e^{v_i}\Phi_iZ_{k,i}\,dy=0
\quad(k=1,2), 
\end{aligned}
\]
where $
c_k=-\frac{3}{16na^2}\sum_i\sigma_i(g_i)_k$.
For $0\le j\le3$,
\[
|D^j\Phi(y)|\le C|g|(1+|y|)^{-1-j}\log(2+|y|).
\]
\end{lemma}

\begin{proof}
Integration by parts gives
\[
\int_{\mathbb R^2}e^{v_i}(g_i\cdot y)Z_{k,i}\,dy
=-2\pi\sigma_i(g_i)_k.
\]
Under stereographic projection, $e^{u_a}\,dy=2\,d\omega$ and
$Z_{k,i}=-2aX_k$, where $X_k$ is the $k$th coordinate on $S^2$.
Hence the Dirichlet pairings of the translation modes are
\[
\sum_i\int_{\mathbb R^2}\nabla Z_{k,i}\cdot\nabla Z_{l,i}\,dy
=4na^2\int_{S^2}\nabla_{S^2}X_k\cdot\nabla_{S^2}X_l\,d\omega
=\frac{32\pi na^2}{3}\delta_{kl}.
\]
These two identities give the coefficients $c_k$.

In the variables $D_\sigma^{1/2}\Phi$, diagonalizing
$D_\sigma^{1/2}AD_\sigma^{1/2}/4$ reduces the equations on $S^2$ to
\[
(-\Delta_{S^2}-2\mu_j)\eta_j=f_j.
\]
By Perron--Frobenius, $1$ is simple and $0<|\mu_j|\le1$. The kernel
therefore consists of the spherical harmonics of degree one in the
eigenspace for $\mu_j=1$. The coefficients $c_k$ remove the two
translation pairings. The pairing with the radial mode vanishes by
angular integration. The forcing is $O(d^{-1})$, where $d$ is the
distance to the omitted pole. Thus $f_j\in L^p(S^2)$ for $1<p<2$.
The Fredholm alternative gives the normalized solution in
$W^{2,p}(S^2)$ and hence in $L^\infty(S^2)$. Rotational invariance
preserves the first angular mode.

Each radial coefficient $G$ of $\Phi$ satisfies
\[
G''+\frac1rG'-\frac1{r^2}G=h,\qquad
h(r)=O(|g|r^{-3})\quad(r\ge1).
\]
Regularity at the origin and boundedness at infinity give
\[
G(r)=-\frac r2\int_r^\infty h(s)\,ds
     -\frac1{2r}\int_0^r s^2h(s)\,ds .
\]
It follows that $|G(r)|+r|G'(r)|\le C|g|r^{-1}\log(2+r)$ for $r\ge1$.
Interior estimates on rescaled annuli give the bounds through order
three. Interior regularity gives the bounds on $|y|\le1$.
\end{proof}

Let $\Phi^t$ and $c_{t,k}$ be the solution and coefficients at $p_t$.
With $Y_t=\lambda_t y_{p_t}$, set
\begin{align}
\widetilde C_i&=\sum_t\varepsilon_t\chi_t\Phi_i^t(Y_t),&
C_i&=\widetilde C_i-\overline{\widetilde C_i},&
\widehat U&=U+C,
\label{II:eq:first-mode-surface}\\
F_{i,t}&=\varepsilon_tV_{i,t}\bigl(\Phi_i^t(Y_t)+g_{i,t}\cdot Y_t\bigr),&
-\Delta_gQ_{i,t}&=F_{i,t}-\overline{F_{i,t}},&
\overline{Q_{i,t}}&=0,\nonumber\\
B_{z,i}&=U^\sharp_{z,i}
       =\widehat U_{z,i}+\sum_t b_{i,t}Q_{i,t},&
E_z^\sharp&=(\operatorname{range}DB_z)^\perp .
\nonumber
\end{align}
Since $V_{i,t}\,dV_g$ is radial in $Y_t$, angular integration gives
$\overline{F_{i,t}}=0$. For $w\in E_z^\sharp$, the amplitude direction
$\partial_{b_{i,t}}B_z=e_i(PW_{i,t}+Q_{i,t})$ gives
\begin{equation}\label{II:eq:amplitude-corrected-anchor}
0=\int_M\langle\nabla_g(PW_{i,t}+Q_{i,t}),\nabla_gw_i\rangle_g\,dV_g
 =\int_M(V_{i,t}+F_{i,t})w_i\,dV_g.
\end{equation}

\begin{lemma}
\label{II:lem:amplitude-projected-size}
The $H^1$ and $L^\infty$ norms of $C=(C_1,\ldots,C_n)$  and $Q_{i,t}$ are
$O(\varepsilon)$, including their normalized derivatives in position
and logarithmic scale through order two.
\end{lemma}

\begin{proof}
Denote $C^t$ for the mean-zero part of
$\varepsilon_t\chi_t\Phi^t(Y_t)$ and $Z_k^t$ for the mean-zero part of
$\chi_t Z_k(Y_t)$. By conformal invariance and the bounds in
Lemma~\ref{II:prop:m4-first-mode-sphere},
\[
\|C^t\|_{H^1\cap L^\infty}\le C\varepsilon_t,
\qquad
\|Z_k^t\|_{H^1\cap L^\infty}\le C.
\]
These bounds also hold after the indicated derivatives. On the
profile, a normalized position derivative gives a translation
derivative, and a logarithmic scale derivative gives a radial
derivative $Y_t\cdot\nabla$.

The equation for $\Phi^t$ gives
\[
-\Delta_g\left(\sum_j a^{ij}C_j^t
         +\varepsilon_t\sum_{k=1}^2c_{t,k}Z_{k,i}^t\right)
=F_{i,t}+R_{i,t}.
\]
The remainder is supported where $\chi_t$ varies. On that annulus
$|Y_t|\asymp\lambda_t$, and the decay of $\Phi^t$ and $Z_k$ gives
\[
\|R_{i,t}\|_{L^p(M)}\le C_p\varepsilon^2L
\quad(p>2),
\]
with the same bound after the indicated derivatives. Both sides of
the equation have mean zero, so $\overline{R_{i,t}}=0$. The mean-zero
Poisson estimate and $W^{2,p}(M)\hookrightarrow H^1(M)\cap L^\infty(M)$
now give
\[
Q_{i,t}
=\sum_j a^{ij}C_j^t
 +\varepsilon_t\sum_{k=1}^2c_{t,k}Z_{k,i}^t
 +O_{H^1\cap L^\infty}(\varepsilon^2L).
\]
The same remainder bound holds after the indicated derivatives.
Since $\varepsilon_t\asymp\varepsilon$, this proves the estimate.
\end{proof}

\subsection{Surface masses and the normal residual}

Work on a chart with separated centers, comparable scales, and
$L|b|\le C$. For $\mathcal U=U$ or $\widehat U$, use $\widetilde C=0$
or the correction in \eqref{II:eq:first-mode-surface}, respectively.
Set
\begin{align*}
Z_i
&=\int_M h_i\exp\left[\sum_l s_{i,l}(\mathcal T_l+\mathfrak e_l)
                         +\widetilde C_i\right]dV_g,\\
\mathfrak h_{i,t}(x)
&=e^{\varphi_{p_t}(x)}h_i(x)
 \exp\left(8\pi s_{i,t}H^g(x,p_t)
           +8\pi\sum_{l\ne t}s_{i,l}G^g(x,p_l)\right),\\
J_s&=\frac{\pi}{a^2(2s-1)},\qquad
k_s(Y)=J_s^{-1}(1+a^2|Y|^2)^{-2s}.
\end{align*}
By \eqref{II:eq:m4-explicit-profile}, polar integration gives
\[
\sigma_{i,t}(b)=\frac{\sigma_i}{2s_{i,t}-1}
 \exp\!\left(b_{i,t}\sum_j a^{ji}\log(2\sigma_j a^2)\right).
\]
Using $\sum_j a^{ji}=\sigma_i/4$ in
\eqref{II:def:full-local-mass}, we obtain
\[
\mathfrak M_{i,t}
=a^{4s_{i,t}}\lambda_t^{4s_{i,t}-2}
  \mathfrak h_{i,t}(p_t)J_{s_{i,t}}.
\]

We first expand the denominator and the normalized density near each
center.

\begin{proposition}\label{II:prop:m4-global-denominator}
The following estimates hold: 
\begin{equation}\label{II:eq:m4-surface-denominator-exact}
 H_i(\mathcal U)
 =e^{\sum_t s_{i,t}\gamma_t-\overline{\widetilde C_i}}Z_i,
 \qquad
 |D_p^j\partial_\ell^k\partial_b^l\eta_i|
 \le C\varepsilon^2 L^{1+l},
 \quad j+k+l\le2.
\end{equation}
On the rescaled inner disk at $p_t$, set $s=s_{i,t}$ and
$d_{i,t}=\nabla_x\log\mathfrak h_{i,t}(p_t)$. Then
\begin{equation}\label{II:eq:m4-normalized-density-error}
 \frac{d\mu_i(\mathcal U)}{dY}
 =\Omega_{i,t}k_s(Y)
 \left[1+\varepsilon_t\bigl(d_{i,t}\cdot Y+\Phi_i^t(Y)\bigr)
       +O\bigl(\varepsilon^2(L+|Y|^2)\bigr)\right],
\end{equation}
where $\Phi_i^t=0$ for $\mathcal U=U$. Outside the inner disks,
$d\mu_i(\mathcal U)/dV_g=O(\varepsilon^2)$.
\end{proposition}

\begin{proof}
The formula for $H_i(\mathcal U)$ follows from
\eqref{II:eq:m4-mean-free-scalar-model} and the definition of
$\widetilde C_i$. For $|Y|<r_0/(2\varepsilon_t)$, with
$x=y_{p_t}^{-1}(\varepsilon_tY)$, the integrand for $Z_i$ in the
rescaled coordinates is
\[
\mathfrak M_{i,t}k_s(Y)
\frac{\mathfrak h_{i,t}(x)}{\mathfrak h_{i,t}(p_t)}
\exp\left(\sum_l s_{i,l}\mathfrak e_l(x)
               +\varepsilon_t\Phi_i^t(Y)\right).
\]
Taylor expansion of the last two factors gives
\[
1+\varepsilon_t\bigl(d_{i,t}\cdot Y+\Phi_i^t(Y)\bigr)
  +O\bigl(\varepsilon^2(1+|Y|^2)\bigr).
\]
Both terms of order $\varepsilon_t$ have zero angular integral.
Since $|s-1|L\le C$, for $R_t=r_0/(2\varepsilon_t)$ we have
\[
 k_s\le Ck_1\quad\text{on }B_{R_t},\qquad
 \int_{|Y|>R_t}k_s\,dY=O(\varepsilon^2),\qquad
 \int_{B_{R_t}} |Y|^2 k_s\,dY=O(L).
\]
The exterior contributes $O(1)$ to $Z_i$, while
$\mathfrak S_i\asymp\varepsilon^{-2}$. Hence
\[
Z_i=\mathfrak S_i\bigl(1+O(\varepsilon^2L)\bigr).
\]

After recentering, up to two derivatives in $p$ preserve the Taylor
remainder since $h_i\in C^4$. Derivatives in $\ell$ preserve the tail
bound, including the contribution from the moving boundary.
A derivative in $b$ gives either a factor $O(L)$ from the prefactor
or a logarithm from the kernel. The latter is controlled by
\[
\int_{B_{R_t}} |Y|^2(1+\log(2+|Y|))^l k_s(Y)\,dY
\le C L^{1+l}.
\]
Differentiating $(Z_i-\mathfrak S_i)/\mathfrak S_i$ gives the bound
for $\eta_i$. Dividing the local integrand by $Z_i$ gives
\eqref{II:eq:m4-normalized-density-error}.
\end{proof}

The density expansion and the orthogonality conditions on
$E_z^\sharp$ give the following estimate for the normal residual.

\begin{proposition}\label{II:prop:m4-amplitude-corrected-sharp-residual}
Let $P_z^\sharp$ be the Dirichlet orthogonal projection onto
$E_z^\sharp$. Uniformly for $L|b|\le C$,
\[
\|P_z^\sharp\nabla J_{\rho^*}(B_z)\|\le C R_L,
\qquad R_L=|b|+\varepsilon^2L^{3/2}.
\]
For $w\in E_z^\sharp$,
\begin{equation}\label{II:eq:m4-sharp-normal-density-mean}
\left|\int_M w_i\,d\mu_i(B_z)\right|\le C R_L\|w\|.
\end{equation}
\end{proposition}

\begin{proof}
Take $w\in E_z^\sharp$ with $\|w\|=1$. For a first Fourier coefficient
$f$ of the rescaled $w_i$,
\[
 |f(r)|^2
 =2\int_0^r f(s)f'(s)\,ds
 \le\int_0^r\left(s|f'(s)|^2+\frac{|f(s)|^2}{s}\right)ds
 \le C.
\]
The identity follows first for smooth tests and then by density.
Since $F_{i,t}$ has first angular frequency and
\[
\int_{\mathbb R^2}e^{v_i(Y)}
 \bigl(|\Phi_i^t(Y)|+|g_{i,t}\cdot Y|\bigr)\,dY<\infty,
\]
the corrected orthogonality condition
\eqref{II:eq:amplitude-corrected-anchor} gives
\[
\left|\int_M V_{i,t}w_i\,dV_g\right|
=\left|\int_M F_{i,t}w_i\,dV_g\right|\le C\varepsilon_t.
\]
Thus Lemma~\ref{II:lem:normal-model} applies. In particular,
\[
\int_{B_{R_t}} k_1(Y)|Y|^2|w_i|\,dY
\le C\sum_{j\le C L}\sqrt{1+j}\le C L^{3/2},
\]
and the first and logarithmic moments are bounded.

On the inner disk we have
\[
d_{i,t}-g_{i,t}=O(|b|),\qquad
|k_s-k_1|\le C|b|(1+\log(2+|Y|))k_1.
\]
Combining these bounds with \eqref{II:eq:m4-normalized-density-error}
and the weighted moments yields
\[
\left|\int_M
 \left[D_i(\widehat U)
 -\sum_t N\Omega_{i,t}(V_{i,t}+F_{i,t})\right]w_i\,dV_g\right|
\le C R_L .
\]
By Lemma~\ref{II:lem:amplitude-projected-size},
$\|B_z-\widehat U\|_{L^\infty}=O(|b|\varepsilon)$, so the normalized
density changes by a relative error of the same order. Each integral
against $V_{i,t}+F_{i,t}$ vanishes by
\eqref{II:eq:amplitude-corrected-anchor}. This proves
\eqref{II:eq:m4-sharp-normal-density-mean}.

Modulo constants, the sources in the quadratic part of the first
variation are
\[
\begin{aligned}
-\sum_j a^{ij}\Delta_g B_{z,j}
={}&\sum_t(V_{i,t}+F_{i,t})
   +\sum_{j,t}a^{ij}b_{j,t}(V_{j,t}+F_{j,t})\\
 &-\sum_{t,k}c_{t,k}\varepsilon_t(-\Delta_g Z^t_{k,i})
   +O_{L^2(M)}(\varepsilon^2 L).
\end{aligned}
\]
The first sum vanishes against $w_i$. The terms containing $b$ give
$O(|b|)$ by the weighted moments. In centered coordinates,
$\partial_{(p_t)_k}y_{p_t}=-e_k+O(|y_{p_t}|)$. The gradient error on
the inner disk is bounded by $C\varepsilon/(|y_{p_t}|+\varepsilon)$,
whose squared integral is $O(\varepsilon^2L)$. The smooth and cutoff
terms give $O(\varepsilon)$, the factors $s_{i,t}-1$ give $O(|b|)$,
and the projected corrections satisfy
Lemma~\ref{II:lem:amplitude-projected-size}. Hence
\[
\left\|\lambda_t^{-1}\partial_{(p_t)_k}B_z+Z_k^t\right\|
\le C\bigl(|b|+\varepsilon\sqrt L\bigr).
\]
Since $w$ is orthogonal to $\partial_{(p_t)_k}B_z$, the translation
terms contribute at most
$C\varepsilon(|b|+\varepsilon\sqrt L)$. Together with
\eqref{II:eq:m4-sharp-normal-density-mean}, these bounds give
\[
|DJ_{\rho^*}(B_z)[w]|
\le C\left[R_L+|b|+\varepsilon^2L
              +\varepsilon(|b|+\varepsilon\sqrt L)\right]
\le C R_L .
\]
Taking the supremum over unit vectors in $E_z^\sharp$ proves the
residual estimate.
\end{proof}

	\subsection{The normal graph and tangent pairings}
	
	The preceding estimates give a uniformly invertible normal Hessian and hence a
	finite-dimensional normal graph.
	
	\begin{theorem}\label{II:thm:m4-verified-normal-inverse}
		Fix a compact separated cell with comparable scales and bounded $L|b|$.
		There exist $\delta_*,L_*,r_*,C>0$ such that, for
		$
		L\ge L_*,
		\Delta_\Omega\le\delta_*,
		$
		the compression of $D^2J_{\rho^*}(B_z)$ to $E_z^\sharp$ is invertible and
		\[
		\left\|
		\left(
		P_z^\sharp D^2J_{\rho^*}(B_z)\big|_{E_z^\sharp}
		\right)^{-1}
		\right\|
		\le C.
		\]
		Moreover, for $v_1,v_2\in E_z^\sharp$ with $\|v_j\|\le r_*$,
		\[
		\left\|
		P_z^\sharp
		\left(
		D^2J_{\rho^*}(B_z+v_1)-D^2J_{\rho^*}(B_z+v_2)
		\right)\big|_{E_z^\sharp}
		\right\|
		\le C\|v_1-v_2\|.
		\]
		The equation
		$P_z^\sharp\nabla J_{\rho^*}(B_z+\bar v)=0$
		has a unique solution in this ball, and
		\begin{equation}\label{II:eq:m4-verified-normal-solution}
			\|\bar v\|\le CR_L,
			\qquad
			G_z:=B_z+\bar v(z).
		\end{equation}
		The same conclusions hold for $J_\rho$ and a normal solution $\bar v_\rho$
		when $|\rho-\rho^*|$ is sufficiently small. The normal graph is $C^2$ in
		finite-scale local coordinates.
	\end{theorem}
	
	\begin{proof}
		The normalized tangent frames of $B_z$ and $U_z$ differ by
		$O(\varepsilon)$; see \cite[Lemmas~A.2--A.4]{AFHZI}. Hence
		\[
		\|P_z^\sharp-P_z\|=O(\varepsilon),
		\qquad
		P_z^\sharp\big|_{E_z}:E_z\to E_z^\sharp
		\]
		is a uniformly bounded isomorphism.

Pulling the normal Hessian at $B_z$ back to $E_z$ gives
		\[
		\left\|
		\left(P_z^\sharp\big|_{E_z}\right)^*
		D^2J_{\rho^*}(B_z)
		\left(P_z^\sharp\big|_{E_z}\right)
		-
		D^2J_{\rho^*}(U_z)\big|_{E_z}
		\right\|
		\le C\varepsilon L=o(1),
		\]
        where $D^2J_{\rho^*}(U_z)|_{E_z}$ denotes the restriction
of the Hessian as a bilinear form to $E_z\times E_z$. 
		Indeed,
		\[
		\|D^2J_{\rho^*}(U_z)\|=O(L),
		\qquad
		\|B_z-U_z\|_{L^\infty}=O(\varepsilon),
		\]
		and the corresponding probability densities differ by a relative
		$O(\varepsilon)$ term. Lemma~\ref{II:lem:normal-model} gives a uniformly bounded inverse
for the limiting operator. Since the perturbation is $o(1)$ in
operator norm, the Neumann series gives the asserted inverse bound.

		The corrected anchors and Lemma~\ref{II:lem:normal-model} also give the
		uniform Lipschitz estimate for the normal Hessian. Set
		$
		F(v):=P_z^\sharp\nabla J_{\rho^*}(B_z+v),
		H:=DF(0).
		$
		Since the normal Hessian is uniformly Lipschitz, there exists
		$C_{\mathrm{Lip}}>0$ such that
		\[
		\|DF(v_1)-DF(v_2)\|
		\le C_{\mathrm{Lip}}\|v_1-v_2\|,
		\qquad
		\|v_1\|,\|v_2\|\le r_*.
		\]
		Choose $r_*>0$ so that
		\[
		\|H^{-1}\|C_{\mathrm{Lip}}r_*\le\frac12.
		\]
		By Proposition~\ref{II:prop:m4-amplitude-corrected-sharp-residual},
		\[
		\|H^{-1}F(0)\|\le CR_L\le\frac{r_*}{2}
		\]
		for large $L$. Therefore
		\[
		v\mapsto v-H^{-1}F(v)
		\]
		is a contraction on $\overline{B}_{r_*}(0)$, and its fixed point satisfies
		$
		\|\bar v\|\le 2\|H^{-1}F(0)\|\le CR_L.
		$
		
		Replacing $\rho^*$ by $\rho$ changes the normal residual and Hessian by
		\[
		O(|\rho-\rho^*|R_L),
		\qquad
		O(|\rho-\rho^*|),
		\]
		respectively; see \eqref{II:eq:m4-sharp-normal-density-mean}. The same
		contraction applies. Finally, $B_z$ is $C^3$ at each finite scale, while
		the normal projection and its local fiber trivialization are $C^2$.
		The implicit function theorem gives the asserted $C^2$ dependence.
	\end{proof}
	
	We next compute the derivatives of the reduced functional along the
	amplitude, scale, and position directions.
	
	\begin{theorem}\label{II:thm:m4-corrected-leading-pairings}
		Under the hypotheses of Theorem~\ref{II:thm:m4-verified-normal-inverse}, set
		\[
		A_{i,t}:=\ell_t^{-1}\partial_{b_{i,t}}B_z,
		\qquad
		S_t:=\partial_{\ell_t}B_z,
		\qquad
		P_{t,k}:=\lambda_t^{-1}\partial_{(p_t)_k}B_z,
		\qquad
		\delta_z:=|b|+\Delta_\Omega.
		\]
		Then
		\begin{align}
			\left|
			\bigl(DJ_{\rho^*}(G_z)-DJ_{\rho^*}(B_z)\bigr)[A_{i,t}]
			\right|
			&\le
			C\left(\frac{\|\bar v\|}{L}+\|\bar v\|^2\right),
			\label{II:eq:m4-verified-amplitude-normal-pairing}\\
			\left|
			\bigl(DJ_{\rho^*}(G_z)-DJ_{\rho^*}(B_z)\bigr)[T]
			\right|
			&\le
			C\left[
			(\delta_z+\varepsilon\sqrt L)\|\bar v\|+\|\bar v\|^2
			\right],
			\quad T=S_t,P_{t,k}.
			\label{II:eq:m4-verified-geometric-normal-pairing}
		\end{align}
		Moreover,
		\begin{align*}
			DJ_{\rho^*}(G_z)[A_{i,t}]
			&=
			2\pi\sigma_i
			\left((K\xi_t)_i+\sigma_i\theta_{i,t}\right)
			+O\left(
			\frac{\delta_z+\varepsilon}{L}+\|\bar v\|^2
			\right),\\
			DJ_{\rho^*}(G_z)[S_t]
			&=
			4\pi\sum_i\sigma_i\theta_{i,t}
			+O\left(
			\delta_z^2+\varepsilon^2L+\|\bar v\|^2
			\right),\\
			DJ_{\rho^*}(G_z)[P_{t,k}]
			&=
			-\frac{1}{N\lambda_t}\partial_{(p_t)_k}\ff(p)
			+O\left(
			\delta_z^2+\varepsilon^2L+\|\bar v\|^2
			\right).
		\end{align*}
	\end{theorem}
	
	\begin{proof}
		Let $w\in E_z^\sharp$, $\|w\|=1$. For the amplitude direction,
		\[
		\ell_t^{-1}
		\sum_j a^{ij}
		\bigl(\nabla w_j,\nabla(PW_{i,t}+Q_{i,t})\bigr)_{L^2}
		=O(L^{-1}).
		\]
		Indeed, the $Q_{i,t}$ term is $O(\varepsilon/L)$. For the remaining term,
		the corrected anchors and
		$
		V_{i,t}=\frac{\sigma_i}{\sigma_j}V_{j,t}
		$
		give the same bound. In the covariance term, subtract $\sigma_i$ on the
		$t$-th core and zero on the other cores. The remaining profile has weighted
		$L^2$ norm $O(L^{-1})$, while the constant part is controlled by
		Proposition~\ref{II:prop:m4-amplitude-corrected-sharp-residual}. Thus
		\[
		\left|D^2J_{\rho^*}(B_z)[w,A_{i,t}]\right|
		\le \frac{C}{L}.
		\]
		
		For $T=S_t$ or $P_{t,k}$, the leading core profile is
		$Z_0+2\mathbf1$ or $-Z_k$. The scalar model and the position estimate give
		\[
		\left|D^2J_{\rho^*}(B_z)[w,T]\right|
		\le C(\delta_z+\varepsilon\sqrt L).
		\]
		Here the kernel and weight errors are controlled by
		\eqref{II:eq:m4-normalized-density-error}, the covariance means by
		\eqref{II:eq:m4-sharp-normal-density-mean}, and the Jacobi terms cancel by
		the corrected anchors. Taylor expansion at $B_z$ yields
		\[
		\begin{aligned}
			\bigl(DJ_{\rho^*}(G_z)-DJ_{\rho^*}(B_z)\bigr)[T]
			&=
			D^2J_{\rho^*}(B_z)[\bar v,T]
			+O(\|\bar v\|^2),
		\end{aligned}
		\]
		which proves
		\eqref{II:eq:m4-verified-amplitude-normal-pairing} and
		\eqref{II:eq:m4-verified-geometric-normal-pairing}.
		
		For the base pairings, let $s_t=(s_{i,t})_i$ and define
		\begin{align}
			E^0_{tt}
			&:=
			32\pi\ell_t+32\pi\log a-16\pi+64\pi^2R^g(p_t),
			\qquad
			E^0_{ts}
			:=
			64\pi^2G^g(p_t,p_s),
			\qquad t\ne s,\notag
			\\
			\mathcal J_0
			&:=
			\frac12\sum_{t,s}s_t^TKs_sE^0_{ts}
			-\sum_i\rho_i^*\log\mathfrak S_i.
			\label{II:eq:m4-explicit-coarse-energy}
		\end{align}
		The scalar projections satisfy
		$ \langle P_t,P_s\rangle
		=
		E^0_{ts}+O(\varepsilon^2L),$
		with the same estimate after one center or log-scale derivative. Indeed,
		\[
		\int_{\mathbb R^2}e^{u_a}=8\pi,
		\qquad
		\int_{\mathbb R^2}e^{u_a}
		\log(1+a^2|y|^2)\,dy=8\pi,
		\]
		the odd Taylor terms vanish, and the quadratic terms are
		$O(\varepsilon^2L)$. Consequently,
		\[
		\partial_\tau
		\left(J_{\rho^*}(U_z)-\mathcal J_0\right)
		=
		O(\varepsilon^2L),
		\qquad
		\partial_\tau\in
		\left\{
		\partial_{p_t},
		\partial_{\ell_t},
		\ell_t^{-1}\partial_{b_{i,t}}
		\right\}.
		\]
		
		Direct differentiation gives
		\[
		\partial_{\ell_t}\mathcal J_0
		=
		16\pi s_t^TKs_t
		-\sum_i\rho_i^*\Omega_{i,t}(4s_{i,t}-2).
		\]
		Since $K\mathbf1=\sigma/4$, the terms linear in $b$ cancel. Moreover,
		\[
		\ell_t^{-1}\partial_{b_{i,t}}\mathcal J_0
		=
		2\pi\sigma_i
		\left((K\xi_t)_i+\sigma_i\theta_{i,t}\right)
		+O\left(\frac{\delta_z}{L}\right).
		\]
		At $b=0$ and $N\Omega_{i,t}=1$, the remaining free and logarithmic terms
		both equal
		\[
		\frac{\pi\sigma_i^2}{\ell_t}
		\left(
		2\log a-1+4\pi\widetilde G_t(p_t;p)
		\right)
		\]
		and hence cancel. Similarly,
		\[
		\lambda_t^{-1}\partial_{(p_t)_k}\mathcal J_0
		=
		-\frac{1}{N\lambda_t}\partial_{(p_t)_k}\ff(p)
		+O(\varepsilon\delta_z).
		\]
		
		It remains to estimate the correction. Denote
		$
		E_{ts}:=\langle P_t,P_s\rangle
		$
		and let $M_t$ be the truncated mass of the scalar source. Since the
		uncentered first angular mode is orthogonal to each radial source,
		\eqref{II:eq:m4-surface-denominator-exact} gives
		\begin{align}
			J_{\rho^*}(\widehat U)-\mathcal J_0
			={}&
			\frac12\sum_{t,s}s_t^TKs_s(E_{ts}-E^0_{ts})
			-\sum_i\rho_i^*
			\left(
			\sum_t s_{i,t}\gamma_t+\log(1+\eta_i)
			\right)
			\nonumber\\
			&+
			\sum_i
			\left(
			\rho_i^*-\sum_{j,t}a^{ij}s_{j,t}M_t
			\right)
			\overline{\widetilde C_i}
			+\frac12\sum_{i,j}a^{ij}\langle C_i,C_j\rangle .
			\label{II:eq:corrected-energy-identity}
		\end{align}
		Using
		$
		M_t=8\pi+O(\varepsilon^2),
		\rho_i^*-\sum_{j,t}a^{ij}s_{j,t}M_t
		=O(|b|+\varepsilon^2),
		$
		and Lemma~\ref{II:lem:amplitude-projected-size}, the last line contributes
		\[
		O(\varepsilon/L)
		\quad\text{in an amplitude direction},
		\qquad
		O(\varepsilon|b|+\varepsilon^2)
		\quad\text{in a geometric direction}.
		\]
		
		Finally, the term
		$
		W_i:=\sum_t b_{i,t}Q_{i,t}
		$
		contributes $O(|b|\varepsilon)$ for a fixed test. The corresponding changes
		of the tangent vectors have norms
		\[
		O(\varepsilon/L)
		\quad\text{and}\quad
		O(|b|\varepsilon)
		\]
		in the amplitude and geometric directions, respectively. Since
		\[
		\|\nabla J_{\rho^*}(B_z)\|
		\le
		C\left(
		\sqrt L\,\delta_z+\varepsilon+\varepsilon^2L^{3/2}
		\right),
		\]
		these errors are absorbed by the stated remainders. Combining the base
		pairings with
		\eqref{II:eq:m4-verified-amplitude-normal-pairing},
		\eqref{II:eq:m4-verified-geometric-normal-pairing}, and
		\eqref{II:eq:m4-verified-normal-solution} proves the result.
	\end{proof}

	The following identity transfers derivatives from the base family to its
	normal graph.
	
	\begin{lemma}\label{II:lem:ABL-first-order-lift}
		Let $B$ be a $C^2$ immersion into $\mathcal H$, with tangent frame
		$(T_a)$ and Gram matrix
		\[
		\mathbb G:=\bigl((T_a,T_b)\bigr)_{a,b}\ge g_*I,
		\qquad g_*>0.
		\]
		Suppose that
		\[
		G=B+\bar v,
		\qquad
		(\bar v,T_a)=0,
		\qquad
		\nabla J_{\rho^*}(G)\in\operatorname{span}\{T_a\}.
		\]
		For a parameter vector $Y$, set
		$
		d_a:=DJ_{\rho^*}(G)[T_a],
		h_a:=(\bar v,D_YT_a).
		$
		Then
		\begin{align}
			DJ_{\rho^*}(G)[DG[Y]]
			&=
			DJ_{\rho^*}(G)[DB[Y]]
			-d^T\mathbb G^{-1}h,\nonumber
			\\
			\left|DJ_{\rho^*}(G)[D_Y\bar v]\right|
			&\le
			g_*^{-1}|d|\,\|\bar v\|
			\left(\sum_a\|D_YT_a\|^2\right)^{1/2}.\label{II:eq:ABL-first-order-lift-bound}
		\end{align}
	\end{lemma}
	
	\begin{proof}
		Since
		\[
		\nabla J_{\rho^*}(G)=\sum_a(\mathbb G^{-1}d)_aT_a,
		\qquad
		(D_Y\bar v,T_a)=-(\bar v,D_YT_a)=-h_a,
		\]
		we obtain
		$
		DJ_{\rho^*}(G)[D_Y\bar v]
		=
		-d^T\mathbb G^{-1}h.
		$
		The two conclusions follow.
	\end{proof}
	
	
	We now evaluate the reduced energy and its common-scale derivative along a
	fixed balanced ray.
	
	\begin{theorem}\label{II:thm:fixed-balanced-ray}
		Fix a scale-balanced pair $(p,r^*)$ and set
		$
		\delta_t:=\sum_\alpha e_t^\alpha r_\alpha^*,
		\lambda_t:=e^{L+\delta_t},
		b:=0,
		$
		\[
		B(L):=U^\sharp_{p,\ell,0},
		\qquad
		G(L):=B(L)+\bar v(p,\ell,0),
		\qquad
		C_4:=-2\pi N\Big(\sum_i\sigma_i\Big)
		\bigl(1+\log(\pi N)\bigr).
		\]
		For all sufficiently large $L$,
		\begin{align*}
			J_{\rho^*}(G(L))
			&=
			C_4-\frac1N\mathcal F_{\rho^*,h}^N(p)
			-\frac{\pi N}{2}\mathcal L_\infty(p,r^*)Le^{-2L}
			+O(e^{-2L}),
			\\
			\frac{d}{dL}J_{\rho^*}(G(L))
			&=
			\pi N\mathcal L_\infty(p,r^*)Le^{-2L}
			+O(e^{-2L}).
		\end{align*}
		The same expansions hold for $U_{p,\ell,0}$ and $B(L)$, uniformly on compact
		sets of balanced pairs. If $\mathcal L_\infty(p,r^*)\ne0$, then
		\[
		\operatorname{sign}
		\frac{d}{dL}J_{\rho^*}(G(L))
		=
		\operatorname{sign}\mathcal L_\infty(p,r^*)
		\]
		for large $L$. Moreover,
		\[
		\|G(L)\|_{\mathcal H}^2=32\pi nNL+O(1).
		\]
	\end{theorem}
	
	\begin{proof}
		Set
		$
		U:=U_{p,\ell,0},
		d_t:=(a\lambda_t)^{-1},
		j_t:=d_t^2\log\frac1{d_t}.
		$
		All derivatives below are taken with respect to $L$, with $p,r^*$ and $b=0$
		fixed. On the $t$-th core,
		\[
		P_t
		=
		-2\log(|y_{p_t}|^2+d_t^2)
		+8\pi H^g(\cdot,p_t)
		+\gamma_t+\mathfrak e_t,
		\]
		whereas outside the cutoff disk,
		$
		P_t=8\pi G^g(\cdot,p_t)+\gamma_t+\mathfrak e_t.
		$
		Radial integration and the mean-zero Poisson estimate give
		\begin{equation}\label{II:eq:fixed-ray-scalar-refinement}
			\gamma_t
			=
			2\int_M\chi_t
			\log\left(1+\frac{d_t^2}{|y_{p_t}|^2}\right)dV_g
			=
			4\pi j_t+O(e^{-2L}),
			\qquad
			\|\mathfrak e_t\|_{C^2}
			+\|\partial_L\mathfrak e_t\|_{C^2}
			\le Ce^{-2L}.
		\end{equation}
		
		The radial identities
		\[
		\int_{\mathbb R^2}
		\frac{8d_t^2}{(d_t^2+|y|^2)^2}\,dy
		=
		8\pi+O(d_t^2),
		\quad
		\int_{\mathbb R^2}
		\frac{8d^2[-2\log(|y|^2+d^2)]}
		{(d^2+|y|^2)^2}\,dy
		=
		32\pi\log\frac1d-16\pi
		\]
		and the quadratic Taylor terms of the smooth part yield
		\begin{equation}\label{II:eq:fixed-ray-energy-matrix}
			\langle P_t,P_s\rangle
			=
			E^0_{ts}+32\pi^2(j_t+j_s)+O(e^{-2L}).
		\end{equation}
		The same remainder holds after one $L$-derivative.
		
		Set
		$
		f_{i,t}:=\rf_{p,t}^i(p_t),
		T_{i,t}:=T_{i,t}^{\rm raw}(p).
		$
		The denominator expansion is
		\begin{equation*}
			Z_i
			=
			\sum_t\frac{\pi f_{i,t}}{d_t^2}
			+\frac{\pi}{2}
			\sum_t f_{i,t}T_{i,t}\log\frac1{d_t}
			+O(1).
		\end{equation*}
		Indeed,
		\[
		\Delta\log
		\left(
		e^{\widehat\varphi}h_i e^{8\pi\widetilde G_t}
		\right)(p_t)
		+
		\left|
		\nabla\log
		\left(
		e^{\widehat\varphi}h_i e^{8\pi\widetilde G_t}
		\right)(p_t)
		\right|^2
		=
		T_{i,t},
		\]
		and
		\[
		\int_{|y|<r_0}
		\frac{|y|^2}{(|y|^2+d_t^2)^2}\,dy
		=
		2\pi\log\frac1{d_t}+O(1).
		\]
		The linear term vanishes by radial symmetry; the cubic, exterior, and
		projection errors are $O(1)$. These estimates remain valid after one
		$L$-derivative.
		
		The balance condition  gives
		$
		\frac{f_{i,t}}{d_t^2}
		=
		\frac{f_{i,s}}{d_s^2}$ for all $s,t.$
		Therefore
		\[
		\log Z_i
		=
		\log\left(
		\sum_t\frac{\pi f_{i,t}}{d_t^2}
		\right)
		+\frac1{2N}\sum_tT_{i,t}j_t
		+O(e^{-2L}),
		\]
		with the same differentiated remainder. The contributions of
		\eqref{II:eq:fixed-ray-energy-matrix} and
		\eqref{II:eq:fixed-ray-scalar-refinement} proportional to
		$\sum_tj_t$ cancel. Hence
		\[
		J_{\rho^*}(U)
		=
		C_4-\frac1N\ff(p)
		-\pi\sum_{i,t}\sigma_iT_{i,t}j_t
		+O(e^{-2L}).
		\]
		Since
		$
		e^{\fJ_i}=\frac{2\sigma_i}{a^2},
		\mathcal L_\infty(p,r^*)
		=
		\frac4N\sum_{i,t}
		e^{-2\delta_t^*}L_{i,t}^{\rm corr}(p),$
		we obtain
		\[
		-\pi\sum_{i,t}\sigma_iT_{i,t}j_t
		=
		-\frac{\pi N}{2}
		\mathcal L_\infty(p,r^*)Le^{-2L}
		+O(e^{-2L}).
		\]
		Differentiation gives
		\[
		\frac{d}{dL}J_{\rho^*}(U)
		=
		\pi N\mathcal L_\infty(p,r^*)Le^{-2L}
		+O(e^{-2L}).
		\]
		The constant term follows from
		\[
		\mathbf1^TA^{-1}\mathbf1
		=
		\frac14\sum_i\sigma_i,
		\qquad
		\frac1N\sum_t\log f_{i,t}
		=
		\log\left(f_{i,t}e^{2\delta_t}\right),
		\]
		together with the double counting of the off-diagonal Green interactions.
		
		We next compare $U$ and $B=U+C$. By
		Lemma~\ref{II:lem:amplitude-projected-size},
		\[
		\|C\|_{H^1\cap L^\infty}
		+
		\|\partial_LC\|_{H^1\cap L^\infty}
		=
		O(e^{-L}).
		\]
		The angular orthogonality and the decay of $\Phi$ imply
		\[
		J_{\rho^*}(B)-J_{\rho^*}(U)=O(e^{-2L}),
		\qquad
		\frac{d}{dL}
		\left(
		J_{\rho^*}(B)-J_{\rho^*}(U)
		\right)
		=
		O(e^{-2L}).
		\]
On the balanced ray,
\[
\Omega_{i,t}=\frac1N,\qquad
\delta_z=0,\qquad
\|\bar v\|=O(e^{-2L}L^{3/2})
\]
by Theorem~\ref{II:thm:m4-verified-normal-inverse}. Hence,
$J_{\rho^*}(G)-J_{\rho^*}(B)=O(e^{-4L}L^3).$
Since $\partial_LB=\sum_tS_t$,
\eqref{II:eq:m4-verified-geometric-normal-pairing} gives
\begin{align*}
\left|
\bigl(DJ_{\rho^*}(G)-DJ_{\rho^*}(B)\bigr)[\partial_LB]
\right|
&\le
C\left(
e^{-L}\sqrt L\,\|\bar v\|+\|\bar v\|^2
\right)=
O(e^{-3L}L^2+e^{-4L}L^3).
\end{align*}

Let
$
d_a:=DJ_{\rho^*}(G)[T_a].
$
The amplitude, scale, and position pairings in
Theorem~\ref{II:thm:m4-corrected-leading-pairings} give
\[
|d|=O(e^{-L}).
\]
Since the normalized frame has bounded common-scale derivative,
\eqref{II:eq:ABL-first-order-lift-bound} yields
\[
\left|DJ_{\rho^*}(G)[\partial_L\bar v]\right|
\le C|d|\,\|\bar v\|
=O(e^{-3L}L^{3/2}).
\]
Therefore
\begin{align*}
\left|
\frac{d}{dL}J_{\rho^*}(G(L))
-\frac{d}{dL}J_{\rho^*}(B(L))
\right|
&=
O(e^{-3L}L^2+e^{-4L}L^3+e^{-3L}L^{3/2})=o(e^{-2L}).
\end{align*}
This proves the energy and derivative expansions for $G(L)$. Finally,
the scalar Gram formula gives
\[
\|U\|_{\mathcal H}^2=32\pi nNL+O(1),
\qquad
\|G-U\|=O(e^{-L}),
\]
and hence
\[
\|G\|_{\mathcal H}^2=32\pi nNL+O(1).
\]
	\end{proof}
	
	\subsection{The amplitude and relative-scale center}
	
	
	Fix a compact, separated location set and compact windows for $r$ and
	$\beta^0$. In both mass regimes, denote
	\begin{equation}\label{II:eq:bounded-amplitude-reconstruction}
		\delta_t=\sum_\alpha e_t^\alpha r_\alpha,\qquad
		\ell_t=L+\delta_t,\qquad
		b_{i,t}
		=
		\frac1L\left(
		\beta_i^0+\sum_\alpha e_t^\alpha c_i^\alpha
		\right).
	\end{equation}
	Set $H_{i,t}:=\log\rf_{p,t}^i(p_t)$.
	
	The relative mass equations can be solved uniformly for the variables
	$c^\alpha$.
	
	\begin{theorem}\label{II:thm:bounded-complete-mass-inverse}
		Fix $0<\delta_0<1/2$ and suppose that
		$\sum_t\theta_{i,t}=0$ and
		$\max_{i,t}|\theta_{i,t}|\le\delta_0$. For all sufficiently large $L$,
		there is a unique solution
		$c=\mathscr C_L(p,\beta^0,r,\theta)$ near
		\begin{equation*}
			(\mathscr C_\infty)_i^\alpha
			=
			\frac1{\sigma_i}
			\left[
			\sum_t e_t^\alpha\log(1-\theta_{i,t})
			-H_i^\alpha(p)-q_ir_\alpha
			\right]
		\end{equation*}
		such that the complete masses in
		\eqref{II:def:full-local-mass} have centered defects $\theta$.
		
		On the fixed compact parameter sets, $\mathscr C_L$ is $C^2$ and
		\begin{align*}
			\|\mathscr C_L-\mathscr C_\infty\|_{C^2}
			&\le \frac{C}{L},&
			\|\partial_L\mathscr C_L\|_{C^1}
			&\le \frac{C}{L^2},&
			|\partial_L^2\mathscr C_L|
			&\le \frac{C}{L^3}.
		\end{align*}
		The $C^2$ norm is taken with respect to $(p,\beta^0,r,\theta)$ in fixed
		finite-dimensional coordinates. In particular, at $\theta=0$,
		\begin{equation*}
			H^\alpha(p)+qr_\alpha+D_\sigma c^\alpha=O(L^{-1}).
		\end{equation*}
	\end{theorem}
	
	\begin{proof}
		Set $u:=L^{-1}$. Taking the relative logarithmic modes in
		\eqref{II:def:full-local-mass} cancels $q_iL$, $\sigma_i\beta_i^0$,
		and the common prefactor. The remaining equations are
		\begin{equation*}
			H_i^\alpha+q_ir_\alpha+\sigma_i c_i^\alpha
			+uE_i^\alpha(u,p,\beta^0,r,c)
			=
			\sum_t e_t^\alpha\log(1-\theta_{i,t}).
		\end{equation*}
		The relative logarithms and the condition $\sum_t\Omega_{i,t}=1$
		determine the positive fractions uniquely.
		
		The error $E$ extends smoothly to $u=0$, with uniformly bounded parameter
		derivatives. Indeed, the powers contribute
		$\sigma_i\delta_t(\beta_i^0+\sum_\gamma e_t^\gamma c_i^\gamma)$, while
		the $k$-th derivatives of the logarithmic prefactors are bounded by
		\[
		C(1+|y|)^{-m_i+\eta}
		\bigl(1+\log^k(2+|y|)\bigr),
		\qquad
		0<\eta<\widehat m-2.
		\]
		At $u=0$, the Jacobian with respect to $c$ is $D_\sigma$ in each relative
		mode. The parameter-dependent implicit function theorem gives the stated
		bounds. The $L$-derivative estimates follow from
		$\partial_L=-u^2\partial_u$.
	\end{proof}

	Let $c:=\mathscr C_L(p,\beta^0,r,\theta)$ in
	\eqref{II:eq:bounded-amplitude-reconstruction}. Since
	$\sum_t\theta_t=0$, denote
	$\theta_t=\sum_\alpha e_t^\alpha\theta^\alpha$. The leading covectors are
	\begin{equation}\label{II:eq:principal-centered-covectors}
		\begin{aligned}
			C^0=\frac1L KD_\sigma\beta^0,\qquad
			C^\alpha=\frac1L KD_\sigma c^\alpha+D_\sigma\theta^\alpha,\qquad
			R^\alpha=q^TD_\sigma\theta^\alpha,
			\qquad 1\le\alpha<N.
		\end{aligned}
	\end{equation}
	
	The following result locates their common zero and gives a uniform estimate
	away from it.
	
	\begin{theorem}\label{II:thm:finite-scale-covector-center}
		Fix a compact separated location set contained in one coordinate chart. Set
		\begin{align*}
			r_{\infty,\alpha}(p)
			&:=
			-\frac{q^TKH^\alpha(p)}{q^TKq},&
			c_\infty^\alpha(p)
			&:=
			-D_\sigma^{-1}
			\bigl(H^\alpha(p)+qr_{\infty,\alpha}(p)\bigr).
		\end{align*}
		Choose bounded windows for $r$ and $c$ containing these limiting vectors
		with positive margins.
		
		For all sufficiently large $L$, the covectors in
		\eqref{II:eq:principal-centered-covectors} have a unique common zero near
		$(\beta^0,r,\theta)=(0,r_\infty,0)$. Denote it as
		$(0,r_\circ(p,L),\theta_\circ(p,L))$ and set
		$c_\circ^\alpha
		:=\mathscr C_L(p,0,r_\circ,\theta_\circ)^\alpha$. Then
		\begin{align}
			r_\circ
			&=r_\infty+O(L^{-1}),&
			c_\circ
			&=c_\infty+O(L^{-1}),&
			\theta_\circ^\alpha
			&=
			-\frac1L D_\sigma^{-1}KD_\sigma c_\infty^\alpha
			+O(L^{-2}).
			\label{II:eq:principal-center-asymptotics}
		\end{align}
		The center is $C^1$ at finite scale and satisfies
		\[
		|\partial_Lr_\circ|
		+|\partial_Lc_\circ|
		+|\partial_L\theta_\circ|
		\le CL^{-2}.
		\]
		
		On a fixed sufficiently small convex neighborhood of this zero, with
		$\beta^0$ in a fixed bounded ball, set
		\begin{equation*}
			d_L
			:=
			\frac{|\beta^0|}{L}
			+\frac{|r-r_\circ|}{L}
			+|\theta-\theta_\circ|,
			\qquad
			|\mathcal C|
			:=
			|C^0|+\sum_\alpha
			\bigl(|C^\alpha|+|R^\alpha|\bigr).
		\end{equation*}
		There exist $c,C>0$, independent of $L$, such that
		\begin{equation}\label{II:eq:finite-scale-relative-coercivity}
			cd_L\le|\mathcal C|\le Cd_L,
			\qquad
			|b-b_\circ|+|\theta-\theta_\circ|
			\le C|\mathcal C|,
		\end{equation}
		where
		$b_{\circ,t}:=L^{-1}\sum_\alpha e_t^\alpha c_\circ^\alpha$. If $p$ is scale-balanced, then
		$
		r_\circ=r^*(p),
		c_\circ=0,
		\theta_\circ=0.
		$
	\end{theorem}
	
	\begin{proof}
		Set $u:=L^{-1}$. In each relative mode, the exact mass equation is
		\begin{equation*}
			H^\alpha+qr_\alpha+D_\sigma c^\alpha+e_L^\alpha
			=
			\sum_t e_t^\alpha\log(1-\theta_t),
		\end{equation*}
		where Theorem~\ref{II:thm:bounded-complete-mass-inverse} gives
		$\|e_L\|_{C^2}=O(u)$. Moreover,
		\begin{align*}
			D_rc
			&=-D_\sigma^{-1}qI+O(u),&
			D_{\beta^0}c
			&=O(u),&
			D_\theta c
			&=O(1).
		\end{align*}
		
		The equation $C^0=0$ gives $\beta^0=0$. For fixed $r$, the equation
		$C^\alpha=0$ becomes
		\[
		\theta^\alpha
		=
		-uD_\sigma^{-1}KD_\sigma c^\alpha.
		\]
		It has a unique solution in $|\theta|\le Cu$ by contraction. Along this
		solution,
		\begin{equation*}
			\frac{R^\alpha}{u}
			=
			-q^TKD_\sigma c^\alpha
			=
			q^TK\bigl(H^\alpha+qr_\alpha\bigr)+O(u)
		\end{equation*}
		in $C^1(r)$. The limiting Jacobian with respect to $r_\alpha$ is
		$(q^TKq)I$, where
		$
		q^TKq
		=
		4\mathbf1^TK\mathbf1
		>0.
		$
		The implicit function theorem at $u=0$ gives the unique common zero and
		\eqref{II:eq:principal-center-asymptotics}. Since
		$\partial_L=-u^2\partial_u$, its bounded $u$-derivatives give the stated
		$L$-derivative estimate.
		
		In the variables
		$\bigl(\beta^0/L,(r-r_\circ)/L,\theta-\theta_\circ\bigr)$, the Jacobian
		of $(C^0,C^\alpha,R^\alpha)$ is
		\begin{equation*}
			\begin{pmatrix}
				KD_\sigma&0&0\\
				0&-Kq&D_\sigma\\
				0&0&q^TD_\sigma
			\end{pmatrix}
			+O(L^{-1}),
		\end{equation*}
		with one relative block for each $\alpha$. The common block
		$KD_\sigma$ is invertible. If $(\tau,v)$ belongs to the kernel of a
		relative block, then
		$
		D_\sigma v=Kq\,\tau,
		q^TD_\sigma v=0,
		$
		and hence $(q^TKq)\tau=0$. Thus $\tau=0$ and $v=0$.
		
		Integrating the Jacobian along the segment from its zero in the convex
		coordinate neighborhood gives
		\eqref{II:eq:finite-scale-relative-coercivity}. If $p$ is scale-balanced,
		the exact masses agree at $b=0$ and $r=r^*(p)$. Uniqueness then gives
		$r_\circ=r^*(p)$, $c_\circ=0$, and $\theta_\circ=0$.
	\end{proof}

	\begin{remark} The amplitude center need not vanish away from a balanced location. 
		Let $M$ be the flat unit torus, $n=N=2$,
		\[
		A_{11}=A_{22}=0,
		\qquad
		A_{12}=A_{21}=1,
		\qquad
		\sigma=(4,4).
		\]
		Take $p_2=p_1+(1/2,0)$ and $h_{1,2}=e^{\pm\psi}$, where
		$\psi(x+(1/2,0))=-\psi(x)$ and $0<|\psi(p_1)|\ll1$. Translation
		invariance of the Green data gives a nonzero relative logarithmic weight
		proportional to $(1,-1)$. Hence
		\[
		r_\infty=0,
		\qquad
		c_\infty\ne0,
		\qquad
		|b_\circ|\asymp L^{-1}
		\]
		at the common zero of the covectors $\mathcal C=
\bigl(C^0,(C^\alpha,R^\alpha)_{\alpha=1}^{N-1}\bigr)$ defined in 
\eqref{II:eq:principal-centered-covectors}.
	\end{remark}
	

\section{Transverse critical points and common-scale derivatives}\label{sec:4}

\subsection{Transverse critical points at the borderline exponent}

Assume $m_i=4$ for every $i$ and \textup{(A1)}. Fix $\rho^*$ and solve
the transverse equations on a slice of constant mean log-scale. Set
\[
f(p,r,b,L):=J_{\rho^*}(G(p,r,b,L)),
\qquad
G=B_z+\bar v,
\qquad
B_z=U_z^\sharp.
\]
Fix a balanced critical pair $(p^*,r^*)$ such that $D^2\ff(p^*)$ is
invertible. Work in an ordered separated chart and denote $x:=p-p^*$ in
local coordinates.

For each $(p,L)$, let $(r_\circ,\theta_\circ)$ be the center given by
Theorem~\ref{II:thm:finite-scale-covector-center}. Define
\begin{equation*}
\zeta
:=
\left(
\frac{\beta^0}{L},
\frac{r-r_\circ}{L},
\theta-\theta_\circ
\right),
\qquad
\mathcal D_L
:=
\left\{
|x|\le\varepsilon L^{5/4},
\quad
|\zeta|\le\varepsilon L^{-1/2}
\right\}.
\end{equation*}
Here $b$ is reconstructed by
Theorem~\ref{II:thm:bounded-complete-mass-inverse}. For $N=1$, the
relative blocks are absent.

The next result gives the transverse critical points and the sign of the
remaining common-scale derivative.

\begin{theorem}
\label{II:thm:slice-signed}
For all sufficiently large $L$, the interior of $\mathcal D_L$ contains a
zero of $D_{p,r,b}f$. Every such zero satisfies
\[
|p-p^*|=O(\varepsilon L),
\qquad
|\zeta|=O(\varepsilon/L).
\]
Moreover, uniformly on $\mathcal D_L$,
\begin{equation}\label{II:eq:slice-signed}
\partial_Lf(p,r,b,L)
=
\pi N\mathcal L_\infty(p^*,r^*)L\varepsilon^2
+O(\varepsilon^2L^{3/4}),
\end{equation}
where $p,r,b$ are fixed in $\partial_Lf$. Hence, if
$\mathcal L_\infty(p^*,r^*)\ne0$, then
\[
\operatorname{sign}(\partial_Lf)
=
\operatorname{sign}\mathcal L_\infty(p^*,r^*)
\]
at every transverse zero and for large $L$. At such a zero, all lifts of
unit common-scale velocity have the same pairing.
\end{theorem}

\begin{proof}
Let
$\mathcal C=(C^0,(C^\alpha,R^\alpha)_\alpha)$ be the covectors in
\eqref{II:eq:principal-centered-covectors}. The proof of
Theorem~\ref{II:thm:finite-scale-covector-center} gives
\[
D_pr_\circ=O(1),
\qquad
D_pc_\circ=O(1),
\qquad
D_p\theta_\circ=O(L^{-1}).
\]
At $p=p^*$, one has
$r_\circ=r^*$ and $c_\circ=\theta_\circ=0$.

By \eqref{II:eq:finite-scale-relative-coercivity} and
\eqref{II:eq:bounded-amplitude-reconstruction}, throughout $\mathcal D_L$,
\begin{equation}\label{II:eq:slice-delta}
\delta:=|b|+|\theta|
\le
C\left(|\zeta|+\frac{|x|}{L}\right)
=
O(\varepsilon L^{1/4}).
\end{equation}
Moreover,
\begin{equation*}
|r-r^*|
\le C\bigl(|x|+L|\zeta|\bigr)
=
O(\varepsilon L^{5/4}),
\qquad
\|\bar v\|
=
O(\varepsilon L^{1/4}).
\end{equation*}
Thus $\mathcal D_L$ lies in the region of small mass defect and bounded
$L|b|$.

At fixed $(p,L)$, integration of the Jacobian in
Theorem~\ref{II:thm:finite-scale-covector-center} gives
$\mathcal C=M\zeta+O(|\zeta|/L)$, where $M$ is invertible. Its common
block is $KD_\sigma$, and each relative block is
\[
(u,v)\longmapsto
\bigl(-Kq\,u+D_\sigma v,\ q^TD_\sigma v\bigr).
\]

Apply Theorem~\ref{II:thm:m4-corrected-leading-pairings} to the fixed
directions at $B_z$. Divide the amplitude rows by $2\pi\sigma_i$, take
their mean and relative modes, and divide the relative-scale rows by
$2\pi$. The resulting covector is
\begin{equation}\label{II:eq:slice-amperror}
\mathcal C
+
O\left(
\frac{\delta+\varepsilon}{L}
+\delta^2+\varepsilon^2L+\|\bar v\|^2
\right)
=
\mathcal C+O(\varepsilon L^{-3/4}).
\end{equation}
The normalized position rows are
\[
-\frac1{N\lambda_t}\partial_{p_t}\ff(p)
+O(\varepsilon^2L),
\]
while the individual leading scale rows are $O(|\zeta|)$.

The full vector of normalized frame pairings has norm $O(\varepsilon)$.
For every normalized amplitude, scale, or position direction $Y$,
$
\|D_YT_a\|=O(1);
$
see \cite[Lemma~A.2]{AFHZI} and
Lemma~\ref{II:lem:amplitude-projected-size}. Therefore
Lemma~\ref{II:lem:ABL-first-order-lift} gives a graph-lift error
$
O(\varepsilon\|\bar v\|)
=
O(\varepsilon^2L^{1/4}).
$
After multiplication by $-N\lambda_t$, the lifted position rows become
$
\partial_{p_t}\ff(p)+O(\varepsilon L).
$

Let $H:=D^2\ff(p^*)$. After the preceding changes of basis, the
transverse equations are equivalent to the zero set of
\begin{equation*}
(x,\zeta)\mapsto
\left(
Hx+O(|x|^2+\varepsilon L),
\;
M\zeta+
O\left(\frac{|\zeta|}{L}+\varepsilon L^{-3/4}\right)
\right).
\end{equation*}
On the two boundary faces, the errors are respectively
$o(\varepsilon L^{5/4})$ and $o(\varepsilon L^{-1/2})$. The homotopy to
$(Hx, M\zeta)$ has no boundary zero, and its degree is
$\operatorname{sign}(\det H\det M)\ne0$. Hence $\mathcal D_L$ contains a
transverse zero.

At every such zero, invertibility of $H$ gives
$|x|=O(\varepsilon L)$. Equations
\eqref{II:eq:slice-delta} and \eqref{II:eq:slice-amperror}, together with
the graph-lift estimate, give
\[
\delta=O(\varepsilon),
\qquad
|\zeta|=O(\varepsilon/L).
\]

We now estimate the common-scale derivative. Here
$\partial_L=\sum_t\partial_{\ell_t}$ with $p,r,b$ fixed. A superscript
$*$ denotes evaluation at $(p^*,r^*,0,L)$.

In \eqref{II:eq:m4-explicit-coarse-energy}, denote
$s_t=\mathbf1+\xi_t/4$. Summing the log-scale derivatives gives
\begin{equation}\label{II:eq:slice-quadratic}
\partial_L\mathcal J_0
=
\pi\sum_t\xi_t^TK\xi_t
+
2\pi\sum_t\xi_t^TD_\sigma\theta_t.
\end{equation}
Indeed,
\[
K\mathbf1=\frac{\sigma}{4},
\qquad
\rho_i^*\Omega_{i,t}
=
2\pi\sigma_i(1-\theta_{i,t}),
\qquad
\sum_t\theta_{i,t}=0.
\]
Thus \eqref{II:eq:slice-quadratic} vanishes at the reference ray and is
$O(\varepsilon^2L^{1/2})$ on $\mathcal D_L$.

We next use \eqref{II:eq:corrected-energy-identity}. The common-scale
derivative of its last line is
\[
O\bigl((|b|+\varepsilon^2)\varepsilon+\varepsilon^2\bigr).
\]
Equations \eqref{II:eq:fixed-ray-scalar-refinement} and
\eqref{II:eq:fixed-ray-energy-matrix} hold uniformly before imposing
balance. Hence the change from the reference ray in the differentiated
leading terms is
\[
O\left(
\varepsilon^2L\bigl(|b|+|r-r^*|\bigr)
\right),
\]
and the differentiated remainders are $O(\varepsilon^2)$.

The mixed derivative estimates in
\eqref{II:eq:m4-surface-denominator-exact} give
\begin{equation*}
\left|
\partial_L\log(1+\eta_i)
-
\bigl(\partial_L\log(1+\eta_i)\bigr)^*
\right|
\le
C\varepsilon^2
\left[
L\bigl(|x|+|r-r^*|\bigr)+L^2|b|
\right]
=
o(\varepsilon^2).
\end{equation*}

Finally, denote $W_i:=\sum_t b_{i,t}Q_{i,t}$, so that
$B_z=\widehat U+W$. Lemma~\ref{II:lem:amplitude-projected-size} gives
\[
\|W\|_{H^1\cap L^\infty}
+
\|\partial_LW\|_{H^1}
=
O(|b|\varepsilon).
\]
Since $\|DJ_{\rho^*}(B_z)\|=O(\sqrt L)$, the correction $W$ changes the
common-scale derivative by $O(\varepsilon^2L^{3/4})$ on $\mathcal D_L$.

Summing \eqref{II:eq:m4-verified-geometric-normal-pairing} over the
log-scale directions gives
\begin{equation*}
\left|
\bigl(
DJ_{\rho^*}(G)-DJ_{\rho^*}(B_z)
\bigr)[\partial_LB_z]
\right|
\le
C\left[
(\delta+\varepsilon\sqrt L)\|\bar v\|
+\|\bar v\|^2
\right]
=
O(\varepsilon^2L^{3/4}).
\end{equation*}
The graph-lift error is $O(\varepsilon^2L^{1/4})$. Therefore
\[
\partial_Lf-(\partial_Lf)^*
=
O(\varepsilon^2L^{3/4}).
\]
Theorem~\ref{II:thm:fixed-balanced-ray} gives
\[
(\partial_Lf)^*
=
\pi N\mathcal L_\infty(p^*,r^*)L\varepsilon^2
+O(\varepsilon^2),
\]
which proves \eqref{II:eq:slice-signed}. At a transverse zero, two lifts
of unit common-scale velocity differ by a transverse vector, so their
pairings agree.
\end{proof}

\subsection{The common-scale derivative for unequal decay exponents}

Fix an entire profile with $2<\widehat m<4$ and a compact separated
parameter cell with comparable scales. We use the original base $U_z$ and
its normal fiber $E_z$. Along a fixed balanced ray, denote
$U:=U_{p,\ell,0}$ and $E:=E_z$.

For $2<m_i<4$, choose a partition $(\Omega_t)$ such that $\Omega_t$
contains a coordinate disk about $p_t$ and no other center. Define
for $2<m_i<4$
\begin{equation*}
D_{i,t}
:=
\lim_{d\downarrow0}
\left[
\int_{\Omega_t\setminus U_d(p_t)}
\frac{
h_i(x)\exp\bigl(2\pi m_i\sum_sG^g(x,p_s)\bigr)
}{
\rf_{p,t}^i(p_t)
}\,dV_g
-
\frac{2\pi}{q_i}d^{-q_i}
\right].
\end{equation*}
Set
$F_i(p,x):=h_i(x)\exp(2\pi m_i\sum_sG^g(x,p_s))$. At a balanced pair,
\[
k_i^*
:=
\rf_{p,t}^i(p_t)e^{q_i\delta_t^*}
\]
is independent of $t$. For $0<q_i<2$ and disjoint disks of radius $d$,
\begin{equation*}
\sum_t e^{-q_i\delta_t^*}D_{i,t}
=
\frac1{k_i^*}
\int_{M\setminus\bigcup_tU_d(p_t)}
F_i(p,x)\,dV_g
-
\frac{2\pi}{q_i}d^{-q_i}
\sum_t e^{-q_i\delta_t^*}
+
O(d^{2-q_i}).
\end{equation*}
Indeed, on the $t$-th disk,
\[
F_i(p,x)\,dV_g
=
|y|^{-2-q_i}
\left(
\rf_{p,t}^i(p_t)+a_{i,t}\cdot y+O(|y|^2)
\right)dy.
\]
The linear term has zero angular integral, while the remainder contributes
$O(d^{2-q_i})$. The constants are uniform on compact separated sets of
balanced pairs. Hence the weighted sum is independent of the partition
and continuous in the parameters.

By \cite[Proposition~3.3]{AFHZI},
\begin{equation}\label{II:eq:mixed-coefficient}
\mathcal L_\infty(p,r^*)
=
\frac{\widehat m-2}{\pi N}
\sum_{i\in\widehat I}e^{\fJ_i}
\sum_t
e^{-(\widehat m-2)\delta_t^*}D_{i,t}.
\end{equation}

We first derive the common-scale expansion on a moving parameter region.

\begin{theorem}\label{II:thm:moving-mixed}
Assume
\[
2<\widehat m<4,
\qquad
\frac{\widehat m-2}{2}
<
\alpha
<
\min\{1,\widehat m-2\}.
\]
Let $(p,r,b,L)$ lie in a compact separated chart with comparable scales and
satisfy
\begin{equation*}
|b|+\Delta_\Omega\le C_0e^{-\alpha L},
\qquad
(p,r)\to (p^*,r^*)
\quad\text{as }L\to\infty.
\end{equation*}
Then $(p^*,r^*)$ is balanced, and the normal graph
$U_z+\bar v_z$ exists for large $L$. For
\[
f(z)=J_{\rho^*}(U_z)
\qquad\text{or}\qquad
f(z)=J_{\rho^*}(U_z+\bar v_z),
\]
one has
\begin{equation*}
\partial_Lf(p,r,b,L)
=
\pi N\mathcal L_\infty(p^*,r^*)
e^{-(\widehat m-2)L}
+
o\left(e^{-(\widehat m-2)L}\right),
\end{equation*}
where $p,r,b$ are fixed in $\partial_Lf$.

If $\mathcal L_\infty(p^*,r^*)\ne0$, then in a sufficiently small
neighborhood of $(p^*,r^*)$,
\begin{equation*}
\operatorname{sign}(\partial_Lf)
=
\operatorname{sign}\mathcal L_\infty(p^*,r^*),
\qquad
|\partial_Lf|
\ge
ce^{-(\widehat m-2)L},
\end{equation*}
uniformly under the preceding defect bound and for large $L$.
\end{theorem}

\begin{proof}
Denote
\[
\psi_i:=\sum_j a^{ij}v_j,
\qquad
v_{i,t}^b:=v_i+b_{i,t}\psi_i,
\qquad
a_{i,t}:=m_i+\xi_{i,t},
\]
and
$
K_{i,t}^b
:=
a_{i,t}+Y\cdot\nabla v_{i,t}^b.
$
The radial equations give
\[
K_{i,t}^b
=
O\left((1+|Y|)^{-(\widehat m-2)}\right),
\]
and integration by parts yields
\begin{equation*}
\int_{\mathbb R^2}e^{v_{i,t}^b}K_{i,t}^b\,dY
=
2\pi\sigma_{i,t}(b)(q_i+\xi_{i,t}).
\end{equation*}
Since $L|b|\to0$, the ratios of the complete masses imply balance of
$(p^*,r^*)$.

Set
\[
k_i:=\frac{\mathfrak S_i}{\rho_i^*},
\qquad
M_{i,t}:=\rho_i^*\Omega_{i,t},
\qquad
d\nu_i:=\frac{h_ie^{U_{z,i}}}{k_i}\,dV_g.
\]
On the $t$-th inner disk, \cite[Lemma~A.1]{AFHZI} gives
\begin{equation}\label{II:eq:mixed-inner}
d\nu_i
=
\frac{M_{i,t}}{2\pi\sigma_{i,t}(b)}
e^{v_{i,t}^b(Y)}
A_{i,t}(Y/\lambda_t)
e^{H_{i,t,L}(x)}\,dY.
\end{equation}
Here $A_{i,t}$ is the smooth centered factor obtained from
$e^{\widehat\varphi_{p_t}}\rf_{p,t}^i\rg_{p,b}^i$, and
\[
\|A_{i,t}\|_{C^2}\le C,
\qquad
\|H_{i,t,L}\|_{C^1}
=
O(\varepsilon^{\widehat m-2}).
\]

Set
\[
T_i:=\sum_tM_{i,t}H_{i,t,L}(p_t),
\qquad
K_{i,L}^b
:=
\sum_t\chi_tK_{i,t}^b(\lambda_ty_{p_t}).
\]
For $i\in\widehat I$, denote 
\[
d_i
:=
e^{\fJ_i}
\sum_t
e^{-(\widehat m-2)\delta_t^*}D_{i,t}(p^*),
\]
and set $d_i:=0$ otherwise. Then
\begin{align*}
\int_Md\nu_i
&=
\rho_i^*+T_i
+\varepsilon^{\widehat m-2}d_i
+o(\varepsilon^{\widehat m-2}),\\
\int_MK_{i,L}^b\,d\nu_i
&=
\sum_tM_{i,t}(q_i+\xi_{i,t})
+q_iT_i
+o(\varepsilon^{\widehat m-2}).
\end{align*}

To prove these expansions, split the integrals at a fixed coordinate radius
$d$. The linear Taylor term of $A_{i,t}$ has zero angular integral. Since
$L|\xi_{i,t}|=o(1)$,
$|Y|^{-\xi_{i,t}}\to1$ uniformly for
$1\le|Y|\le d\lambda_t$. The omitted radial mass is
\begin{equation*}
\frac{
2\pi
e^{\fJ_i+b_{i,t}\sum_j a^{ij}\fJ_j}
}{
q_i+\xi_{i,t}
}
(d\lambda_t)^{-q_i-\xi_{i,t}}
+
O_d\left(
\varepsilon^{q_i+\widehat m-2}e^{O(L|b|)}
\right).
\end{equation*}
For the minimal components, the interior quadratic remainder is
$
O(
\varepsilon^{\widehat m-2}
d^{2-(\widehat m-2)}
).
$
If $q_i>\widehat m-2$, it is bounded by
\begin{equation*}
C_d\varepsilon^2
\left(
1+\int_1^{d\lambda_t}r^{1-q_i}\,dr
\right)
=
\begin{cases}
O_d(\varepsilon^{q_i}),&q_i<2,\\
O_d(L\varepsilon^2),&q_i=2,\\
O_d(\varepsilon^2),&q_i>2.
\end{cases}
\end{equation*}
Each term is $o(\varepsilon^{\widehat m-2})$.

The variation of $H_{i,t,L}$ about $p_t$ contributes
$o(\varepsilon^{\widehat m-2})$, while its square contributes
$O(\varepsilon^{2(\widehat m-2)})$. On the fixed exterior,
\[
U_{z,i}
\to 
2\pi m_i\sum_sG^g(\cdot,p_s^*)
\]
uniformly, and
\[
e^{-q_iL}k_i
\to e^{-\fJ_i}
\rf_{p^*,t}^i(p_t^*)e^{q_i\delta_t^*}.
\]

All preceding estimates are obtained for fixed \(d>0\); we first let
\(L\to\infty\) and then let \(d\downarrow0\).

The exterior contribution minus the omitted tails is therefore $d_i$.
The additional decay in the second integral makes its exterior and tail
terms $o(\varepsilon^{\widehat m-2})$. This proves the two expansions.

Let $S:=\partial_LU_z$. Replacing $S_i$ by $K_{i,L}^b$ in the difference
between the free and normalized sources changes the pairing by
$o(\varepsilon^{\widehat m-2})$. Indeed, the test error is
$O(\varepsilon^{\widehat m-2})$, while the source difference has total
variation
\[
O\bigl(|b|+\Delta_\Omega+\varepsilon^\alpha\bigr).
\]

Denote $a_t=(a_{i,t})_i$. On each whole plane, the free contribution is
$\pi a_t^TKa_t$. The constant part gives $2\pi a_t^TKa_t$, while radial
integration gives
\[
-\pi
\left[
\bigl(r(v_t^b)'\bigr)^T
K\bigl(r(v_t^b)'\bigr)
\right]_0^\infty
=
-\pi a_t^TKa_t.
\]
The discarded tails are
$O(\varepsilon^{\widehat m-2+\min_iq_i})$. Dividing the two preceding
mass expansions cancels $T_i$, and hence
\begin{equation*}
DJ_{\rho^*}(U_z)[S]
=
\pi\sum_ta_t^TKa_t
-
\sum_{i,t}M_{i,t}(q_i+\xi_{i,t})
+
\varepsilon^{\widehat m-2}\sum_iq_id_i
+
o(\varepsilon^{\widehat m-2}).
\end{equation*}
Using resonance, $Km=\sigma$, and $\sum_t\theta_{i,t}=0$, the first two
terms reduce to
\[
\pi\sum_t\xi_t^TK\xi_t
+
2\pi\sum_t\xi_t^TD_\sigma\theta_t
=
O(\varepsilon^{2\alpha}).
\]
Moreover, \eqref{II:eq:mixed-coefficient} gives
\[
\sum_iq_id_i
=
\pi N\mathcal L_\infty(p^*,r^*).
\]
This proves the expansion for $J_{\rho^*}(U_z)$.

It remains to pass to the normal graph. Equation
\eqref{II:eq:mixed-inner}, the separate amplitude anchors, and
$\alpha<\widehat m-2$ give
\[
\|P_z\nabla J_{\rho^*}(U_z)\|
\le C\varepsilon^\alpha.
\]
Lemma~\ref{II:lem:normal-model} and the contraction argument in
Theorem~\ref{II:thm:m4-verified-normal-inverse} give a unique solution
$\bar v_z$ satisfying
\[
P_z\nabla J_{\rho^*}(U_z+\bar v_z)=0,
\qquad
\|\bar v_z\|\le C\varepsilon^\alpha.
\]
The solution is $C^1$ in the finite-scale parameters.

For every $w\in E_z$ with $\|w\|=1$,
\[
\left|
D^2J_{\rho^*}(U_z)[S,w]
\right|
\le C\varepsilon^\alpha.
\]
Indeed, the free scale source differs from $V_{i,t}Z_{0,i}$ by $O(|b|)$,
while the covariance terms
$V_{i,t}(Z_{0,i}+q_i)$ cancel by the anchors. Taylor expansion gives
\[
\left|
\bigl(
DJ_{\rho^*}(U_z+\bar v_z)-DJ_{\rho^*}(U_z)
\bigr)[S]
\right|
\le
C\left(
\varepsilon^\alpha\|\bar v_z\|+\|\bar v_z\|^2
\right)
=
O(\varepsilon^{2\alpha}).
\]
For an ambient unit test, its $V_{i,t}$-mean is $O(\sqrt L)$, whereas its
centered part satisfies the anchored estimates. Hence
\begin{equation}\label{II:eq:mixed-residual}
\|P_z\nabla J_{\rho^*}(U_z)\|
\le C\varepsilon^\alpha,
\qquad
\|\nabla J_{\rho^*}(U_z)\|
\le C\sqrt L\,\varepsilon^\alpha.
\end{equation}
Lemma~\ref{II:lem:ABL-first-order-lift} gives a graph-lift error
\[
O(L\varepsilon^{2\alpha})
=
o(\varepsilon^{\widehat m-2}),
\]
because $2\alpha>\widehat m-2$. This proves the expansion on the normal
graph. Continuity of $\mathcal L_\infty$ and compactness of a sufficiently
small closed neighborhood give the uniform sign estimate.
\end{proof}

We now specialize the preceding result to a fixed balanced ray.

\begin{theorem}\label{II:thm:mixed-ray}
Assume $2<\widehat m<4$. Fix a balanced pair $(p,r^*)$ and set
\[
U:=U_{p,\ell,0},
\qquad
\ell_t:=L+\delta_t^*,
\qquad
E:=E_{(p,\ell,0)}.
\]
Let $P_E$ denote the Dirichlet projection onto $E$. For sufficiently large
$L$, the compression of $D^2J_{\rho^*}(U)$ to $E$ is uniformly invertible,
and
\[
P_E\nabla J_{\rho^*}(U+\bar v)=0
\]
has a unique solution in a fixed small normal ball. For every
$0<\alpha<\min\{1,\widehat m-2\}$,
\begin{equation*}
\|\bar v\|\le C_\alpha\varepsilon^\alpha,
\qquad
\left|
J_{\rho^*}(U+\bar v)-J_{\rho^*}(U)
\right|
\le C_\alpha\varepsilon^{2\alpha}.
\end{equation*}
Set $G_0(L):=U$ and $G_1(L):=U+\bar v(L)$. Then, for $j=0,1$,
\begin{equation}\label{II:eq:mixed-signed}
\frac{d}{dL}J_{\rho^*}(G_j(L))
=
\pi N\mathcal L_\infty(p,r^*)
e^{-(\widehat m-2)L}
+
o\left(e^{-(\widehat m-2)L}\right).
\end{equation}
The normal solution is $C^1$ in the finite-scale parameters.
\end{theorem}

\begin{proof}
Along the balanced ray,
\[
b=0,
\qquad
N\Omega_{i,t}=1,
\qquad
\frac{\mathfrak S_i}{\rho_i^*}
=
e^{-\fJ_i}\rf_{p,t}^i(p_t)\lambda_t^{q_i}.
\]
Theorem~\ref{II:thm:moving-mixed} therefore applies for
$ \frac{\widehat m-2}{2}
<
\alpha
<
\min\{1,\widehat m-2\} $
and gives \eqref{II:eq:mixed-signed} for $G_0$ and $G_1$. The normal
estimate for smaller $\alpha>0$ follows by lowering the exponent.
Finally, \eqref{II:eq:mixed-residual} and the bounded normal Hessian give
\[
\left|
J_{\rho^*}(U+\bar v)-J_{\rho^*}(U)
\right|
\le
C\left(
\varepsilon^\alpha\|\bar v\|+\|\bar v\|^2
\right)
\le
C_\alpha\varepsilon^{2\alpha}.
\]
The estimates are uniform on compact subsets of balanced pairs.
\end{proof}

Integrating the common-scale derivative gives the energy expansion along
the projected and reduced rays.

\begin{corollary}
Define 
\begin{equation*}
C(v,N)
:=
\frac N2\sum_i\int_{\mathbb R^2}e^{v_i}v_i
+
\pi N\sum_i\sigma_i\fJ_i
-
\sum_i\rho_i^*\log(2\pi N\sigma_i).
\end{equation*}
For $j=0,1$,
\begin{equation*}
J_{\rho^*}(G_j(L))
=
C(v,N)-\frac1N\ff(p)
-
\frac{\pi N}{\widehat m-2}
\mathcal L_\infty(p,r^*)
e^{-(\widehat m-2)L}
+
o\left(e^{-(\widehat m-2)L}\right).
\end{equation*}
\end{corollary}

\begin{proof}
Set
\[
H_{i,s}:=\log\rf_{p,s}^i(p_s),
\qquad
\delta_s:=\ell_s-L,
\qquad
\overline H_i:=\frac1N\sum_sH_{i,s}.
\]
Uniformly on a separated chart with comparable scales and $b=0$,
\begin{align}
J_{\rho^*}(U_{p,\ell,0})
={}&
C(v,N)-\frac1N\ff(p)
\nonumber\\
&-
\sum_i\rho_i^*
\log\left[
\frac1N
\sum_s
e^{H_{i,s}+q_i\delta_s-\overline H_i}
\right]
+
O(e^{-\eta L}).
\label{II:eq:coarse-energy-capture}
\end{align}
Indeed, testing the free term against its Poisson sources gives, for the
$i$-th component and the $s$-th center,
\[
\int_{\mathbb R^2}e^{v_i}v_i
+
2\pi\sigma_i
\left(
m_i\ell_s-\fJ_i
+2\pi m_i\widetilde G_s(p_s;p)
\right)
\]
before division by two, up to an exponentially small error. The discarded
tails and the variation of the smooth weights have the same order. The
core denominator expansion gives
\[
H_i(U)
=
2\pi\sigma_i e^{-\fJ_i}e^{q_iL}
\sum_s
e^{H_{i,s}+q_i\delta_s}
\left(1+O(e^{-\eta L})\right).
\]
Substitution and $m^T\sigma=4\sum_i\sigma_i$ yield
\eqref{II:eq:coarse-energy-capture}.

At a balanced pair,
$
H_{i,s}+q_i\delta_s=\kappa_i\, \text{for all }s.
$
Since $\sum_s\delta_s=0$,
$
\kappa_i
=
\frac1N\sum_s(H_{i,s}+q_i\delta_s)
=
\overline H_i.
$
Hence
$
\frac1N\sum_s
e^{H_{i,s}+q_i\delta_s-\overline H_i}
=1,$ 
and the logarithmic term in
\eqref{II:eq:coarse-energy-capture} vanishes. Therefore
\[
\lim_{L\to\infty}J_{\rho^*}(G_j(L))
=
C(v,N)-\frac1N\ff(p),
\qquad j=0,1,
\]
where Theorem~\ref{II:thm:mixed-ray} gives
\[
J_{\rho^*}(G_1(L))-J_{\rho^*}(G_0(L))
=
O(\varepsilon^{2\alpha})
=o(1).
\]
Integrating \eqref{II:eq:mixed-signed} yields
\begin{align*}
J_{\rho^*}(G_j(L))
&=
\lim_{T\to\infty}J_{\rho^*}(G_j(T))
-\int_L^\infty
\frac{d}{ds}J_{\rho^*}(G_j(s))\,ds\\
&=
C(v,N)-\frac1N\ff(p)
-\frac{\pi N}{\widehat m-2}
\mathcal L_\infty(p,r^*)
e^{-(\widehat m-2)L}
+o\left(e^{-(\widehat m-2)L}\right).
\end{align*}

If $\widehat m=4$, then
\[
\frac{1}{2\pi\sigma_i}
\int_{\mathbb R^2}e^{v_i}v_i
=
\log(2\sigma_i a^2)-2,
\qquad
\fJ_i=\log\frac{2\sigma_i}{a^2},
\qquad
\rho_i^*=2\pi N\sigma_i.
\]
Consequently,
\begin{align*}
C(v,N)
&=
\pi N\sum_i\sigma_i
\left[
\log(2\sigma_i a^2)-2
+\log\frac{2\sigma_i}{a^2}
-2\log(2\pi N\sigma_i)
\right]\\
&=
-2\pi N\sum_i\sigma_i
\bigl(1+\log(\pi N)\bigr)
=
C_4.
\end{align*}
\end{proof}

	\section{The pseudogradient and its critical points at infinity}\label{sec:5}
	\subsection{The terminal field in ray coordinates}

	We first realize the signed common-scale direction as an ambient
	descending field near each balanced ray.

	\begin{theorem}
\label{II:thm:ambient-rays}
Assume \ref{II:item:a1}--\ref{II:item:a2}. Let $(p^*,r^*)$ be a
balanced critical pair and set
\begin{equation*}
\mathcal L:=\mathcal L_\infty(p^*,r^*),\qquad
\gamma(L):=U_{p^*,\ell,0},\qquad
\ell_t=L+\delta_t^*.
\end{equation*}
Up to simultaneous permutation, there exist $L_0,c,C>0$, a tube
$\mathcal T$ containing $\gamma([L_0,\infty))$, and a locally Lipschitz
field $W$ on $\mathcal T$ such that
\begin{equation*}
\|W(u)\|\le C,\qquad
DJ_{\rho^*}(u)[W(u)]
\le-c\Psi_{\widehat m}(e^{L(u)}),\qquad
W(\gamma(L))
=-\operatorname{sign}(\mathcal L)\gamma'(L).
\end{equation*}
The tubes associated with distinct pairs are disjoint for large $L$.

Set 
\begin{equation*}
u=\gamma(L)+h,\qquad
(h,\gamma'(L))_{\mathcal H}=0,
\end{equation*}
and define
\begin{equation*}
\mathcal T_{\rm in}
:=
\left\{
u\in\mathcal T:
L>L_0+1,\quad e^{8L}\|h\|^2<1
\right\}.
\end{equation*}
If $\mathcal L<0$, then $\mathcal T_{\rm in}$ is forward invariant.
Every solution $\dot u=W(u)$ with $u(0)\in\mathcal T_{\rm in}$ satisfies
\begin{equation*}
L(u(s))=L(u(0))+s\to \infty,\qquad
\|h(s)\|\to 0,
\qquad 
J_{\rho^*}(u(s))
\to
C(v,N)-\frac1N\ff(p^*).
\end{equation*}
If $\mathcal L>0$, then
\begin{equation*}
\frac{d}{ds}L(u(s))=-1,
\end{equation*}
until the trajectory reaches the face $L=L_0+1$.
\end{theorem}

	\begin{proof}
		Fix a representative $(p^*,r^*)$ and set
		$\gamma(L):=U_{p^*,\ell,0}$, $\ell_t=L+\delta_t^*$, and
		$\mathcal L:=\mathcal L_\infty(p^*,r^*)$. Theorems~\ref{II:thm:fixed-balanced-ray}
		and~\ref{II:thm:mixed-ray} give
		\begin{equation*}
		DJ_{\rho^*}(\gamma)[\gamma']=\pi N\mathcal L\,\Psi_{\widehat m}(e^L)+o(\Psi_{\widehat m}(e^L)).
		\end{equation*}
		The positive Gram limit and the projected derivative estimates
		\cite[Lemmas~A.2--A.4]{AFHZI} imply
		$c\le\|\gamma'\|\le C$ and $\|\gamma''\|\le C$.
		By \cite[(A.4)]{AFHZI}, $\|\gamma(L)\|^2=2\pi NL\sum_i m_i^2+O(1)$.
		On the unit ball about $\gamma(L)$, \eqref{II:eq:ambient-moments} gives $\|D^2J_{\rho^*}\|\le CL$ and $\|DJ_{\rho^*}\|\le C\sqrt L$.
		
		Use coordinates $u=\gamma(L)+h$ with $(h,\gamma'(L))_{\mathcal H}=0$.
		These coordinates are single-valued for large $L$ and $\|h\|<2e^{-4L}$.
		Indeed,
		$
		\partial_L(u-\gamma(L),\gamma'(L))
		=-\|\gamma'\|^2+(h,\gamma''),
		$
		which is uniformly negative. The implicit function theorem gives local
		uniqueness. By \cite[Lemma~A.5]{AFHZI}, two representations at scales
		$L,L'\to\infty$ have the same centers, up to permutation, and
		$|L-L'|=o(1)$. Balance gives the same $r^*$, so local uniqueness applies.
		The finiteness and uniform separation of the pairs follow from
		\cite[Proposition~5.1]{AFHZI}. Since $\|\gamma(L)\|\to\infty$, the tubes
		do not accumulate at finite norm.
		
		Let $s=-\operatorname{sign}\mathcal L$. In the outer tube define
		\begin{equation}\label{II:eq:ambient-tube-field}
			S=\left(1-\frac{(h,\gamma'')}{\|\gamma'\|^2}\right)\gamma',
			\qquad DL(u)[Z]=\frac{(Z,\gamma')}{\|\gamma'\|^2-(h,\gamma'')}.
		\end{equation}
		Differentiating $(h,\gamma')=0$ gives the formula for $DL$. Define
		\begin{equation*}
			V_s=sS-6h .
		\end{equation*}
		The field is bounded and locally Lipschitz, with $DL[V_s]=s$. Moreover,
		$
		DJ_{\rho^*}(u)[V_s]
		=sDJ_{\rho^*}(\gamma)[\gamma']+O(L\|h\|).
		$
		Since $Le^{-4L}=o(\Psi_L)$, its descent is at least $c\Psi_L$.
		Because $S\perp h$ and $\gamma'\perp h$,
		\begin{equation*}
		\tfrac12\tfrac d{dt}\|h\|^2=-6\|h\|^2,\qquad
		\tfrac d{dt}(e^{8L}\|h\|^2)=(8s-12)e^{8L}\|h\|^2.
		\end{equation*}
		
		Set
		$
		\mathcal T_{\rm in}:=\{L>L_0+1:\ e^{8L}\|h\|^2<1\}.
		$
		If $s=1$, $\mathcal T_{\rm in}$ is forward invariant and $L\to\infty$.
		If $s=-1$, its lateral boundary remains invariant until $L=L_0+1$.
		On a negative tail,
		$
		|J_{\rho^*}(u)-J_{\rho^*}(\gamma)|
		\le C\sqrt L\,e^{-4L},
		$
		and the ray energy converges. In the mixed case,
		$\|\nabla J_{\rho^*}(\gamma)\|\to0$ by
		\eqref{II:eq:mixed-residual}. If $\widehat m=4$, then
		$\|\nabla J_{\rho^*}(B)\|=O(e^{-L})$ at balance. Since
		$\|B-\gamma\|=O(e^{-L})$ and $\|D^2J_{\rho^*}\|\le CL$,
		$$\|\nabla J_{\rho^*}(\gamma)\|=O(Le^{-L}).$$
        The Hessian bound gives
		the same limit on $\mathcal T_{\rm in}$. Hence these tails lie in the
		small-gradient charts, and \cite[Lemma~2.3]{AFHZI} identifies their
		limiting parameters.
	\end{proof}
	
	\subsection{The locally Lipschitz field outside the bubbling neighborhoods}
	\label{II:sec:exterior}
	Set 
	\begin{equation*}
	P(u)=\sum_i\rho_i^*\log H_i(u),\qquad
	J_t=J_{\rho^*}+tP=J_{(1-t)\rho^*},\qquad 0\le t\le\kappa .
	\end{equation*}
	Fix $0<\kappa<1$ independently of the bubble scales.
	We use Lucia's deformation method \cite{Lucia2007}, as in
	\cite[Section~2]{MBL2017}. The next lemma gives the compactness needed for
	the indefinite quadratic form in \eqref{II:def:E-functional}.
	
	\begin{lemma}
		\label{II:lem:exterior-mass-compact}
		If $0\le t_k\le\kappa$, $\sup_k P(u^k)<\infty$ and
		$\|\nabla J_{t_k}(u^k)\|_{\mathcal H}\to0$, then, after extraction,
		$t_k\to t$ and $u^k\to u$ in $\mathcal H$, with $DJ_t(u)=0$.
	\end{lemma}
	\begin{proof}
		Jensen's inequality and the fixed component means give
		$H_i(u^k)\ge\min_M h_i>0$, so the upper bound on $P$ bounds every
		$H_i$ above. Denote $z^k=A\nabla J_{t_k}(u^k)$ in the positive
		Dirichlet metric. Then $z^k\to0$ in $\mathcal H$ and, using the
		mean-zero Green function of Section~\ref{II:sec:resonant-data}, exactly
		\begin{equation*}
		u_i^k=z_i^k+G^g*\nu_i^k,\qquad
		d\nu_i^k=\sum_j a_{ij}(1-t_k)\rho_j^*
		\frac{h_j e^{u_j^k}}{H_j(u^k)}\,dV_g .
		\end{equation*}
      
By \ref{II:assump:H1},
\begin{equation*}
\nu_i^k\ge0,
\qquad
\nu_i^k(M)
=(1-t_k)\sum_j a_{ij}\rho_j^*
\le C.
\end{equation*}
After extraction,
\begin{equation*}
\nu_i^k\rightharpoonup\nu_i,
\qquad
z_i^k\to0
\quad\text{a.e. on }M.
\end{equation*}
Since $G^g(x,\cdot)$ is lower semicontinuous and bounded below,
$$G^g*\nu_i(x)
\le
\liminf_{k\to\infty}G^g*\nu_i^k(x).$$
Fatou's lemma yields
\begin{equation*}
\int_M h_i e^{G^g*\nu_i}\,dV_g
\le
\liminf_{k\to\infty}H_i(u^k)
<\infty.
\end{equation*}
If $\nu_i(\{y\})=a\ge4\pi$, then, near $y$,
\begin{equation*}
G^g*\nu_i(x)
\ge
-\frac{a}{2\pi}\log d_g(x,y)-C,
\end{equation*}
and hence
\begin{equation*}
h_i(x)e^{G^g*\nu_i(x)}
\ge
C\,d_g(x,y)^{-a/(2\pi)}
\notin L^1(M).
\end{equation*}
Therefore
$
\nu_i(\{y\})<4\pi,\, 
\text{for every }y\in M.
$
	
Since $\nu_i(\{x\})<4\pi$ for every $x\in M$, there exist
coordinate disks $\{B_r^\alpha\}_{\alpha=1}^J$ covering $M$ and
$m<4\pi$ such that  $ \nu_i(\overline{B_{2r}^\alpha})<m
\, 
\text{for all }i,\alpha.  $
By weak convergence, for all large $k$,
\begin{equation*}
\nu_i^k(B_{2r}^\alpha)\le m
\qquad
\text{for all }i,\alpha.
\end{equation*}
Choose $p>1$ such that $pm<4\pi$. On $B_r^\alpha$,
\begin{equation*}
G^g*\nu_i^k
=
G_{2r}*\bigl(\nu_i^k|_{B_{2r}^\alpha}\bigr)+O(1),
\end{equation*}
where $G_{2r}$ is the Dirichlet Green function on $B_{2r}^\alpha$.
Indeed, the two kernels have the same logarithmic singularity, while
the sources outside $B_{2r}^\alpha$ contribute a bounded term.
The Brezis--Merle estimate \cite[Theorem~1]{BM1991} gives
\begin{equation*}
\sup_k
\left\|
e^{G^g*\nu_i^k}
\right\|_{L^p(B_r^\alpha)}
<\infty.
\end{equation*}
Moreover, Moser--Trudinger and $z_i^k\to0$ in $\mathcal H$ imply
\begin{equation*}
\sup_k\|e^{z_i^k}\|_{L^q(M)}<\infty
\qquad
\text{for every }q<\infty.
\end{equation*}
Hence, for every $1<p'<p$, Hölder's inequality and
$H_j(u^k)\ge\min_M h_j$ yield
\begin{equation*}
\sup_k
\left\|
\frac{d\nu_i^k}{dV_g}
\right\|_{L^{p'}(M)}
<\infty.
\end{equation*}
Therefore
\begin{equation*}
\|G^g*\nu_i^k\|_{W^{2,p'}(M)}\le C,
\qquad
G^g*\nu_i^k\to  G^g*\nu_i
\quad\text{in }H^1(M).
\end{equation*}
Since $z^k\to0$ in $\mathcal H$,
\begin{equation*}
u^k\to  u
\quad\text{in }\mathcal H.
\end{equation*}
After also taking $t_k\to t$, continuity of the first variation gives
$ DJ_t(u)=0.$
\end{proof}
	
	The compactness lemma yields a joint descending field for the endpoint
	functionals.

	\begin{proposition}\label{II:prop:exterior-gluing}
		Set $d_\kappa=\min_{0\le t\le\kappa}\|\nabla J_t\|$ and $E_\kappa=\{d_\kappa=0\}$. On its open complement there is a locally Lipschitz field $W_{\rm ext}$ of norm at most one such that
		\begin{equation}\label{II:eq:BL-exterior}
			-DJ_{\rho^*}[W_{\rm ext}]\ge d_\kappa/2,\qquad -DJ_\kappa[W_{\rm ext}]\ge d_\kappa/2.
		\end{equation}
	Every trajectory starting outside an open neighborhood of $E_\kappa$ either
	enters that neighborhood or reaches any prescribed sublevel of
	$J_{\rho^*}$ in finite time.
	\end{proposition}
	\begin{proof}
Fix $u$ with $d_\kappa(u)>0$ and set
\begin{equation*}
g(u)
:=
\operatorname*{argmin}_{\xi\in
[\nabla J_{\rho^*}(u),\,\nabla J_\kappa(u)]}
\|\xi\|.
\end{equation*}
Then $\|g(u)\|=d_\kappa(u)$, and the projection identity gives
\begin{equation*}
(g(u),v-g(u))\ge0
\qquad
\text{for every }
v\in[\nabla J_{\rho^*}(u),\nabla J_\kappa(u)].
\end{equation*}
Hence, for
\begin{equation*}
W_u:=-\frac{g(u)}{\|g(u)\|},
\end{equation*}
one has
\begin{equation*}
-DJ_{\rho^*}(u)[W_u],
\;
-DJ_\kappa(u)[W_u]
\ge d_\kappa(u).
\end{equation*}
By continuity, these inequalities hold with $d_\kappa/2$ on a
neighborhood of $u$. A locally finite Lipschitz partition of unity,
as in \cite[Proposition~2.4]{HL2009}, gives $W_{\rm ext}$ satisfying
\eqref{II:eq:BL-exterior}.

Let $u(s)$ solve $\dot u=W_{\rm ext}(u)$. If
$J_{\rho^*}(u(s))\ge a$, then
\begin{equation*}
J_\kappa(u(s))\le J_\kappa(u(0)),
\qquad
P(u(s))
=
\frac{J_\kappa(u(s))-J_{\rho^*}(u(s))}{\kappa}
\le
\frac{J_\kappa(u(0))-a}{\kappa}.
\end{equation*}
For every open neighborhood $\mathcal U\supset E_\kappa$, set
$C_a:=\frac{J_\kappa(u(0))-a}{\kappa}.$
Then
\begin{equation*}
\delta
:=
\inf\left\{
d_\kappa(v):
v\notin\mathcal U,\ P(v)\le C_a
\right\}
>0.
\end{equation*}
Indeed, otherwise there exist $v_k\notin\mathcal U$ and
$t_k\in[0,\kappa]$ such that
\begin{equation*}
P(v_k)\le C_a,
\qquad
\|\nabla J_{t_k}(v_k)\|
=d_\kappa(v_k)\to0.
\end{equation*}
Lemma~\ref{II:lem:exterior-mass-compact} gives
$v_k\to v\in E_\kappa\subset\mathcal U$, a contradiction.

Therefore, while $u(s)\notin\mathcal U$ and
$J_{\rho^*}(u(s))\ge a$,
\begin{equation*}
-\frac{d}{ds}J_{\rho^*}(u(s))
\ge\frac{\delta}{2}.
\end{equation*}
Thus, the trajectory enters $\mathcal U$ or reaches
$J_{\rho^*}\le a$ within time
$T\le
{2\bigl(J_{\rho^*}(u(0))-a\bigr)}/{\delta}.$
The bound $\|W_{\rm ext}\|\le1$ gives a limit at every finite maximal
time, and local existence extends the trajectory unless it has entered
$\mathcal U$.
	\end{proof}
	\subsection{Estimates for crossing a modulation neighborhood}
	
	Fix an ordered balanced Morse pair $(p^*,r^*)$. In the original
	minimizing chart denote $u=U_z+w$, $w\in E_z$, and set
	\begin{equation*}
	x=|p-p^*|,\qquad X=|r-r^*|+L|b|+\|w\|,\qquad \varepsilon=e^{-L}.
	\end{equation*}
	Use a smooth position chart whose Euclidean norm is equivalent to the
	product surface distance. All estimates below hold for $x+X\le h$, with
	$h>0$ fixed and small and $L$ large, in both mass regimes. Here $t$ denotes
	the parameter in $J_t$, while $s$ denotes a bubble index.
	Set $y_{i,s}=1-(1-t)N\Omega_{i,s}$. Thus the mean of $y_{i,s}$ over
	$s$ is $t$. Define $H^\alpha$, $r_\infty$ and $c_\infty$ as in
	Theorems~\ref{II:thm:bounded-complete-mass-inverse}
	and~\ref{II:thm:finite-scale-covector-center}, with $q_i=m_i-2$.
	In particular $r_\infty(p^*)=r^*$ and $c_\infty(p^*)=0$.
	
	Here $L$ is the modulation scale; when needed, the ray coordinate in
	\eqref{II:eq:ambient-tube-field} is denoted by $L_{\rm ray}$. The
	small-gradient condition defining $V^\eta_{\rho^*}$ is required only at
	entry. The Gram estimates and their first derivatives
	\cite[Lemma~2.3 and Appendix~A]{AFHZI} extend the modulation to the box
	$x+X\le h$, independently of $\|DJ_{\rho^*}(u)\|$. Moreover,
	$L|b|\le h$ implies
	$
	\|U_z-U_{p,\ell,0}\|=O(hL^{-1/2}).
	$
	The separation result \cite[Lemma~A.5]{AFHZI} and the local inverse give
	uniqueness on each ordered sheet. Thus a continuous curve can be followed
	to every parameter face, including the gradient-cutoff face.

	The next lemma converts a small tilted gradient into estimates for all
	modulation variables.
	
	\begin{lemma}
		\label{II:lem:collar-modulation}
		For sufficiently small $h,\kappa>0$ there are $C,d_0,L_*,\eta_1>0$ such that whenever $x+X\le h$, $L\ge L_*$ and $d=d_\kappa(u)\le d_0$, every minimizing $t$ satisfies
		\begin{equation}\label{II:eq:collar-tilt}
			\begin{aligned}
				\|w\|+|y|&\le C\bigl(d+|b|+e^{-\eta_1 L}\bigr),\\
				0\le t&\le C\bigl(d+|b|^2+e^{-\eta_1 L}\bigr).
			\end{aligned}
		\end{equation}
		For a sequence in this chart with $L\to\infty$, $p\to p_0$ and $d_\kappa=o(L^{-2})$,
		\begin{equation*}
		\|w\|\to0,\quad Lb^0\to0,\quad
		r\to r_\infty(p_0),\quad c^\alpha\to c_\infty^\alpha(p_0),
		\quad \theta\to0 .
		\end{equation*}
		If $d_\kappa=o(e^{-L})$, then $\nabla\ff(p_0)=0$.
	\end{lemma}
\begin{proof}
Choose $t\in[0,\kappa]$ such that $d:=d_\kappa(u)=\|\nabla J_t(u)\|.$
Throughout the proof, $E_L$ denotes a quantity bounded by
$Ce^{-\eta L}$ for some $\eta>0$; decreasing $\eta$ absorbs any fixed
power of $L$. Set 
\begin{equation*}
Y=\lambda_s y_{p_s}(x),\qquad
V_{i,s}\,dV_g=\chi_s e^{v_i(Y)}\,dY,\qquad
\psi_i:=\sum_j a^{ij}v_j.
\end{equation*}
The normal constraint gives
\begin{equation*}
\int_M V_{i,s}w_i\,dV_g=0.
\end{equation*}

On the $s$-th core, \cite[Lemma~A.1]{AFHZI} gives
\begin{equation*}
h_i e^{U_{z,i}}\,dV_g
=
\lambda_s^{q_i+\xi_{i,s}}
e^{v_i(Y)+b_{i,s}\psi_i(Y)}
a_{i,s}(Y/\lambda_s)
\bigl(1+O(E_L)\bigr)\,dY.
\end{equation*}
The functions $a_{i,s}^{\pm1}$ have uniformly bounded $C^1$ norms.
After decreasing $\eta$,
\begin{equation*}
\int_{D_s}h_i e^{U_{z,i}}\,dV_g
=
\mathfrak M_{i,s}\bigl(1+O(E_L)\bigr),
\end{equation*}
with the same relative error after one derivative in $\ell_s$. Moreover,
on anchored single and double tests,
\begin{equation*}
\frac{e^{v_i+b_{i,s}\psi_i}}
     {\int_{\mathbb R^2}e^{v_i+b_{i,s}\psi_i}}
=
\frac{e^{v_i}}{2\pi\sigma_i}
+
O\!\left(
|b|(1+\log(2+|Y|))e^{v_i}
\right).
\end{equation*}
Thus, with $D_i^t:=(1-t)D_i$,
\begin{equation*}
D_i^t(U_z)
=
\sum_s(1-y_{i,s})V_{i,s}
+
O_{\rm test}(|b|+E_L).
\end{equation*}
Consequently, the normal Hessian differs from the model Hessian by
\begin{equation*}
O(|y|+|b|+E_L).
\end{equation*}
Since $|y|\le C(h+\kappa)+o(1)$, the constrained inverse remains
uniformly bounded for small $h,\kappa$. The normal equation therefore
gives
\begin{equation}\label{II:eq:collar-normal}
\|w\|\le C(d+|b|+E_L).
\end{equation}

We next test $DJ_t(U_z+w)$ against $e_iPW_{i,s}/\ell_s$. On $D_s$,
\begin{equation*}
\frac{PW_{i,s}}{\ell_s}
=
\sigma_i
+
O\!\left(\frac{1+\log(2+|Y|)}{L}\right),
\end{equation*}
while it is $O(L^{-1})$ on the other cores. Hence
\begin{align*}
DJ_t(U_z+w)\left[\frac{e_iPW_{i,s}}{\ell_s}\right]
={}&
2\pi\sigma_i
\bigl((K\xi_s)_i+\sigma_i y_{i,s}\bigr)\\
&+
O\!\left(
E_L+\frac{|b|+|y|+\|w\|}{L}
+|b|\|w\|+\|w\|^2
\right).
\end{align*}
Since the test has a uniformly bounded norm and
$\|\nabla J_t(u)\|=d$, we obtain
\begin{equation}\label{II:eq:collar-amp}
|K\xi_s+D_\sigma y_s|
\le
C\left(
d+E_L+\frac{|b|+|y|+\|w\|}{L}
+|b|\|w\|+\|w\|^2
\right).
\end{equation}

For the scale test, set
\begin{equation*}
F_{i,s}:=v_i+b_{i,s}\psi_i,\qquad
\mathfrak a_s:=m+\xi_s,\qquad
S_s:=\partial_{\ell_s}U_z.
\end{equation*}
By \cite[Lemma~A.1]{AFHZI},
\begin{equation*}
S_s
=
\mathfrak a_s+Y\cdot\nabla F_s+O(E_L)
\quad\text{in  }D_s,
\qquad
S_s=O(E_L)
\quad\text{in  }D_s^c.
\end{equation*}
The radial identity in Theorem~\ref{II:thm:moving-mixed} gives
\begin{align*}
DJ_t(U_z)[S_s]
={}&
\pi\mathfrak a_s^TK\mathfrak a_s
-\sum_i(1-t)\rho_i^*\Omega_{i,s}
   (q_i+\xi_{i,s})
+O(E_L)\\
={}&
2\pi q^TD_\sigma y_s
+O(|b|^2+|b||y|+E_L),
\end{align*}
where resonance cancels the constant term and $Km=\sigma$ cancels
the term linear in $b$.

The normal transfer contributes
\begin{equation*}
O\!\left(
(|y|+|b|+E_L)\|w\|+\|w\|^2
\right).
\end{equation*}
Indeed, the anchor removes the constant $q_i$, while the Jacobi terms
cancel. Therefore
\begin{equation}\label{II:eq:collar-scales}
|q^TD_\sigma y_s|
\le
C\left(
d+E_L+|b|^2+|b||y|
+(|y|+|b|)\|w\|+\|w\|^2
\right).
\end{equation}

Combining \eqref{II:eq:collar-normal} and
\eqref{II:eq:collar-amp}, and absorbing $C|y|/L$, gives
\begin{equation*}
\|w\|+|y|
\le
C(d+|b|+E_L).
\end{equation*}
Since $\frac1N\sum_s y_s=t\mathbf1,
\,
q^TD_\sigma\mathbf1>0,$
averaging \eqref{II:eq:collar-scales} over $s$ yields
\begin{equation*}
0\le t
\le
C(d+|b|^2+E_L).
\end{equation*}
This proves \eqref{II:eq:collar-tilt}.

Now suppose $L\to\infty$, $p\to p_0$, and
$d_\kappa=o(L^{-2})$. Since $|b|\le h/L$,
\begin{equation*}
\|w\|+|y|=O(L^{-1}),
\qquad
t=O(L^{-2}).
\end{equation*}
Equations \eqref{II:eq:collar-amp} and
\eqref{II:eq:collar-scales} give
\begin{equation*}
KD_\sigma b_s+D_\sigma y_s=O(L^{-2}),
\qquad
q^TD_\sigma y_s=O(L^{-2}).
\end{equation*}
Taking the mean and the relative modes yields
\begin{equation*}
Lb^0=O(L^{-1}),
\qquad
q^TKD_\sigma c^\alpha=O(L^{-1}).
\end{equation*}
The logarithmic mass equations give
\begin{equation*}
H^\alpha(p)+qr_\alpha+D_\sigma c^\alpha=o(1).
\end{equation*}
Since $q^TKq>0$, the last two relations imply
\begin{equation*}
r\to r_\infty(p_0),
\qquad
c^\alpha\to c_\infty^\alpha(p_0).
\end{equation*}
Moreover, $y=t\mathbf1+(1-t)\theta$, so
\begin{equation*}
\theta\to0,\qquad Lb^0\to0,\qquad \|w\|\to0.
\end{equation*}

It remains to identify the limiting centers. Consider a sequence such that
\begin{equation*}
d_\kappa(u^k)=o(e^{-L_k}),
\qquad
L_k\to\infty,
\qquad
p^k\to p_0.
\end{equation*}
Let $t_k$ minimize $d_\kappa(u^k)$ and set
\begin{equation*}
z^k:=A\nabla J_{t_k}(u^k),
\qquad
d\nu_i^k:=D_i^{t_k}(u^k)\,dV_g.
\end{equation*}
Then, $\nu_i^k
\rightharpoonup
2\pi\sigma_i\sum_{s=1}^N\delta_{p_{0,s}}.$
The core comparison, the anchored exponential estimate, and Jensen's
inequality give
\begin{equation*}
\left\|
\frac{d\nu_i^k}{dV_g}
\right\|_{L^2(M)}
\le Ce^{L_k}.
\end{equation*}

Fix a center $p_{0,s}$ and choose a smooth vector field $\mathcal V$
supported near $p_{0,s}$ and constant in centered isothermal coordinates
near that point. From
\begin{equation*}
u_i^k
=
z_i^k+\sum_j a_{ij}G^g*\nu_j^k
\end{equation*}
and
$\int_M\operatorname{div}_g(D_i^{t_k}(u^k)\mathcal V)\,dV_g=0$,
symmetrization gives
\begin{align*}
0={}&
\sum_i\int_M
\left(
\operatorname{div}_g\mathcal V
+\mathcal V(\log h_i)
+\mathcal V(z_i^k)
\right)d\nu_i^k\\
&+
\frac12\sum_{i,j}a_{ij}
\iint_{M^2}
(\mathcal V_x+\mathcal V_y)G^g(x,y)\,
d\nu_i^k(x)d\nu_j^k(y).
\end{align*}
Since
\begin{equation*}
\left|
\int_M\mathcal V(z_i^k)\,d\nu_i^k
\right|
\le
Ce^{L_k}d_\kappa(u^k)
=o(1),
\end{equation*}
this term vanishes in the limit.

In local coordinates,
\begin{equation*}
G^g(x,y)
=
-\frac1{2\pi}\log|x-y|+H(x,y),
\end{equation*}
where $H$ is smooth and symmetric. The symmetrized logarithmic kernel is
bounded because $\mathcal V(x)-\mathcal V(y)=O(|x-y|).$
It is continuous at every atom of the limiting product measures.
Passing to the limit and using $A\sigma=m$ gives
\begin{equation*}
\sum_i\sigma_i
\mathcal V\!\left(\log\rf_{p_0,s}^i\right)(p_{0,s})
=0.
\end{equation*}
Since $\mathcal V(p_{0,s})$ is arbitrary,
\begin{equation*}
\sum_i\sigma_i
\nabla_x\log\rf_{p_0,s}^i(p_{0,s})
=0
\qquad
\text{for }s=1,\ldots,N.
\end{equation*}
The identity for $\nabla\ff$ in
Section~\ref{II:sec:resonant-data}  yields  $\nabla\ff(p_0)=0.$
\end{proof}

	The preceding estimates give quantitative costs for leaving a fixed
	modulation neighborhood.

	\begin{corollary}\label{II:cor:crossing-cost}
		For each sufficiently small $b_0>0$, one may choose $a_0>0$ arbitrarily
		small such that, in the box $x\le2a_0$, $X\le2b_0$ and for large $L$,
		\begin{equation*}
		d_\kappa\ge c_{a_0,b_0}e^{-L}\quad(a_0\le x\le2a_0),
		\qquad
		d_\kappa\ge c_{a_0,b_0}L^{-2}\quad(b_0\le X\le2b_0).
		\end{equation*}
		Let a field on the box satisfy $\|W\|\le B$ and $-DJ_{\rho^*}[W]\ge\mu d_\kappa$ with fixed $B,\mu>0$. A crossing of the first strip along $W$ costs at least $c_{a_0,b_0}\mu/B$ in energy. A crossing of the second with $L\le H$ costs at least $c_{a_0,b_0}(\mu/B)H^{-5/2}$.
	\end{corollary}
	\begin{proof}
The balance condition  gives
$r_\infty(p^*)=r^*,
\,
c_\infty^\alpha(p^*)=0.$ By continuity, choose $a_0>0$ such that $p^*$ is the only critical point
of $\ff$ in $\{|p-p^*|\le2a_0\}$ and
\[
\sup_{|p-p^*|\le2a_0}
\left(
|r_\infty(p)-r^*|
+\max_{i,t}\left|
\sum_\alpha e_t^\alpha c_{\infty,i}^\alpha(p)
\right|
\right)
<\frac{b_0}{2}.
\]

If the first lower bound failed, there would exist a sequence with
\[
a_0\le x_k\le2a_0,
\qquad
L_k\to\infty,
\qquad
e^{L_k}d_\kappa(u_k)\to0.
\]
After extraction, $p_k\to p_0$. Lemma~\ref{II:lem:collar-modulation}
would give $\nabla\ff(p_0)=0$, contrary to the choice of $a_0$.
Similarly, failure of the second bound would give
\[
b_0\le X_k\le2b_0,
\qquad
L_k\to\infty,
\qquad
L_k^2d_\kappa(u_k)\to0.
\]
The same lemma then yields
\[
X_k\to
|r_\infty(p_0)-r^*|
+\max_{i,t}\left|
\sum_\alpha e_t^\alpha c_{\infty,i}^\alpha(p_0)
\right|
<\frac{b_0}{2},
\]
a contradiction.

Let $V$ be the velocity along a crossing. Differentiating the modulation
orthogonality conditions and using \cite[Appendix~A]{AFHZI}, we obtain
\begin{equation*}
|\dot p|\le Ce^{-L}\|V\|,
\qquad
|\dot L|+|\dot r|+\sqrt L\,|\dot b|+\|\dot w\|
\le C\|V\|.
\end{equation*}
Consequently,
\begin{equation*}
|\dot x|\le Ce^{-L}\|V\|,
\qquad
|\dot X|
\le |\dot r|+|\dot L|\,|b|+L|\dot b|+\|\dot w\|
\le C\sqrt L\,\|V\|.
\end{equation*}

Suppose that $\|V\|\le B$ and
$-DJ_{\rho^*}[V]\ge\mu d_\kappa$. On the position strip,
\begin{equation*}
-\frac{d}{ds}J_{\rho^*}(u(s))
\ge c\mu e^{-L}
\ge \frac{c\mu}{B}|\dot x|.
\end{equation*}
On the transverse strip, if $L\le H$, then
\begin{equation*}
-\frac{d}{ds}J_{\rho^*}(u(s))
\ge c\mu L^{-2}
\ge c\mu H^{-5/2}\sqrt L
\ge \frac{c\mu}{B}H^{-5/2}|\dot X|.
\end{equation*}
Integration across strips of widths $a_0$ and $b_0$, respectively,
gives the two stated energy bounds.
\end{proof}
	
	\subsection{Pseudogradient directions in the bubbling coordinates}
	
	The following parameter fields are defined on a compact separated location chart with bounded scale ratios; see \cite[Proposition~8]{MBL2017} for the scalar construction.

To control the amplitudes and the mass defects on the full bubbling chart,
we use localized amplitude directions. In the original coordinates
$u=U_z+w$, assume that $\sqrt L\,|b|+\|w\|$ is small; no bound on
$L|b|$ is imposed. Set
\[
D_s:=U_{r_0/2}(p_s),\qquad
\omega_{i,s}:=\mu_i(u)(D_s),\qquad
T_{i,s}:=\frac{e_iPW_{i,s}}{\sigma_i\ell_s},
\]
where $e_i$ is the $i$th standard vector of $\mathbb R^n$. Uniformly for
$0\le t\le\kappa$,
\begin{align}
DJ_t(u)[T_{i,s}]
&=
2\pi\left[
(KD_\sigma b_s)_i
+\sigma_i\bigl(1-N(1-t)\omega_{i,s}\bigr)
\right]
+O\left(\frac{1+\|w\|}{L}\right),
\label{II:eq:disk-amplitude-test}\\
|\omega_{i,s}-\Omega_{i,s}|
&\le
C\left(e^{-\eta L}+|b|\,\|w\|+\|w\|^2\right)
\label{II:eq:disk-fraction-test}
\end{align}
for some $\eta>0$.

We first prove \eqref{II:eq:disk-fraction-test}. By
\cite[Lemma~A.1]{AFHZI}, after division by $\mathfrak M_{i,s}$, the
density on $D_s$ at $w=0$ has radial part
\[
\exp\left(v_i+b_{i,s}\sum_j a^{ij}v_j\right)
\]
and a smooth factor
\[
1+O\left(e^{-\eta L}(1+|Y|)^\eta\right).
\]
After decreasing $\eta$, the weighted radial tails and all fixed
logarithmic moments are integrable. The anchor
$\int_MV_{i,s}w_i\,dV_g=0$ cancels the leading term linear in $w_i$.
The change of the radial density contributes $O(|b|\,\|w\|)$, while
Lemma~\ref{II:lem:normal-model} gives an $O(\|w\|^2)$ bound for the
exponential remainder. Hence
\begin{equation*}
\int_{D_s}h_ie^{U_{z,i}+w_i}\,dV_g
=
\mathfrak M_{i,s}
\left[
1+O\left(e^{-\eta L}+|b|\,\|w\|+\|w\|^2\right)
\right].
\end{equation*}
Moreover,
\begin{equation*}
\int_{M\setminus\bigcup_lD_l}h_ie^{U_{z,i}+w_i}\,dV_g
\le Ce^{-\eta L}\mathfrak S_i.
\end{equation*}
Division by the total mass gives \eqref{II:eq:disk-fraction-test}.

On $D_s$,
\[
(T_{i,s})_i
=
1+O\left(\frac{1+\log(2+|Y|)}{L}\right),
\]
whereas $(T_{i,s})_i=O(L^{-1})$ away from $D_s$. The preceding moment
bounds therefore give
\begin{equation*}
(1-t)\rho_i^*
\int_M(T_{i,s})_i\,d\mu_i(u)
=
(1-t)\rho_i^*\omega_{i,s}
+O\left(\frac{1+\|w\|}{L}\right).
\end{equation*}
The free sources in component $i$ are
\[
\sum_{j,l}\bigl(\delta_{ij}+a^{ij}b_{j,l}\bigr)V_{j,l}
\]
modulo constants. 
The free pairing is
\begin{equation*}
\sum_{j=1}^n a^{ij}
\int_M
\nabla_g(U_{z,j}+w_j)\cdot
\nabla_g\left(\frac{PW_{i,s}}{\sigma_i\ell_s}\right)\,dV_g
=
2\pi\left[\sigma_i+(KD_\sigma b_s)_i\right]
+O\left(\frac{1+\|w\|}{L}\right).
\end{equation*}

Since $\rho_i^*=2\pi N\sigma_i$, subtraction yields
\eqref{II:eq:disk-amplitude-test}.

	The common amplitude is controlled outside the bounded logarithmic window.

	\begin{proposition}\label{II:prop:large-common}
		Under the preceding assumptions, put $g_i=(2\pi)^{-1}\sum_s DJ_{\rho^*}(u)[T_{i,s}]$. Then
		\begin{equation*}
			g=NK D_\sigma b^0+O((1+\|w\|)/L).
		\end{equation*}
		For $|b^0|\ge M/L$ with $M$ sufficiently large there is a bounded locally Lipschitz field keeping $p,\ell$ fixed with $-DJ_{\rho^*}[W]\ge c|b^0|$.
	\end{proposition}
	\begin{proof}
Summing \eqref{II:eq:disk-amplitude-test} over $s$ at $t=0$ gives
\begin{equation*}
\frac1{2\pi}\sum_sDJ_{\rho^*}(u)[T_{i,s}]
=
N(KD_\sigma b^0)_i
+O\left(\frac{1+\|w\|}{L}\right).
\end{equation*}
For $Y_{i,s}:=\frac1{\sigma_i\ell_s}\partial_{b_{i,s}},$
the corresponding horizontal lift is
$\widetilde T_{i,s}
:=
T_{i,s}+(D_{Y_{i,s}}P_z)w.$
The normalized frame estimates in
\cite[Lemmas~A.2--A.4]{AFHZI} give
\begin{equation*}
DJ_t(u)[\widetilde T_{i,s}]
=
DJ_t(u)[T_{i,s}]
+O\left(\frac{\|w\|}{L}\right),
\qquad 0\le t\le\kappa.
\end{equation*}
Indeed, the geometric and amplitude coefficients of
$(D_{Y_{i,s}}P_z)w$ are
\[
O\left(\frac{\|w\|}{L}\right),
\qquad
O\left(\frac{\|w\|}{L^{3/2}}\right),
\]
respectively.

For $b^0\ne0$, set $W
:=
-\frac1{2\pi}
\sum_{i,s}
\frac{(KD_\sigma b^0)_i}{|KD_\sigma b^0|}
\,\widetilde T_{i,s}.$
The normalized frame bounds imply $\|W\|\le C$, and
\begin{align*}
-DJ_{\rho^*}(u)[W]
&\ge
N|KD_\sigma b^0|
-C\frac{1+\|w\|}{L}\ge c|b^0|
\end{align*}
whenever $|b^0|\ge M/L$ and $M$ is sufficiently large.
	\end{proof}
	
	We next construct the normal descent field in both mass regimes.

	\begin{proposition}\label{II:prop:normal-field}
		Use $B_z=U_z^\sharp$, $\Pi_z=P_z^\sharp$ when $\widehat m=4$, and $B_z=U_z$, $\Pi_z=P_z$ when $2<\widehat m<4$. On a cell with $L|b|\le C$ and $\Delta_\Omega\le\delta_*$ sufficiently small, there is a normal graph $G_z=B_z+\bar v_z$ for large $L$. In both regimes, for $h\in\operatorname{ran}\Pi_z$ of sufficiently small norm,
		\begin{equation*}
		c\|h\|\le\|\Pi_z\nabla J_{\rho^*}(G_z+h)\|\le C\|h\|.
		\end{equation*}
		The locally Lipschitz vector field
		\begin{equation*}
		W_\nu=-\frac{\Pi_z\nabla J_{\rho^*}(G_z+h)}
		{\max\{\Psi_L,\|\Pi_z\nabla J_{\rho^*}(G_z+h)\|\}}
		\end{equation*}
		has norm at most one and keeps every modulation parameter fixed. For $\|h\|\ge\Psi_L$ its descent is at least $c'\|h\|$.
	\end{proposition}
	\begin{proof}
For $\widehat m=4$, the existence of the normal graph follows from
Theorem~\ref{II:thm:m4-verified-normal-inverse}. Suppose
$2<\widehat m<4$. On the $t$th core, the leading potential is
$N\Omega_{i,t}V_{i,t}$, and
\[
\int_MV_{i,t}w_i\,dV_g=0
\qquad
(w\in E_z).
\]
The radial shape, spatial variation, and free amplitude terms therefore
give
\begin{equation*}
\left\|P_z\nabla J_{\rho^*}(U_z)\right\|
\le C\left(|b|+e^{-\eta L}\right)
\end{equation*}
for some $\eta>0$. By Lemma~\ref{II:lem:normal-model}, the compressed
Hessian is uniformly invertible when $\Delta_\Omega$ is sufficiently
small, and it is uniformly Lipschitz on a fixed normal ball. The
contraction argument then gives a unique normal solution $\bar v_z$ with
\begin{equation*}
P_z\nabla J_{\rho^*}(U_z+\bar v_z)=0,
\qquad
\|\bar v_z\|
\le C\left(|b|+e^{-\eta L}\right).
\end{equation*}

Thus, in both regimes, $\Pi_z\nabla J_{\rho^*}(G_z)=0.$
For $h\in\operatorname{ran}\Pi_z$, Taylor expansion in the fixed normal
fiber gives
\begin{equation*}
\Pi_z\nabla J_{\rho^*}(G_z+h)
=
\Pi_zD^2J_{\rho^*}(G_z)h
+O(\|h\|^2).
\end{equation*}
The operator $\Pi_zD^2J_{\rho^*}(G_z)
\big|_{\operatorname{ran}\Pi_z}$
has uniform upper and inverse bounds. Hence, for $\|h\|$ sufficiently
small,
\begin{equation*}
c\|h\|
\le
\left\|\Pi_z\nabla J_{\rho^*}(G_z+h)\right\|
\le
C\|h\|.
\end{equation*}

Set
   $R:=\left\|\Pi_z\nabla J_{\rho^*}(G_z+h)\right\|.$
For the normal field $W_\nu$,
\begin{equation*}
-DJ_{\rho^*}(G_z+h)[W_\nu]
=
\frac{R^2}{\max\{\Psi_L,R\}}.
\end{equation*}
If $\|h\|\ge\Psi_L$, the two-sided estimate implies
\[
\frac{R^2}{\max\{\Psi_L,R\}}
\ge c'\|h\|.
\]
Finally, $\Psi_L>0$ at every finite scale, and the normal projection and
gradient depend smoothly on the finite-scale parameters. Thus $W_\nu$
is locally Lipschitz.
\end{proof}

	At the borderline exponent, the tangent pairings determine the transverse
	and common-scale fields.

	\begin{proposition}\label{II:prop:slice-field}
		Suppose $\widehat m=4$. On the corrected graph of Theorem~\ref{II:thm:m4-verified-normal-inverse} put $f=J_{\rho^*}\circ G$. There is a reduced transverse field with $\dot L=0$ and bounded normalized parameter velocities such that
		\begin{equation}\label{II:eq:pg-transverse}
			\begin{gathered}
				-df[W_T]\ge c A_0-C \mathcal E_0,\\
				A_0=|\mathcal C|+\sum_t\lambda_t^{-1}|\nabla_{p_t}\ff|,\qquad
				\mathcal E_0=(\delta+\varepsilon)/L+\delta^2+\varepsilon^2 L .
			\end{gathered}
		\end{equation}
		where $\delta=|b|+\Delta_\Omega$ and $\mathcal C$ denotes the covectors in \eqref{II:eq:principal-centered-covectors}. In particular $-df[W_T]\ge c' A_0$ when $A_0\ge M \mathcal E_0$ for a sufficiently large fixed $M$. In that sector one can instead choose a field with the same lower bound up to a fixed factor and $\dot\ell_t\le-1$ for all $t$. On the box $\mathcal D_L$ in Theorem~\ref{II:thm:slice-signed}, if $\mathcal L_\infty(p^*,r^*)\ne0$, there is a reduced field with
		\begin{equation*}
		-df[W]\ge c_{p^*}(|a|+\Psi_L),\qquad
		\tfrac12\le -\operatorname{sign}(\mathcal L_\infty)\,\dot\ell_t\le\tfrac32\quad\text{for every }t .
		\end{equation*}
		Here $a$ is the actual transverse derivative vector specified below. These are parameter fields on the graph. In the signed box, their images under $DG$ have uniformly bounded Dirichlet norm.
	\end{proposition}
\begin{proof}
Let $(Y_j)_j$ be the transverse parameter directions obtained from
\begin{equation*}
\frac{1}{2\pi\sigma_i\ell_t}\partial_{b_{i,t}},
\qquad
\frac{1}{2\pi}\sum_t e_t^\alpha\partial_{\ell_t},
\qquad
\lambda_t^{-1}\partial_{p_{t,k}},
\end{equation*}
where the amplitude directions are decomposed into their mean and
zero-sum modes. Set
\[
a_j:=df[Y_j],
\qquad
a:=(a_j)_j.
\]
Theorem~\ref{II:thm:m4-corrected-leading-pairings} gives, in these
coordinates,
\begin{equation*}
a
=
\left(
\mathcal C,
\left(-\frac{1}{N\lambda_t}\nabla_{p_t}\ff\right)_t
\right)
+O(\mathcal E_0)
\end{equation*}
before differentiation along the normal graph.

We estimate the graph correction using
Lemma~\ref{II:lem:ABL-first-order-lift}. With $d_a$ denoting the
pairings with the normalized tangent frame,
\begin{align*}
|d_{\rm amp}|
&\le
C\left(\sqrt L\,\delta+\varepsilon+\varepsilon^2L^{3/2}\right),\qquad
|d_{\rm geo}|
\le
C\left(\delta+\varepsilon+\varepsilon^2L\right).
\end{align*}
Moreover, \cite[Lemma~A.2]{AFHZI} and
Lemma~\ref{II:lem:amplitude-projected-size} give
\begin{align*}
\|D_{Y_{\rm amp}}T_a\|
&\le CL^{-1},\qquad
\|D_{Y_{\rm geo}}T_a^{\rm amp}\|
\le CL^{-1/2},
\qquad
\|D_{Y_{\rm geo}}T_a^{\rm geo}\|
\le C,
\end{align*}
and
\begin{equation*}
\|\mathbb G_{\rm amp,geo}\|
+
\|(\mathbb G^{-1})_{\rm amp,geo}\|
\le CL^{-1/2}.
\end{equation*}
Since
$
L|b|\le C,
\,
\|\bar v\|
\le C\left(|b|+\varepsilon^2L^{3/2}\right),
$
Lemma~\ref{II:lem:ABL-first-order-lift} yields
\begin{equation*}
\left|
df[Y_j]-DJ_{\rho^*}(G_z)[DB_z(Y_j)]
\right|
\le C\mathcal E_0.
\end{equation*}
Consequently,
\begin{equation*}
|a|
\ge
c\left(
|\mathcal C|
+\sum_t\lambda_t^{-1}|\nabla_{p_t}\ff|
\right)
-C\mathcal E_0
=
cA_0-C\mathcal E_0.
\end{equation*}

Define
\begin{equation*}
W_T
:=
-\sum_j
\frac{a_j}{\max\{|a|,\Psi_L\}}Y_j.
\end{equation*}
Then
\begin{equation*}
-df[W_T]
=
\frac{|a|^2}{\max\{|a|,\Psi_L\}}
\ge |a|-\Psi_L.
\end{equation*}
Since $\Psi_L\le\mathcal E_0$,
\begin{equation*}
-df[W_T]
\ge cA_0-C\mathcal E_0,
\end{equation*}
which proves \eqref{II:eq:pg-transverse}. The coefficients of $W_T$
have absolute value at most one. The finite-scale $C^2$ regularity of
$f$ from Theorem~\ref{II:thm:m4-verified-normal-inverse} therefore
implies that $W_T$ is locally Lipschitz.

The scale part of $W_T$ has zero mean. For some fixed $B_0>0$,
\begin{equation*}
\sum_tW_T\ell_t=0,
\qquad
|W_T\ell_t|\le B_0.
\end{equation*}
Summing the scale pairings in
Theorem~\ref{II:thm:m4-corrected-leading-pairings} and using
$\sum_t\theta_{i,t}=0$ gives
\begin{equation*}
|\partial_Lf|
\le C\left(\delta^2+\varepsilon^2L\right)
\le C\mathcal E_0.
\end{equation*}
Set
$\widetilde W_T:=W_T-(B_0+1)\partial_L.$
Then $\widetilde W_T\ell_t
=
W_T\ell_t-(B_0+1)
\le-1$
and
$-df[\widetilde W_T]
\ge
cA_0-C\mathcal E_0.$
Thus, if $A_0\ge M\mathcal E_0$ and $M$ is sufficiently large,
\begin{equation*}
-df[\widetilde W_T]\ge c'A_0,
\qquad
\widetilde W_T\ell_t\le-1.
\end{equation*}

On $\mathcal D_L$, set
$s:=-\operatorname{sign}\mathcal L_\infty(p^*,r^*).$
By Theorem~\ref{II:thm:slice-signed},
\begin{equation*}
-df[s\partial_L]
=
-s\,\partial_Lf
\ge c_{p^*}\Psi_L.
\end{equation*}
Choose $\tau>0$ such that
$|\tau W_T\ell_t|\le\frac12
\, \text{for every }t,$
and define
$W:=s\partial_L+\tau W_T.$
Its scale velocities satisfy
\begin{equation*}
\frac12
\le
-\operatorname{sign}\bigl(\mathcal L_\infty(p^*,r^*)\bigr)
\,W\ell_t
\le
\frac32.
\end{equation*}
Furthermore,
\begin{equation*}
-df[W]
=
-s\,\partial_Lf
+
\tau\frac{|a|^2}{\max\{|a|,\Psi_L\}}
\ge
c\left(|a|+\Psi_L\right).
\end{equation*}
Indeed, the second term controls $|a|$ when $|a|\ge\Psi_L$, while the
first controls both $\Psi_L$ and $|a|$ when $|a|<\Psi_L$.

Finally, Proposition~\ref{II:prop:ambient-lift} with $h=0$ gives
$\|DG_z[W]\|_{\mathcal H}\le C$,
throughout the signed box.
\end{proof}

	\subsection{Amplitude bounds, ambient lifts and scale signs}

	The small-gradient condition controls the relative amplitudes on the full
	bubbling chart.

	\begin{proposition}\label{II:prop:full-relative}
		Fix separation and scale comparability as in Section~\ref{II:sec:bubbling-chart}.
		There exist $C$, $L_0$ and $\epsilon_*>0$ such that for
		$u\in V^\eta_{\rho^*}(N,\epsilon_0)$, $0<\epsilon_0\le\epsilon_*$, $L\ge L_0$,
		its original modulation satisfies
		\begin{equation}\label{II:eq:full-relative}
			\max_{i,t}L|b_{i,t}-b_i^0|\le C,\qquad
			\Delta_\Omega\le C(\epsilon_0^2+L^{-1}+\epsilon_0 L^{-1/2}).
		\end{equation}
		In particular, outside the sector $|b^0|> M/L$, the entire logarithmic amplitude $L|b|$ is bounded by $C+M$. The constants may depend on the chart margins.
	\end{proposition}
	
\begin{proof}
By \cite[Lemma~2.3]{AFHZI},
$\sqrt L\,|b|+\|w\|\le C\epsilon_0$. Since
$\|T_{i,s}\|\le CL^{-1/2}$ and
$\|\nabla J_{\rho^*}(u)\|<\epsilon_0$,
\[
|DJ_{\rho^*}(u)[T_{i,s}]|
\le C\epsilon_0L^{-1/2}.
\]
Equation~\eqref{II:eq:disk-amplitude-test} at $t=0$ gives
\begin{align*}
\sigma_i|1-N\omega_{i,s}|
&\le |(KD_\sigma b_s)_i|
+C\epsilon_0L^{-1/2}
+C\frac{1+\|w\|}{L}\le C\left(\epsilon_0L^{-1/2}+L^{-1}\right).
\end{align*}
Using \eqref{II:eq:disk-fraction-test}, we obtain
\begin{align*}
|1-N\Omega_{i,s}|
&\le |1-N\omega_{i,s}|+N|\omega_{i,s}-\Omega_{i,s}|\\
&\le C\left(
\epsilon_0L^{-1/2}+L^{-1}
+e^{-\eta L}+|b|\,\|w\|+\|w\|^2
\right)\\
&\le C\left(
\epsilon_0^2+L^{-1}+\epsilon_0L^{-1/2}
\right).
\end{align*}
Hence
$\Delta_\Omega\le
C(\epsilon_0^2+L^{-1}+\epsilon_0L^{-1/2})$.
For $\epsilon_*$ small and $L_0$ large,
$N\Omega_{i,s}\in[1/2,3/2]$ for every $i,s$.

Since $\Omega_{i,t}=\mathfrak M_{i,t}/\mathfrak S_i$,
\[
\log\frac{\mathfrak M_{i,t}}{\mathfrak M_{i,s}}
=
\log\frac{\Omega_{i,t}}{\Omega_{i,s}}
=
O(1).
\]
The definition of $\mathfrak M_{i,t}$ also gives
\begin{align*}
\log\frac{\mathfrak M_{i,t}}{\mathfrak M_{i,s}}
&=
\sigma_iL(b_{i,t}-b_{i,s})
+q_i(\delta_t-\delta_s)
+\sigma_i(b_{i,t}\delta_t-b_{i,s}\delta_s)
+O(1).
\end{align*}
The $q_iL$ terms cancel, while the remaining scale and prefactor terms
are bounded. Thus $\sigma_iL(b_{i,t}-b_{i,s})=O(1)$ uniformly in
$i,s,t$. Since $b_i^0=N^{-1}\sum_s b_{i,s}$,
\[
L|b_{i,t}-b_i^0|
=
\frac1N\left|
\sum_sL(b_{i,t}-b_{i,s})
\right|
\le C.
\]
This proves \eqref{II:eq:full-relative}. If $|b^0|\le M/L$, then
$L|b|\le L|b^0|+\max_{i,t}L|b_{i,t}-b_i^0|\le M+C$.
\end{proof}

	The following result lifts a reduced field to the ambient space without
	changing its parameter velocities.

	\begin{proposition}\label{II:prop:ambient-lift}
		On any of the small-defect bounded-$L|b|$ cells of Proposition~\ref{II:prop:normal-field}, let $B_z$ be its base, $v_z$ the normal solution, $G_z=B_z+v_z$, $f=J_{\rho^*}\circ G$ and $u=G_z+h$, $\Pi_z h=h$, with the normal projection of Proposition~\ref{II:prop:normal-field}. Suppose $Y$ is a locally Lipschitz parameter field made of bounded combinations of normalized position and log-scale directions and $\ell_t^{-1}\partial_{b_{i,t}}$. If $-df[Y]\ge c A(z)$ on the sector, there is a bounded locally Lipschitz ambient lift with the same parameter velocities and
		\begin{equation*}
		-DJ_{\rho^*}[W]\ge c A(z)+c'\|h\|^2
		\end{equation*}
		for all sufficiently small $\|h\|$. If $A(z)\ge c_0\Psi_L$, the lift can instead satisfy $-DJ_{\rho^*}[W]\ge c_1(A(z)+\|h\|)$. The constants are uniform on the indicated sector.
	\end{proposition}
	
\begin{proof}
Let $(T_a)$ be the Dirichlet-normalized tangent frame, decomposed into
amplitude and geometric columns. The corresponding rows satisfy
\[
|DJ_{\rho^*}(u)[T_a^{\rm geo}]|\le C,
\qquad
|DJ_{\rho^*}(u)[T_a^{\rm amp}]|\le C\sqrt L.
\]
Indeed, the geometric tests have bounded Dirichlet and supremum norms,
while the amplitude tests have supremum norm $O(\sqrt L)$. The corrected
terms $C$ and $Q_{i,t}$ have bounded Dirichlet norm.

The frame estimates in \cite[Lemmas~A.2--A.4]{AFHZI} give
\[
\|D_YT_a^{\rm amp}\|\le CL^{-1/2},
\qquad
\|D_YT_a^{\rm geo}\|\le C,
\qquad
\|(\mathbb G^{-1})_{\rm amp,geo}\|\le CL^{-1/2}.
\]
For $k\in\operatorname{ran}\Pi_z$, differentiation of $\Pi_zk=k$ gives
$(D_Y\Pi_z)k\in\operatorname{ran}(I-\Pi_z)$ and
\[
\bigl((D_Y\Pi_z)k,T_a\bigr)=-(k,D_YT_a).
\]
Thus
\begin{equation*}
(D_Y\Pi_z)k
=
\sum_aO(\|k\|)T_a^{\rm geo}
+
\sum_aO\left(\frac{\|k\|}{\sqrt L}\right)T_a^{\rm amp}.
\end{equation*}
Consequently,
\[
|DJ_{\rho^*}(u)[(D_Y\Pi_z)k]|\le C\|k\|.
\]

Set $g:=\nabla J_{\rho^*}(G_z)$. Normal stationarity gives
$\Pi_zg=0$. Since $D_Y\Pi_z$ is self-adjoint,
\[
\|(D_Y\Pi_z)g\|
=
\sup_{\substack{k\in\operatorname{ran}\Pi_z\\ \|k\|=1}}
|(g,(D_Y\Pi_z)k)|
\le C.
\]

Let $H:=D^2J_{\rho^*}(G_z)$. Differentiating
$\Pi_z\nabla J_{\rho^*}(G_z)=0$ and $\Pi_zv_z=v_z$ gives
\begin{equation*}
(\Pi_zH\Pi_z)\Pi_zD_Yv_z
=
-(D_Y\Pi_z)g
-\Pi_zH\left[D_YB_z+(D_Y\Pi_z)v_z\right].
\end{equation*}
The preceding coefficient bounds imply
$\|(D_Y\Pi_z)v_z\|\le C\|v_z\|$. The anchored moment estimates give
\[
\left\|
(D_Y\Pi_z)g
+\Pi_zH\left[D_YB_z+(D_Y\Pi_z)v_z\right]
\right\|
\le C.
\]
Since $\Pi_zH\Pi_z$ has a uniformly bounded inverse,
$\|\Pi_zD_Yv_z\|\le C$, and hence
\[
\|D_YG_z\|
=
\|D_YB_z+D_Yv_z\|
\le C.
\]
The normal part of $D_YG_z$ has uniformly bounded weighted moments, and
its tangent part has bounded supremum norm.

Transport $h$ by  $\dot h=(D_Y\Pi_z)h$
and set $Z:=D_YG_z+(D_Y\Pi_z)h.$
Since $(g,h)=0$, differentiation along this transported section gives
$H[h,D_YG_z]
+
DJ_{\rho^*}(G_z)[(D_Y\Pi_z)h]
=
0.$
Moreover,
\[
|H[h,(D_Y\Pi_z)h]|\le C\|h\|^2,
\qquad
|D^3J_{\rho^*}(G_z+\tau h)[h,h,D_YG_z]|
\le C\|h\|^2
\]
for $0\le\tau\le1$. Taylor expansion therefore yields
$DJ_{\rho^*}(u)[Z]
=
df[Y]+O(\|h\|^2).$

Set $R:=\Pi_z\nabla J_{\rho^*}(u)$. By
Proposition~\ref{II:prop:normal-field},
\[
c\|h\|\le\|R\|\le C\|h\|.
\]
Choose $K_0$ sufficiently large and define
$\widehat W:=Z-K_0R$. Since $-df[Y]\ge cA(z)$,
\begin{align*}
-DJ_{\rho^*}(u)[\widehat W]
&=
-DJ_{\rho^*}(u)[Z]+K_0\|R\|^2 \ge
cA(z)-C\|h\|^2+K_0c^2\|h\|^2\\
&\ge
cA(z)+c'\|h\|^2.
\end{align*}
Both $R$ and the normal correction have zero parameter velocity, so
$\widehat W$ has the same parameter velocities as $Y$.

For the second estimate, add
$W_\nu:=-\frac{R}{\max\{\Psi_L,\|R\|\}}.$
Then
\[
-DJ_{\rho^*}(u)[W_\nu]
=
\frac{\|R\|^2}{\max\{\Psi_L,\|R\|\}}.
\]
If $\|h\|\ge\Psi_L$, this term is bounded below by $c\|h\|$. If
$\|h\|<\Psi_L$, the assumption $A(z)\ge c_0\Psi_L$ gives
$A(z)\ge c_0\|h\|$. Hence
\[
-DJ_{\rho^*}(u)[\widehat W+W_\nu]
\ge c_1\bigl(A(z)+\|h\|\bigr).
\]
The frame bounds, finite-scale smoothness, and positivity of the
denominator imply that both fields are bounded and locally Lipschitz.
\end{proof}

	At an unbalanced finite-scale center, the amplitude quadratic form fixes the
	common-scale sign.

\begin{theorem}
Assume $m_i=4$ for every $i$. At the finite-scale center of
Theorem~\ref{II:thm:finite-scale-covector-center}, set
$\xi^\alpha:=\frac{D_\sigma c_\circ^\alpha}{L}.$
Then
\begin{equation*}
\partial_L f(p,r_\circ,b_\circ,L)
=
-\pi\sum_\alpha(\xi^\alpha)^TK\xi^\alpha
+o(L^{-2}),
\end{equation*}
where $\partial_L$ is taken with $p,r,b$ fixed at
$(p,r_\circ,b_\circ)$. The remainder is uniform on compact separated
location sets and bounded ranges of the center variables.

If $|c_\infty(p)|\ge\eta>0$ on such a compact set, then, for all
sufficiently large $L$,
\begin{equation*}
-df[-\partial_L]
=
\partial_Lf
\ge c_\eta L^{-2}.
\end{equation*}

Furthermore, on any bounded $L|b|$ cell contained in
$\{|c_\infty(p)|\ge\eta\}$, there exist
$\delta_1,c_1,L_1>0$ such that, whenever
$\Delta_\Omega\le\delta_1$ and $L\ge L_1$, a bounded ambient field $W$
satisfies
\begin{equation*}
\|W\|\le C,\qquad
\dot\ell_t\le-c_1\quad(1\le t\le N),
\qquad
-DJ_{\rho^*}(G_z+h)[W]
\ge
c_1\left(A_0+L^{-2}+\|h\|\right)
\end{equation*}
for every sufficiently small normal displacement
$h\in\operatorname{ran}\Pi_z$. Here $A_0$ is defined in
\eqref{II:eq:pg-transverse}.
\end{theorem}
\begin{proof}
It suffices to consider arbitrary sequences in the stated compact sets
with $L\to\infty$. Set 
\[
B:=B_z,\qquad v:=\bar v_z,\qquad S:=\partial_LB.
\]
Assume first that $|\mathcal C|=O(L^{-2})$ and $L|b|\le C$. The amplitude
rows give
\[
|\theta|=O(L^{-1}),
\qquad
|b|+|\theta|+\|v\|=O(L^{-1}).
\]
After extraction, $L\xi_t\to u_t\in\mathbb R^n$. On the $t$th core set
\[
k(Y):=k_1(Y),\qquad
P_0(Y):=-2\log(1+a^2|Y|^2),
\]
\[
S_0:=4+Y\cdot\nabla P_0=\frac4{1+a^2|Y|^2},
\qquad
Z:=S_0-2,\qquad T:=Y\cdot\nabla .
\]
Then
\[
Tk+2k=kZ,\qquad T P_0=Z-2,\qquad \int_{\mathbb R^2}kS_0=2.
\]

Define the rescaled normal correction by
\[
w_{i,t}^{(L)}(Y)
:=
L\,v_i\!\left(y_{p_t}^{-1}(Y/\lambda_t)\right).
\]
Proposition~\ref{II:prop:m4-amplitude-corrected-sharp-residual} and
Lemma~\ref{II:lem:normal-model} give, after extraction,
\[
w_t^{(L)}\rightharpoonup w_t
\quad\text{in }H^1_{\mathrm{loc}}(\mathbb R^2)^n,
\qquad
w_t^{(L)}\to w_t
\quad\text{in }L^p(k\,dY)^n
\]
for every finite $p$. The convergence also holds with fixed logarithmic
weights. The corrected anchor differs from the original anchor by
$O(\varepsilon)$, so $\int_{\mathbb R^2}kw_{i,t}\,dY=0.$

Fix one core and suppress $t$. Define
\begin{equation*}
\mathscr H(w,\phi)
:=
\sum_{i,j}a^{ij}\int_{\mathbb R^2}
\nabla w_j\cdot\nabla\phi_i
-
\sum_i2\pi\sigma_i\int_{\mathbb R^2}kw_i\phi_i.
\end{equation*}
Set
\[
a_i':=\frac{u_i}{4},
\qquad
\Pi:=P_0-\int_{\mathbb R^2}kP_0.
\]
The limit $w$ satisfies
\begin{equation}\label{II:eq:core-first-response}
\mathscr H(w,\phi)
=
\sum_i2\pi\sigma_i a_i'
\int_{\mathbb R^2}k\Pi\phi_i
\end{equation}
for every $\phi\in H^1(S^2)^n$.

Indeed, Theorem~\ref{II:thm:m4-corrected-leading-pairings} gives
amplitude rows $O(L^{-3/2})$ and geometric rows
$O(L^{-2}+\varepsilon)$. Since
$q^TD_\sigma\theta^\alpha=O(L^{-2})$ and
$\sum_t\theta_{i,t}=0$, normal stationarity and the Gram bounds yield
$L\|\nabla J_{\rho^*}(G_z)\|\to0$. On a fixed rescaled core, the free
and potential densities are
\begin{align*}
2\pi\bigl[\sigma_i+(K\xi_t)_i\bigr]k
&+O_\phi(e^{-\eta L}),\\
2\pi\sigma_i(1-\theta_{i,t})k
\left[
1+(s_{i,t}-1)\Pi+O((s_{i,t}-1)^2)
\right]
&+O_\phi(e^{-\eta L}),
\end{align*}
respectively. The amplitude equations give
\[
(K\xi_t)_i+\sigma_i\theta_{i,t}=O(L^{-2}),
\qquad
L(s_{i,t}-1)\to a_i'.
\]
Moreover, $L\|v\|^2=o(1)$. Multiplication by $L$ and passage to the
limit therefore give \eqref{II:eq:core-first-response} for compactly
supported tests. Since $\Pi$ belongs to every finite
$L^p(S^2)$, point capacity and continuity extend the identity to
$H^1(S^2)^n$.

We next compare the common-scale pairings at $B+v$ and $B$. The preceding
convergences give
\begin{align*}
&L^2\left(
DJ_{\rho^*}(B+v)[S]-DJ_{\rho^*}(B)[S]
\right)\to 
-\sum_t\sum_i2\pi\sigma_i
\int_{\mathbb R^2}k
\left[
a_i'(S_0+Z\Pi)w_{i,t}
+\frac12Zw_{i,t}^2
\right].
\end{align*}
To obtain this limit, use
\[
\int_Mv_i\,d\mu_i(B)=O(L^{-2}),
\qquad
(K\xi_t)_i+\sigma_i\theta_{i,t}=O(L^{-2}),
\qquad
\|C\|+\|Q_{i,t}\|=O(\varepsilon).
\]
The linear term follows by expanding the density against
$s_{i,t}S_0-2$. The quadratic term equals
$-\pi\sigma_i\int kZw_{i,t}^2$. The remaining Taylor term is
$O(\|v\|^3)=o(L^{-2})$. Proposition~\ref{II:prop:m4-global-denominator},
the anchored Gaussian moments, and the $O(\varepsilon)$ bounds for
$C,Q_{i,t}$ give $o(L^{-2})$ for the spatial, cutoff, and exterior
errors.

We show that the displayed limit vanishes. Elliptic regularity gives
$w\in W^{2,p}(S^2)^n$ for every finite $p$, so $Tw$ is admissible in
\eqref{II:eq:core-first-response}. Conformal invariance gives
\[
\sum_{i,j}a^{ij}\int_{\mathbb R^2}
\nabla w_j\cdot\nabla(Tw_i)=0.
\]
Since $T\Pi=Z-2$ and $Tk+2k=kZ$,
\begin{align*}
\frac12\sum_i2\pi\sigma_i\int kZw_i^2
&=
\sum_i2\pi\sigma_i a_i'\int k\Pi\,Tw_i=
-\sum_i2\pi\sigma_i a_i'
\int k(Z\Pi+Z-2)w_i.
\end{align*}
Because $\int kw_i=0$ and $S_0=Z+2$,
\[
\sum_i2\pi\sigma_i
\int k
\left[
a_i'(S_0+Z\Pi)w_i+\frac12Zw_i^2
\right]
=0.
\]
Hence
\[
DJ_{\rho^*}(B+v)[S]
=
DJ_{\rho^*}(B)[S]+o(L^{-2}).
\]

Lemma~\ref{II:lem:ABL-first-order-lift}, the
$O(L^{-3/2}+\varepsilon)$ tangent covector, and $\|v\|=O(L^{-1})$ give
\[
\partial_Lf
=
DJ_{\rho^*}(G_z)[S]+o(L^{-2}).
\]
The corrected energy identity \eqref{II:eq:corrected-energy-identity}
and the $bQ$ estimate in the proof of
Theorem~\ref{II:thm:slice-signed} give
\[
DJ_{\rho^*}(B)[S]
=
\partial_L\mathcal J_0+o(L^{-2}).
\]
By \eqref{II:eq:slice-quadratic} and the amplitude equations,
\[
\partial_L\mathcal J_0
=
-\pi\sum_t\xi_t^TK\xi_t+O(L^{-3}).
\]
At the finite-scale center, $\xi_t=\sum_\alpha e_t^\alpha\xi^\alpha,
\, 
\xi^\alpha=\frac{D_\sigma c_\circ^\alpha}{L}.$
The orthonormality of $(e^\alpha)$ gives $\sum_t\xi_t^TK\xi_t
=
\sum_\alpha(\xi^\alpha)^TK\xi^\alpha.$
Therefore
\[
\partial_Lf(p,r_\circ,b_\circ,L)
=
-\pi\sum_\alpha(\xi^\alpha)^TK\xi^\alpha+o(L^{-2}).
\]

At the center, $q^TK\xi^\alpha=0$. With $x=K\xi^\alpha$,
\cite[Lemma~3.1]{AFHZI} gives
\[
(\xi^\alpha)^TK\xi^\alpha
=
x^TAx
\le-c|K\xi^\alpha|^2.
\]
Since $c_\circ\to c_\infty$ uniformly, $|c_\infty(p)|\ge\eta$ implies
\[
\sum_\alpha|K\xi^\alpha|^2
=
\frac1{L^2}
\sum_\alpha|KD_\sigma c_\circ^\alpha|^2
\ge c_\eta L^{-2}.
\]
Thus, $-df[-\partial_L]
=
\partial_Lf
\ge c_\eta L^{-2}$,
for all sufficiently large $L$.

It remains to construct the field on the entire unbalanced cell. Choose
$M$ as in Proposition~\ref{II:prop:slice-field}. In the outer region
$A_0>M\mathcal E_0$, that proposition gives a bounded field
$W_{\rm out}$ such that
\[
W_{\rm out}\ell_t\le-1,
\qquad
-df[W_{\rm out}]\ge cA_0.
\]
In the inner region $A_0<2M\mathcal E_0$, the amplitude rows give
$\delta\le C(L^{-1}+A_0)$, and hence
\[
A_0
\le
C_M\left(
L^{-2}+\frac{A_0}{L}+A_0^2
+\frac{\varepsilon}{L}+\varepsilon^2L
\right).
\]
For $\Delta_\Omega\le\delta_1$ sufficiently small and $L$ sufficiently
large, absorption gives $A_0\le C_ML^{-2}.$

On $\{|c_\infty(p)|\ge\eta\}\cap\{A_0>M\mathcal E_0\}$, one also has
$A_0+L^{-2}\le C_\eta A_0.$
Otherwise, along a sequence,
$A_0=o(L^{-2})$ and $\mathcal E_0=o(L^{-2})$. Thus
$L|b|\to0$ and $\theta\to0$, and the limiting mass equations give a
balanced location, contrary to $|c_\infty|\ge\eta$.

In the inner region, the mass equations give
\[
H^\alpha+qr_\alpha+D_\sigma c^\alpha=o(1),
\qquad
q^TKD_\sigma c^\alpha=o(1).
\]
The relative block is invertible, so
$c-c_\infty(p)\to0,\, 
r-r_\infty(p)\to0.$
Consequently,
\[
\sum_t|L\xi_t|^2\ge c_\eta,
\qquad
q^TK(L\xi_t)\to0.
\]
The preceding normal-response calculation therefore gives
\[
-df[-\partial_L]
=
\partial_Lf
\ge c_\eta L^{-2}
\]
throughout the inner region, including points with nonzero amplitude
covectors.

Let $\chi:[0,\infty)\to[0,1]$ satisfy
$\chi=0$ on $[0,M]$ and $\chi=1$ on $[2M,\infty)$. Set
\[
W
:=
\left(1-\chi\left(\frac{A_0}{\mathcal E_0}\right)\right)(-\partial_L)
+
\chi\left(\frac{A_0}{\mathcal E_0}\right)W_{\rm out}.
\]
Both fields satisfy $\dot\ell_t\le-1$ on the overlap, and the preceding
bounds give
\[
\|W\|\le C,\qquad
\dot\ell_t\le-1,\qquad
-df[W]\ge c_\eta(A_0+L^{-2}).
\]
Finally, $L^{-2}\ge\Psi_L$ for large $L$. Proposition~\ref{II:prop:ambient-lift}
therefore yields a bounded ambient lift, still denoted by $W$, such that
\[
\dot\ell_t\le-c_1,
\qquad
-DJ_{\rho^*}(G_z+h)[W]
\ge
c_1\left(A_0+L^{-2}+\|h\|\right)
\]
for every sufficiently small normal displacement $h$.
\end{proof}

	The mixed regime admits the same signed descent after adding a normal
	displacement.

	\begin{corollary}
		In the signed region of Theorem~\ref{II:thm:moving-mixed}, there exist
$c,C,h_0,L_0>0$ and a locally Lipschitz ambient field $W$ such that,
for $L\ge L_0$ and $h\in\operatorname{ran}P_z$ with $\|h\|\le h_0$,
\begin{equation*}
\|W\|\le C,\qquad
D\ell_t(u)[W(u)]=-\operatorname{sign}\mathcal L_\infty(p^*,r^*)
\quad(1\le t\le N),
\end{equation*}
and
\begin{equation*}
-DJ_{\rho^*}(G_z+h)[W]
\ge c\left(\Psi_L+\|h\|\right).
\end{equation*}
	\end{corollary}
	\begin{proof}
		Apply Proposition~\ref{II:prop:ambient-lift} to $Y=-\operatorname{sign}(\mathcal L_\infty(p^*,r^*))\partial_L$, with $A=\Psi_L$.
	\end{proof}

The next example shows that, in the mixed regime, a decreasing-scale
field cannot attain the rate $\Psi_L$ away from balanced locations.

\begin{proposition}
There exist smooth data on a flat unit-area torus satisfying
\textup{\ref{II:assump:H1}--\ref{II:assump:H2}} and
\textup{\ref{II:item:a1}--\ref{II:item:a2}}, with $N=1,\, \widehat m=\frac52,$
and points $p^*,p_0\in M$ such that
\[
\nabla\ff(p^*)=0,\qquad
\mathcal L_\infty(p^*)<0,\qquad
\nabla\ff(p_0)\ne0,
\]
where $d_g(p_0,p^*)$ can be arbitrarily small.

For every sufficiently large $L$, the original normal graph contains a
point $u_L$ in the full small-gradient chart with modulation location
$p(u_L)=p_0$. Moreover, for every $W\in\mathcal H$,
\begin{equation}\label{II:eq:mixed-obstruction}
DL(u_L)[W]\le0
\quad\to \quad
-DJ_{\rho^*}(u_L)[W]
\le Ce^{-L}\|W\|.
\end{equation}
Here $L(u)$ denotes the common scale in the original modulation.

In particular, if $\|W_L\|\le C_0$ and $DL(u_L)[W_L]\le0$, then
\[
-DJ_{\rho^*}(u_L)[W_L]
\le CC_0e^{-L}
=o(e^{-L/2})
=o(\Psi_L).
\]
Hence, no bounded field with nonincreasing common scale can satisfy
$-DJ_{\rho^*}[W]\ge c\Psi_L$ uniformly at these points.
\end{proposition}

\begin{proof}
Let
\[
A=K=
\begin{pmatrix}
0&1\\
1&0
\end{pmatrix},
\qquad
\sigma=\left(10,\frac52\right).
\]
Then
\[
m=A\sigma=\left(\frac52,10\right),
\qquad
4\sum_i\sigma_i=50=\sigma^TA\sigma.
\]
Thus $\rho^*$ is resonant, $\widehat m=5/2$, and only the first
component contributes to $\mathcal L_\infty$.

On the flat square torus, set
\[
\phi(x):=2-\cos(2\pi x_1)-\cos(2\pi x_2).
\]
The function $\phi$ is Morse, nonnegative, and has a unique zero
$p^*$. Initially, take $h_2=1$ and $h_1=e^{-B\phi}$, $B>0$. Since the
Robin function is constant,
\begin{align*}
D_{1,1}^{(B)}(p^*)-D_{1,1}^{(0)}(p^*)
&=
\int_M
\left(e^{-B\phi(x)}-1\right)
e^{2\pi m_1(G^g(x,p^*)-R^g(p^*))}\,dV_g.
\end{align*}
For each $B$, the integral is finite because
$\phi(x)=O(d_g(x,p^*)^2)$ and $m_1<4$. Moreover,
\[
D_{1,1}^{(B)}(p^*)\to -\infty
\qquad\text{as }B\to\infty.
\]
Fix $B$ such that $D_{1,1}^{(B)}(p^*)<0$.

Choose $\chi\ge0$, $\chi\not\equiv0$, supported away from the critical
set of $\phi$, and replace $h_1$ by
\[
h_{1,s}:=e^{-B\phi}(1+s\chi).
\]
For every critical point $p$,
\[
D_{1,1}^{(s)}(p)
=
D_{1,1}^{(0)}(p)+sI(p),
\qquad
I(p)>0.
\]
For sufficiently small $|s|$, the critical locations are unchanged.
Choosing $s$ outside a finite set gives
\[
D_{1,1}^{(s)}(p)\ne0
\]
at every critical point while preserving
$D_{1,1}^{(s)}(p^*)<0$. Hence \textup{(A1)--(A2)} hold and
$\mathcal L_\infty(p^*)<0$. Since balance is automatic for $N=1$,
continuity gives a point $p_0\ne p^*$, arbitrarily close to $p^*$, such
that
\[
\nabla\ff(p_0)\ne0,
\qquad
\mathcal L_\infty(p_0)<0.
\]

We next construct an amplitude zero at $(p_0,L)$. Fix
\[
\frac{\widehat m-2}{2}<\alpha<\widehat m-2,
\qquad
s_L:=e^{-\alpha L}.
\]
Let $G(p,L,b)=U_z+\bar v_z$ be the original normal graph and set
$f:=J_{\rho^*}\circ G$. On $L|b|\le1$, the mixed anchored estimates give
\[
\|\bar v_z\|\le C(|b|+s_L),
\qquad
\|\nabla J_{\rho^*}(U_z)\|
\le C\sqrt L\,(|b|+s_L).
\]
For $Y_i:=L^{-1}\partial_{b_i}$,
\begin{equation*}
df[Y_i]
=
2\pi\sigma_i(KD_\sigma b)_i
+
O\left(
s_L+\frac{|b|+s_L}{L}+(|b|+s_L)^2
\right).
\end{equation*}
Indeed, the raw source discrepancy has total mass
\[
2\pi(KD_\sigma b)_i
+
O\left(e^{-(\widehat m-2)L}(1+|b|)\right),
\]
while
\[
\frac{PW_{i,1}}{L}
=
\sigma_i
+
O\left(\frac{1+\log(2+|Y|)}{L}\right)
\]
on the core. The anchored moments give
$O((|b|+s_L)/L)$ for the nonconstant part. The normal transfer satisfies
\[
O\left(\frac{\|\bar v_z\|}{L}\right)
=
O\left(\frac{|b|+s_L}{L}\right).
\]
Finally, Lemma~\ref{II:lem:ABL-first-order-lift}, the
$O(L^{-1})$ frame derivatives, and
$\|\nabla J_{\rho^*}(G_z)\|\le CL(|b|+s_L)$ give an
$O((|b|+s_L)^2)$ graph-lift error.

Define
\[
\mathcal A_L(b)
:=
\left(
\frac{df[Y_i]}{2\pi\sigma_i}
\right)_{i=1}^2.
\]
Then
\[
\mathcal A_L(b)
=
KD_\sigma b+E_L(b),
\qquad
|E_L(b)|
\le
C\left(
s_L+\frac{|b|+s_L}{L}+(|b|+s_L)^2
\right).
\]
Since $KD_\sigma$ is invertible, choose $R>0$ such that, for all large
$L$,
\[
|E_L(b)|<|KD_\sigma b|
\qquad\text{on }|b|=Rs_L.
\]
Therefore
\[
\deg(\mathcal A_L,B_{Rs_L},0)
=
\deg(KD_\sigma,B_{Rs_L},0)
\ne0.
\]
Thus there exists $b_L$ such that
\[
|b_L|\le Rs_L,
\qquad
D_bf(p_0,L,b_L)=0.
\]

Set $u_L:=G(p_0,L,b_L).$
The preceding estimates give
\[
\|\nabla J_{\rho^*}(u_L)\|\to 0,
\qquad
\|u_L-U_{p_0,L,0}\|
=
O(\sqrt L\,s_L)\to 0.
\]
Hence $u_L$ belongs to the full small-gradient chart for all sufficiently
large $L$. Since $b_L=O(e^{-\alpha L})$, Theorem~\ref{II:thm:moving-mixed}
applies and gives
\[
\partial_Lf(p_0,L,b_L)
=
\pi\mathcal L_\infty(p_0)e^{-(\widehat m-2)L}
+o\left(e^{-(\widehat m-2)L}\right)
\le-c\Psi_L.
\]

We now prove \eqref{II:eq:mixed-obstruction}. Transport the isothermal
disks and cutoffs by translations of the flat torus. In the translated
coordinates, $U_z$ and $E_z$ are independent of $p$; only
$h_i(p+\cdot)$ depends on $p$. Normal stationarity therefore gives
\[
|D_pf|
\le
\sum_i\rho_i^*\|\nabla\log h_i\|_{L^\infty(M)}
\le C.
\]

For $W\in\mathcal H$, denote
\[
W
=
DG[\dot p,\dot L,\dot b]+k,
\qquad
k\in E_z.
\]
Differentiation of the modulation orthogonality conditions gives
\[
|\dot p|\le Ce^{-L}\|W\|,
\qquad
\dot L=DL(u_L)[W].
\]
Since $k\in E_z$ and $D_bf(p_0,L,b_L)=0$,
\[
DJ_{\rho^*}(u_L)[W]
=
D_pf[\dot p]+\partial_Lf\,\dot L.
\]
If $DL(u_L)[W]=\dot L\le0$, then
$\partial_Lf\,\dot L\ge0$. Consequently,
\begin{align*}
-DJ_{\rho^*}(u_L)[W]
&=
-D_pf[\dot p]-\partial_Lf\,\dot L
\le
|D_pf|\,|\dot p|\le
Ce^{-L}\|W\|.
\end{align*}
This proves \eqref{II:eq:mixed-obstruction}.

Finally, Proposition~\ref{II:prop:ambient-lift} gives
$\|DG[\partial_L]\|\le C$. Hence
\[
c\Psi_L
\le
|\partial_Lf|
=
|DJ_{\rho^*}(u_L)[DG[\partial_L]]|
\le
C\|\nabla J_{\rho^*}(u_L)\|,
\]
so $\|\nabla J_{\rho^*}(u_L)\|\ge c'\Psi_L$. Since
\[
\frac{e^{-L}}{\Psi_L}
=
e^{-L/2}\to 0,
\]
no bounded field with nonincreasing common scale can have uniform
descent of order $\Psi_L$ at these points.
\end{proof}

	\subsection{A field in the full modulation chart}

	Use the original minimizing chart $u=U_z+w$, $P_zw=w$. In this subsection $P_z$ is the normal projection and $P(u)$ is the logarithmic mass of Section~\ref{II:sec:exterior}. The Euclidean norm $|\theta|$ satisfies $\Delta_\Omega\le |\theta|\le\sqrt{nN}\,\Delta_\Omega$, and $\sum_s\theta_{i,s}=0$ for each component.
	
	Dominant mass imbalance defines a uniform descending direction.

	\begin{lemma}\label{II:lem:tau-flow}
		Fix a compact separated location cell, scale comparability and a small fixed bound on $\sqrt L|b|+\|w\|$. Suppose $|\theta|\le\delta$ for sufficiently small fixed $\delta$. Set
		\begin{equation*}
		B=|b|+\|w\|^2+L^{-1}.
		\end{equation*}
		There are $M,L_0,c,C>0$ such that, for $L\ge L_0$, $N\ge2$ and $|\theta|\ge MB$, the original chart admits a field $Y_\theta$ of norm at most $C$, with fixed centers, scales and common amplitudes, satisfying
		\begin{equation}\label{II:eq:tau-descent}
			-DJ_t[Y_\theta]\ge c\sqrt L\,|\theta|\quad(0\le t\le\kappa),
			\qquad D|\theta|[Y_\theta]\ge c\sqrt L .
		\end{equation}
		Here $0<\kappa\le1/2$. In this region $\|DJ_{\rho^*}\|\le C\sqrt L\,|\theta|$. A segment following this field while all the indicated inequalities persist has length in flow time at most $\delta/(c\sqrt L)$. When $N=1$, the centered defect is identically zero.
	\end{lemma}
	
\begin{proof}
Prescribe
\begin{equation*}
\dot b_{i,s}
=
-\frac{\theta_{i,s}}{\sigma_i\sqrt L\,|\theta|},
\qquad
\dot p=\dot\ell=0,
\qquad
\dot w=(D_{\dot b}P_z)w.
\end{equation*}
Since $\sum_s\theta_{i,s}=0$,
\[
\dot b_i^0
=
\frac1N\sum_s\dot b_{i,s}
=
0.
\]
Moreover, $(D_{\dot b}P_z)w$ is tangent and
\[
P_z(D_{\dot b}P_z)w=0,
\qquad
\frac{d}{d\tau}(P_zw-w)=0.
\]
Thus the centers, scales, common amplitudes, and constraint $P_zw=w$
are preserved. The Gram estimates in \cite[Appendix~A]{AFHZI} and the
horizontal lift bounds give $\|Y_\theta\|\le C.$

Logarithmic differentiation of
\eqref{II:def:full-local-mass} gives
\begin{align*}
(\log\mathfrak M_{i,s})^{\cdot}
&=
\sigma_i\ell_s\dot b_{i,s}+O(|\dot b|)=
-\sqrt L\,\frac{\theta_{i,s}}{|\theta|}
+O(L^{-1/2}),
\end{align*}
because $\ell_s=L+O(1)$. Since
$\Omega_{i,s}=\mathfrak M_{i,s}/\mathfrak S_i$,
\[
\dot\Omega_{i,s}
=
\Omega_{i,s}
\left[
(\log\mathfrak M_{i,s})^{\cdot}
-\sum_l\Omega_{i,l}
(\log\mathfrak M_{i,l})^{\cdot}
\right].
\]
Using $\theta_{i,s}=1-N\Omega_{i,s}$, we obtain
\begin{align*}
\frac12\frac{d}{d\tau}|\theta|^2
&=
\frac{N\sqrt L}{|\theta|}
\sum_i
\left[
\sum_s\Omega_{i,s}\theta_{i,s}^2
-
\left(\sum_s\Omega_{i,s}\theta_{i,s}\right)^2
\right]
+O\left(\frac{|\theta|}{\sqrt L}\right).
\end{align*}
For each $i$,
\begin{align*}
\sum_s\Omega_{i,s}\theta_{i,s}^2
-
\left(\sum_s\Omega_{i,s}\theta_{i,s}\right)^2
&=
\min_{c\in\mathbb R}
\sum_s\Omega_{i,s}(\theta_{i,s}-c)^2\ge
\frac{1-\delta}{N}
\min_{c\in\mathbb R}
\sum_s(\theta_{i,s}-c)^2=
\frac{1-\delta}{N}
\sum_s\theta_{i,s}^2.
\end{align*}
The last identity follows from $\sum_s\theta_{i,s}=0$. Hence
\[
\frac{d}{d\tau}|\theta|
\ge
(1-\delta)\sqrt L-O(L^{-1/2})
\ge c\sqrt L
\]
for large $L$. Since $|\theta|\le\delta$, the time spent in this sector
is at most $\delta/(c\sqrt L)$.

Let $Y_{i,s}:=\frac1{\sigma_i\ell_s}\partial_{b_{i,s}}.$
The parameter part of $Y_\theta$ is
\[
Y_\theta
=
-\sum_{i,s}
\frac{\ell_s\theta_{i,s}}{\sqrt L\,|\theta|}
\,Y_{i,s}.
\]
The horizontal lift of \eqref{II:eq:disk-amplitude-test} has error
$O((1+\|w\|)/L)$ by Proposition~\ref{II:prop:large-common}. Since
$\ell_s=L+O(1)$, replacing $\ell_s$ by $L$ contributes
$O(L^{-1/2})$. Using \eqref{II:eq:disk-fraction-test} and
$\sum_s\theta_{i,s}=0$, we obtain
\begin{equation*}
-DJ_t(u)[Y_\theta]
=
2\pi(1-t)\sqrt L\,
\frac{\sum_{i,s}\sigma_i\theta_{i,s}^2}{|\theta|}
+
O(\sqrt L\,B).
\end{equation*}
Since $0\le t\le\kappa\le1/2$,
\[
2\pi(1-t)
\frac{\sum_{i,s}\sigma_i\theta_{i,s}^2}{|\theta|}
\ge c|\theta|.
\]
Therefore, if $|\theta|\ge MB$ and $M$ is sufficiently large,
\[
-DJ_t(u)[Y_\theta]
\ge c\sqrt L\,|\theta|,
\qquad
0\le t\le\kappa.
\]
This proves \eqref{II:eq:tau-descent}.

It remains to estimate the full gradient. Let
$\varphi\in\mathcal H$ satisfy $\|\varphi\|=1$. On each core, subtract
the $V_{i,s}$-mean of $\varphi_i$. The mean satisfies
$|\overline\varphi_{i,s}|\le C\sqrt L$, while the oscillation has
uniform logarithmically weighted moments by
Lemma~\ref{II:lem:normal-model}. The constant source discrepancy is
$O(|\theta|+B)$, and the oscillatory discrepancy is
\[
O\left(|\theta|+|b|+e^{-\eta L}\right).
\]
The normal tilt contributes $O(\|w\|)$, and the exterior contribution
is exponentially small. Consequently,
\[
|DJ_{\rho^*}(u)[\varphi]|
\le
C\sqrt L\,(|\theta|+B)+C\|w\|.
\]
Since  $\frac{\|w\|}{\sqrt L}
\le
\frac12\left(\|w\|^2+L^{-1}\right)
\le B,$
we obtain
\[
\|\nabla J_{\rho^*}(u)\|
\le
C\sqrt L\,(|\theta|+B)
\le
C\sqrt L\,|\theta|
\]
on $|\theta|\ge MB$. The velocity formulas are locally Lipschitz on
this sector and extend to its threshold faces.
\end{proof}

	The inequality $D|\theta|[Y_\theta]>0$ concerns dominant imbalance on a fixed scale slice. At $\theta=0$ the mass equations still allow $H^\alpha+qr_\alpha=-D_\sigma c^\alpha+O(L^{-1})$. The no-return argument below also controls $c^\alpha$ and hence gives geometric balance.
	
	We combine the tangent and normal rows into a pseudogradient on the full
	modulation chart.

\begin{proposition}\label{II:prop:chart-segment}
Choose the radius of the original modulation chart sufficiently small
and its lower scale sufficiently large. On the complement of $E_\kappa$,
there exists a bounded locally Lipschitz field $W_{\rm ch}$ satisfying
\[
\|W_{\rm ch}(u)\|\le C.
\]
For each normalized position, log-scale, or amplitude direction $T_a$,
set $\widetilde T_a:=T_a+(D_aP_z)w.$
Then
\begin{equation}\label{II:eq:chart-segment}
\begin{aligned}
&\min\left\{
-DJ_{\rho^*}(u)[W_{\rm ch}(u)],
-DJ_\kappa(u)[W_{\rm ch}(u)]
\right\}\\
&\qquad\ge
c\min_{0\le t\le\kappa}
\left(
\sum_a|DJ_t(u)[\widetilde T_a]|^2
+\|P_z\nabla J_t(u)\|^2
\right)^{1/2}
\ge c'd_\kappa(u).
\end{aligned}
\end{equation}

In the dominant-imbalance region,
\[
|\theta|\ge2MB,
\qquad
W_{\rm ch}=Y_\theta.
\]
On $|\theta|<2MB$, the field is constructed from all tangent rows
$DJ_t[\widetilde T_a]$ and the normal row $P_z\nabla J_t$.

The same construction with the single covector $DJ_{\rho^*}$ gives a
bounded locally Lipschitz field $W_{\rho^*}$ on
$\mathcal H\setminus\operatorname{Crit}(J_{\rho^*})$ such that
\[
-DJ_{\rho^*}(u)[W_{\rho^*}(u)]
\ge
c\|DJ_{\rho^*}(u)\|.
\]
\end{proposition}

\begin{proof}
Equip the normalized parameter space with its Euclidean metric and
$E_z$ with the Dirichlet metric. With
$\widetilde T_a:=T_a+(D_aP_z)w$, define
\[
B_u(v,k)
:=
\sum_a v_a\widetilde T_a+k,
\qquad
v=(v_a),\quad k\in E_z.
\]
The vectors $\widetilde T_a$ differentiate $u=U_z+w$ under the
constraint $P_zw=w$. By \cite[Lemma~2.3 and Appendix~A]{AFHZI},
\[
\|D_aP_z\|\le C,
\qquad
C^{-1}\bigl(|v|+\|k\|\bigr)
\le
\|B_u(v,k)\|
\le
C\bigl(|v|+\|k\|\bigr).
\]
Thus $B_u$ is a uniformly bounded isomorphism. Its adjoint satisfies
\begin{equation*}
B_u^*\nabla J_t(u)
=
\left(
\bigl(DJ_t(u)[\widetilde T_a]\bigr)_a,
P_z\nabla J_t(u)
\right),
\end{equation*}
and hence
\[
\|B_u^*\nabla J_t(u)\|
=
\left(
\sum_a|DJ_t(u)[\widetilde T_a]|^2
+\|P_z\nabla J_t(u)\|^2
\right)^{1/2}.
\]

Set
\[
q_0:=B_u^*\nabla J_{\rho^*}(u),
\qquad
q_\kappa:=B_u^*\nabla J_\kappa(u).
\]
Since $\nabla J_t$ is affine in $t$,
\[
[q_0,q_\kappa]
=
\left\{
B_u^*\nabla J_t(u):0\le t\le\kappa
\right\}.
\]
Let $q_*$ be the point of this segment nearest to zero. Then
\[
(q_*,q_0)\ge\|q_*\|^2,
\qquad
(q_*,q_\kappa)\ge\|q_*\|^2.
\]
At the fixed point $u$, define
\[
W_u:=
-B_u\frac{q_*}{\|q_*\|}.
\]
The uniform bound for $B_u$ gives $\|W_u\|\le C$, while
\begin{align*}
-DJ_{\rho^*}(u)[W_u]
&=
\frac{(q_0,q_*)}{\|q_*\|}
\ge\|q_*\|,\qquad
-DJ_\kappa(u)[W_u]
=
\frac{(q_\kappa,q_*)}{\|q_*\|}
\ge\|q_*\|.
\end{align*}
Moreover,
\begin{align*}
\|q_*\|
&=
\min_{0\le t\le\kappa}
\left(
\sum_a|DJ_t(u)[\widetilde T_a]|^2
+\|P_z\nabla J_t(u)\|^2
\right)^{1/2}\ge
c\min_{0\le t\le\kappa}\|\nabla J_t(u)\|
=
cd_\kappa(u).
\end{align*}

Freeze the direction $-q_*/\|q_*\|$ in a local isometric
trivialization. By continuity, both descent estimates remain valid with
half their constants on a neighborhood of $u$. A locally Lipschitz
partition of unity gives a bounded field $W_{\rm seg}$ satisfying
\eqref{II:eq:chart-segment}.

For the single functional $J_{\rho^*}$, the same construction uses
\[
q:=B_u^*\nabla J_{\rho^*}(u),
\qquad
W_u:=-B_u\frac{q}{\|q\|}.
\]
Then
\[
-DJ_{\rho^*}(u)[W_u]
=
\|B_u^*\nabla J_{\rho^*}(u)\|
\ge
c\|\nabla J_{\rho^*}(u)\|.
\]
A locally Lipschitz partition of unity gives $W_{\rho^*}$ on
$\mathcal H\setminus\operatorname{Crit}(J_{\rho^*})$.

By Proposition~\ref{II:prop:full-relative}, shrinking the chart radius
and increasing its lower scale give $|\theta|\le\delta$. For $N\ge2$,
choose a Lipschitz function $\chi:[0,\infty)\to[0,1]$ such that
\[
\chi(r)=0\quad(r\le M),
\qquad
\chi(r)=1\quad(r\ge2M),
\]
and set
\[
W_{\rm ch}
:=
\left(1-\chi\left(\frac{|\theta|}{B}\right)\right)W_{\rm seg}
+
\chi\left(\frac{|\theta|}{B}\right)Y_\theta.
\]
On $MB\le|\theta|\le2MB$, Lemma~\ref{II:lem:tau-flow} gives
\[
-DJ_{\rho^*}[Y_\theta],
\ -DJ_\kappa[Y_\theta]
\ge
c\sqrt L\,|\theta|.
\]
The gradient estimate in the same lemma and the boundedness of $B_u^*$
give
\[
\min_{0\le t\le\kappa}
\|B_u^*\nabla J_t(u)\|
\le
C\|\nabla J_{\rho^*}(u)\|
\le
C\sqrt L\,|\theta|.
\]
Thus both $W_{\rm seg}$ and $Y_\theta$ satisfy
\eqref{II:eq:chart-segment} on the overlap. Their convex combination
does as well. It equals $Y_\theta$ on $|\theta|\ge2MB$.

The direction $Y_\theta$ is invariant under changes of the position
frame and simultaneous permutations of the bubble labels, so these
local fields agree under chart gluing. If $N=1$, then $\theta=0$ and
$W_{\rm ch}=W_{\rm seg}$.
\end{proof}

On the positive-coefficient terminal tubes, we also control the
logarithmic mass $P$. Denote
\[
l:=L_{\rm ray},
\qquad
Q:=\sum_i\rho_i^*q_i>0.
\]
Along each balanced ray,
\begin{equation*}
P(\gamma(l))=Ql+O(1),
\qquad
DP(\gamma(l))[\gamma'(l)]
=
Q+O(e^{-\eta l})
\end{equation*}
for some $\eta>0$. Indeed, the normalized radial mean of the $i$th
scale profile is
\begin{align*}
\frac1{2\pi\sigma_i}
\int_{\mathbb R^2}
e^{v_i}\left(m_i+Y\cdot\nabla v_i\right)\,dY
&=
m_i-2=
q_i.
\end{align*}
The core and exterior errors are exponentially small by
Theorems~\ref{II:thm:fixed-balanced-ray} and
\ref{II:thm:mixed-ray}.

For $u=\gamma(l)+h$ with $\|h\|<2e^{-4l}$,
\eqref{II:eq:ambient-moments} gives
\[
\|DP(u)\|\le C\sqrt l,
\qquad
\|D^2P(u)\|\le Cl.
\]
Since $\|S-\gamma'(l)\|\le C\|h\|$,
\begin{align*}
|DP(u)[S]-DP(\gamma(l))[\gamma'(l)]|
&\le
Cl\|h\|+C\sqrt l\,\|h\|=o(1).
\end{align*}
Hence, $DP(u)[S]=Q+o(1).$

On a positive-coefficient tube,
$V_-=-S-6h$. Therefore
\begin{align*}
DP(u)[V_-]
&=
-DP(u)[S]-6DP(u)[h]=
-Q+o(1)
\le-\frac Q2
\end{align*}
for all sufficiently large $l$. Since
$J_\kappa=J_{\rho^*}+\kappa P$,
\begin{align*}
-DJ_\kappa(u)[V_-]
&=
-DJ_{\rho^*}(u)[V_-]-\kappa DP(u)[V_-]\ge
c\Psi_l+\frac{\kappa Q}{2}.
\end{align*}
Thus $V_-$ decreases both $J_{\rho^*}$ and $J_\kappa$ on the overlap
with the exterior field.

\subsection{Proof of the main result}

\begin{proof}[Proof of Theorem~\ref{II:thm:characterization-cpi}]
Set
\[
\mathcal C_{\rho^*}
:=
\{u\in\mathcal H:DJ_{\rho^*}(u)=0\}.
\]
This set is compact by \cite[Corollary~1.1]{AFHZI}. For a balanced
critical pair $(p^*,r^*)$, denote
\[
c_*:=C(v,N)-\frac1N\ff(p^*).
\]
In the corresponding original modulation chart, let
\[
x:=|p-p^*|,
\qquad
X:=|r-r^*|+L|b|+\|w\|.
\]
For large $L$, define the switching region by
\begin{equation}\label{II:eq:switch-region}
x+X+|J_{\rho^*}(u)-c_*|<L^{-6}.
\end{equation}

We first locate $E_\kappa$. Along $\rho=(1-t)\rho^*$,
\[
\Lambda_{I,N}((1-t)\rho^*)
=
4t(1-t)\sum_i\sigma_i>0,
\qquad 0<t<1.
\]
Let $DJ_t(u)=0$, $t\to0$, and suppose that $u$ blows up. After extraction
and relabeling, \cite[Theorem~3.1, Proposition~3.2 and
Lemma~4.1]{AFHZI} gives a balanced Morse pair $(p^*,r^*)$ and
\[
x+X=o(L^{-k})
\qquad\text{for every fixed }k>0.
\]
Equation~\eqref{II:eq:coarse-energy-capture},
$\|u-U_{p,\ell,0}\|=O(\sqrt L\,|b|+\|w\|)$, and
$\|\nabla J_{\rho^*}\|=O(\sqrt L)$ on the joining segment give
\[
|J_{\rho^*}(u)-c_*|=o(L^{-k})
\qquad\text{for every fixed }k>0.
\]
Moreover, \cite[Theorem~4.1]{AFHZI} gives
\[
\frac{\Lambda_{I,N}((1-t)\rho^*)}
{\Psi_{\widehat m}(e^L)}
=
-\mathcal L_\infty(p^*,r^*)+o(1).
\]
The left-hand side is positive. Assumption~\ref{II:item:a2} therefore
implies
\[
\mathcal L_\infty(p^*,r^*)<0.
\]

Choose bounded neighborhoods $O_0,O_1$ of $\mathcal C_{\rho^*}$ such
that $\overline{O_0}\subset O_1$. If
$\mathcal C_{\rho^*}=\varnothing$, take $O_0=O_1=\varnothing$.
Choose one representative of each class in
$\mathcal K_\infty^-/\mathcal S_N$ and an ordered location cell
containing it. After shrinking these cells, the corresponding high-scale
charts are pairwise disjoint. Choose $L_*$ so that their parts with
$L>L_*$ are disjoint from $\overline{O_1}$.

For sufficiently small $\kappa>0$,
\[
E_\kappa\setminus O_0
\subset
\bigcup_{[(p^*,r^*)]\in\mathcal K_\infty^-/\mathcal S_N}
\left\{
L>L_*+2:
x+X+|J_{\rho^*}(u)-c_*|<\frac14L^{-6}
\right\}.
\]
Otherwise, there would exist $DJ_{t_k}(u_k)=0$, $t_k\to0$, violating
this inclusion. If $(u_k)$ is bounded, then $P(u_k)$ is bounded by
Moser--Trudinger, and
Lemma~\ref{II:lem:exterior-mass-compact} gives a limit in
$\mathcal C_{\rho^*}$. If $(u_k)$ is unbounded, the preceding blow-up
estimates place it in one of the displayed negative regions. Both
alternatives are contradictions.

The switching regions lie in the full bubbling charts. Indeed, balance
and smoothness give
\[
\Delta_\Omega(p,\ell,0)
=
\max_{i,t}|1-N\Omega_{i,t}(p,\ell,0)|
=
O(L^{-6}).
\]
For $\phi\in\mathcal H$ with $\|\phi\|_{\mathcal H}=1$, set
\[
c_{i,t}
:=
\frac{\int_MV_{i,t}\phi_i\,dV_g}
     {\int_MV_{i,t}\,dV_g}.
\]
The anchored moment estimates and the amplitude Gram bound give
\[
\|\phi_i-c_{i,t}\|_{L^2(V_{i,t}dV_g)}\le C,
\qquad
|c_{i,t}|\le C\sqrt L.
\]
The core comparison and the exponentially small exterior mass therefore
yield
\[
\left|DJ_{\rho^*}(U_{p,\ell,0})[\phi]\right|
\le
C\sqrt L\left(\Delta_\Omega+e^{-\eta L}\right).
\]
Hence
\[
\|\nabla J_{\rho^*}(U_{p,\ell,0})\|_{\mathcal H}
\le
C\sqrt L\left(L^{-6}+e^{-\eta L}\right).
\]
Since
\[
\|u-U_{p,\ell,0}\|_{\mathcal H}=O(L^{-6}),
\]
and
\[
\sup_{\xi\in[0,1]}
\left\|
D^2J_{\rho^*}
\bigl(U_{p,\ell,0}+\xi(u-U_{p,\ell,0})\bigr)
\right\|
\le CL,
\]
we obtain
\begin{align*}
\|\nabla J_{\rho^*}(u)\|_{\mathcal H}
&\le
C\sqrt L\left(L^{-6}+e^{-\eta L}\right)
+CL\|u-U_{p,\ell,0}\|_{\mathcal H}=
O(L^{-5}).
\end{align*}
Thus the switching regions lie in the full bubbling charts for all
sufficiently large $L$.

Conversely, the negative inner ray tubes lie in the smaller switching
regions for large $L$. The uniform modulation inverse gives
\[
|L-L_{\rm ray}|
=
O(e^{-4L_{\rm ray}}),
\qquad
x+X
=
O(\sqrt{L_{\rm ray}}\,e^{-4L_{\rm ray}}).
\]
Moreover,
\[
\left|
J_{\rho^*}(u)-J_{\rho^*}(\gamma(L_{\rm ray}))
\right|
=
O(\sqrt{L_{\rm ray}}\,e^{-4L_{\rm ray}}),
\]
whereas the balanced-ray expansion gives
\[
\left|
J_{\rho^*}(\gamma(L_{\rm ray}))-c_*
\right|
=
O(\Psi_{L_{\rm ray}}).
\]
Consequently,
\begin{align*}
x+X+\left|J_{\rho^*}(u)-c_*\right|
&\le
O(\sqrt{L_{\rm ray}}\,e^{-4L_{\rm ray}})
+O(\Psi_{L_{\rm ray}})=
o(L_{\rm ray}^{-6}),
\end{align*}
which proves the claimed inclusion.

We now construct the global field. Choose Lipschitz functions
$\eta_{\rm sw},\eta_{\rm ray},\beta:[0,\infty)\to[0,1]$ such that
\[
\eta_{\rm sw}=1\ \text{on }[0,1/3],
\qquad
\eta_{\rm sw}=0\ \text{on }[1,\infty),
\]
\[
\eta_{\rm ray}=1\ \text{on }[0,1],
\qquad
\eta_{\rm ray}=0\ \text{on }[9/4,\infty),
\]
and
\[
\beta=0\ \text{on }[0,L_*],
\qquad
\beta=1\ \text{on }[L_*+1,\infty).
\]
In each negative chart and ray tube, set
\begin{align*}
\chi
&:=
\beta(L)\,
\eta_{\rm sw}\left(
L^6\left[x+X+|J_{\rho^*}(u)-c_*|\right]
\right),\\
\zeta
&:=
\beta(L_{\rm ray})\,
\eta_{\rm ray}\left(
e^{8L_{\rm ray}}\|h\|^2
\right).
\end{align*}
Extend these functions by zero outside their coordinate neighborhoods.

Choose nested unions of full charts
\[
\overline V_0\subset V_1,
\qquad
\overline V_1\subset V_2,
\]
with successively relaxed margins in the location, separation,
scale-ratio, lower-scale, gradient, and modulation conditions. The
switching regions and signed tubes lie in $V_0$. On $V_2$, take the
field $W_{\rm ch}$ from Proposition~\ref{II:prop:chart-segment}. Choose
a locally Lipschitz function $\vartheta$ such that
\[
\vartheta=1\ \text{on }\overline V_0,
\qquad
\vartheta=0\ \text{on }\mathcal H\setminus V_1.
\]
Away from $E_\kappa$, define
\[
W_b:=\vartheta W_{\rm ch}+(1-\vartheta)W_{\rm ext}.
\]
Propositions~\ref{II:prop:chart-segment} and
\ref{II:prop:exterior-gluing} give constants $c_0,B_0>0$ such that
\begin{equation}\label{II:eq:common-crossing}
\|W_b\|\le B_0,
\qquad
\min\{-DJ_{\rho^*}[W_b],-DJ_\kappa[W_b]\}
\ge c_0d_\kappa.
\end{equation}

Near a negative pair, set
\[
W
=
\chi\bigl(\zeta V_++(1-\zeta)W_{\rho^*}\bigr)
+(1-\chi)W_b.
\]
Near a positive pair, set
\[
W=\zeta V_-+(1-\zeta)W_b.
\]
Near $\overline{O_0}$, use
\[
W_0
:=
-\frac{\nabla J_{\rho^*}}
{\sqrt{1+\|\nabla J_{\rho^*}\|^2}}
\]
and interpolate with $W_b$ inside $O_1$. Elsewhere, set $W=W_b$.

Every active summand strictly decreases $J_{\rho^*}$ away from
$\mathcal C_{\rho^*}$. The coefficient of $W_b$ vanishes near
$E_\kappa$, while $W_{\rho^*}$ is supported at large scale away from
$\mathcal C_{\rho^*}$. Thus $W$ is bounded and locally Lipschitz.
Consequently, its forward flow is complete.

On the part of $V_0$ where $|\theta|\ge2MB$ and the terminal cutoffs
vanish,
\[
W=W_b=W_{\rm ch}=Y_\theta
\]
by Lemma~\ref{II:lem:tau-flow}. On the chart--exterior overlap,
\eqref{II:eq:common-crossing} gives the hypotheses of
Corollary~\ref{II:cor:crossing-cost}. On each negative inner tube,
$\chi=\zeta=1$, so Theorem~\ref{II:thm:ambient-rays} gives forward
invariance. On each positive inner tube,
\[
\dot L_{\rm ray}=-1
\]
until the lower scale face is reached.

We next establish finite-time entrance into the terminal regions. Let
$\mathcal U$ be the union of $O_1$, the negative switching regions with
$L>L_*+1$, and a bounded neighborhood of the lower-scale transition
bands. Outside $\mathcal U$, both $J_{\rho^*}$ and $J_\kappa$ decrease.
As long as $J_{\rho^*}(u(s))>a$,
\[
P(u(s))
=
\frac{J_\kappa(u(s))-J_{\rho^*}(u(s))}{\kappa}
\le
\frac{J_\kappa(u(0))-a}{\kappa}.
\]
On the closed complement of $\mathcal U$ with this bound on $P$,
Lemma~\ref{II:lem:exterior-mass-compact} gives
\[
\inf d_\kappa>0.
\]
Indeed, otherwise a minimizing sequence would converge to a point of
$E_\kappa\subset\mathcal U$.

Where $V_-$ is active, the logarithmic mass estimate gives
\[
-DJ_\kappa[V_-]
=
-DJ_{\rho^*}[V_-]-\kappa DP[V_-]
\ge
c\Psi_{L_{\rm ray}}+\frac{\kappa Q}{2}.
\]
The upper bound on $P$ and
\[
P=QL_{\rm ray}+O(1)
\]
give an upper bound for $L_{\rm ray}$. Hence the descent of
$J_{\rho^*}$ is also bounded below there. Integration shows that every
trajectory either enters $\mathcal U$ in finite time or reaches the
sublevel $\{J_{\rho^*}\le a\}$.

We now consider bounded trajectories. For every $R<\infty$ and
$\delta>0$, there exists $m(R,\delta)>0$ such that
\begin{equation}\label{II:eq:bounded-descent}
-DJ_{\rho^*}(u)[W(u)]
\ge m(R,\delta)
\end{equation}
whenever
\[
\|u\|\le R,
\qquad
\operatorname{dist}_{\mathcal H}
(u,\mathcal C_{\rho^*})\ge\delta.
\]
If \eqref{II:eq:bounded-descent} failed, one active summand of the
glued field would have coefficient bounded below along a sequence
\((u_k)\) with vanishing descent.

Suppose first that this summand is \(W_b\). Then
\eqref{II:eq:common-crossing} gives
\[
d_\kappa(u_k)\to0.
\]
By Lemma~\ref{II:lem:exterior-mass-compact}, after extraction,
\[
u_k\to u_\infty\in E_\kappa.
\]
By construction, \(E_\kappa\) is contained in an open neighborhood on
which the coefficient of \(W_b\) vanishes. Hence this coefficient is
zero for all sufficiently large \(k\), contradicting its positive
lower bound.

For \(W_{\rho^*}\) and \(W_0\), Lemma~\ref{II:lem:exterior-mass-compact}
with $t=0$ gives a limit in \(\mathcal C_{\rho^*}\), contradicting
\(\operatorname{dist}_{\mathcal H}(u_k,\mathcal C_{\rho^*})\ge\delta\).
On a signed ray tube, the norm bound controls \(L_{\rm ray}\) from
above and the cutoff controls it from below; hence
\(c\Psi_{L_{\rm ray}}\) is bounded away from zero. This proves
\eqref{II:eq:bounded-descent}.

Let $B$ bound the speed of $W$. A bounded trajectory has energy bounded
below. If it reaches
\[
\operatorname{dist}_{\mathcal H}
(u,\mathcal C_{\rho^*})\ge2\delta,
\]
then for at least $\delta/B$ units of time its distance is at least
$\delta$. Equation~\eqref{II:eq:bounded-descent} gives
\[
J_{\rho^*}(u(s_0))
-
J_{\rho^*}(u(s_0+\delta/B))
\ge
\frac{\delta}{B}m(R,\delta).
\]
Such visits can occur only finitely many times. Therefore
\[
\operatorname{dist}_{\mathcal H}
(u(s),\mathcal C_{\rho^*})
\to0.
\]
Since $\mathcal C_{\rho^*}$ is compact, the trajectory is relatively
compact. If $\mathcal C_{\rho^*}=\varnothing$, bounded trajectories do
not exist.

Let now $u(s)$ be unbounded and suppose
\[
J_{\rho^*}(u(s))\downarrow c>-\infty.
\]
Choose $R>1$ such that
\[
\|v\|\le R-1
\qquad\text{for every }v\in\mathcal C_{\rho^*}.
\]
Each crossing of
\[
R\le\|u\|\le2R
\]
has distance at least one from $\mathcal C_{\rho^*}$ and costs at least $\frac{R}{B}m(2R,1)$
in energy. Since the trajectory has arbitrarily late visits above
$2R$, it eventually remains above $R$. Repeating this argument with
arbitrarily large $R$ gives
\[
\|u(s)\|\to\infty.
\]
Applying the finite-time entrance property from arbitrarily late times
shows that the trajectory visits negative switching regions
arbitrarily late. If $\mathcal K_\infty^-=\varnothing$, no such
trajectory exists.

Since $\mathcal K_\infty^-/\mathcal S_N$ is finite, one negative class,
represented by $(p^*,r^*)$, is visited arbitrarily late. Along these
visits, $L\to\infty$ and \eqref{II:eq:switch-region} gives $c=c_*$.
Choose $a_0,b_0$ as in Corollary~\ref{II:cor:crossing-cost}, with
$2a_0+2b_0<h$, and choose the corresponding closed box disjoint from
all other terminal neighborhoods. Let
\[
S_*
:=
\left\{
x+X+|J_{\rho^*}(u)-c_*|\le L^{-6}
\right\}
\]
inside this box at large scale. On $S_*$,
\[
d_\kappa(u)
\le
\|\nabla J_{\rho^*}(u)\|
=
O(L^{-5}).
\]

Consider an excursion starting at the last contact time $s_0$ with
$S_*$ and ending before the next contact or the first exit from the
box. Set $L_0:=L(s_0)$. Outside $S_*$, the field equals $W_b$, so
\eqref{II:eq:common-crossing} gives
\begin{align*}
c_0\int_{s_0}^s d_\kappa(u(\tau))\,d\tau
&\le
J_{\rho^*}(u(s_0))-J_{\rho^*}(u(s))
\le L_0^{-6},\\
\kappa\bigl(P(u(s))-P(u(s_0))\bigr)
&\le
J_{\rho^*}(u(s_0))-J_{\rho^*}(u(s))
\le L_0^{-6}.
\end{align*}
Here $J_{\rho^*}(u(s))\ge c_*$ because $c_*$ is the limiting energy. In
this box,
\[
P(u)=QL+O(1),
\qquad Q>0.
\]
The second inequality therefore gives
\[
L(s)\le H:=L_0+C
\]
throughout the excursion.

A first exit from the box must cross either
\[
a_0\le x\le2a_0
\qquad\text{or}\qquad
b_0\le X\le2b_0.
\]
Corollary~\ref{II:cor:crossing-cost} gives an energy loss of at least
\[
c_{a_0,b_0}
\qquad\text{or}\qquad
c_{a_0,b_0}H^{-5/2},
\]
respectively. Both quantities exceed $L_0^{-6}$ for large $L_0$.
Hence no sufficiently late excursion leaves the box.

Letting first $b_0\to0$ and then $a_0\to0$ gives
\[
L\to\infty,
\qquad
p\to p^*,
\qquad
r\to r^*,
\qquad
L|b|+\|w\|\to0.
\]
The first limit follows from $\|u(s)\|\to\infty$, since bounded scales
give bounded bubble families.

It remains to control the gradient cutoff. Set
\[
d(s):=d_\kappa(u(s)).
\]
On each excursion,
\[
\sup_{0\le t\le\kappa}\|D^2J_t(u(s))\|\le CH,
\]
and the speed is bounded. Hence $d(s)$ is $CH$-Lipschitz. Therefore
\[
\frac{d}{ds}d(s)^2\le2CHd(s)
\]
almost everywhere. Using the integral estimate above,
\[
d(s)^2
\le
d(s_0)^2
+
2CH\int_{s_0}^s d(\tau)\,d\tau
=
O(L_0^{-5}).
\]
Since $L(s)\le L_0+C$,
\[
d_\kappa(u(s))=O(L^{-5/2})
\]
along the tail.

Let $t(s)\in[0,\kappa]$ attain the minimum defining
$d_\kappa(u(s))$. Equation~\eqref{II:eq:collar-tilt} and
$L|b|\to0$ give
\[
\|w\|+|y|=o(L^{-1}),
\qquad
t(s)=o(L^{-2}).
\]
Since
\[
y=t(s)\mathbf1+(1-t(s))\theta,
\]
we obtain
\[
L|\theta|\to0.
\]
Finally, $\|DP\|\le C\sqrt L$ gives
\begin{equation}\label{II:eq:untilted-gradient}
\|\nabla J_{\rho^*}(u(s))\|
\le
d_\kappa(u(s))+t(s)\|DP(u(s))\|
=
o(L^{-3/2}).
\end{equation}

The amplitude Gram estimate gives
\[
\|u-U_{p,\ell,0}\|
\le
C\sqrt L\,|b|+\|w\|
\to0.
\]
Together with \eqref{II:eq:untilted-gradient} and the strict separation
and scale-ratio margins, this places the tail in a full chart. Moreover,
\[
\left(
p,r,\sqrt L\,b^0,(c^\alpha)_{\alpha=1}^{N-1},
\theta,w,e^{-L}
\right)
\to
\left(p^*,r^*,0,0,0,0,0\right).
\]
Thus every unbounded forward orbit with energy bounded below determines
a balanced critical pair $(p^*,r^*)$ satisfying
\[
\mathcal L_\infty(p^*,r^*)<0.
\]
Two such orbits represent the same critical point at infinity precisely
when their limiting pairs agree after a simultaneous permutation.
Conversely, Theorem~\ref{II:thm:ambient-rays} provides a
forward-invariant open tube of escaping trajectories for every negative
pair. Hence the critical points at infinity are in bijection with
\[
\mathcal K_\infty^-/\mathcal S_N.
\]
The necessity argument uses the chart field only through bounded speed,
strict descent, and \eqref{II:eq:common-crossing}, which also proves the
assertion concerning admissible replacements.
\end{proof}

The proof also gives the following quantitative asymptotics.

\begin{corollary}
Every unbounded forward trajectory $u(s)$ with
\[
\inf_{s\ge0}J_{\rho^*}(u(s))>-\infty
\]
satisfies
\[
L(s)|\theta(s)|\to 0 ,
\qquad
L(s)^{3/2}
\|\nabla J_{\rho^*}(u(s))\|
\to 0.
\]
Equivalently,
\[
|\theta(s)|=o(L(s)^{-1}),
\qquad
\|\nabla J_{\rho^*}(u(s))\|
=o(L(s)^{-3/2}).
\]
\end{corollary}

\begin{proof}
These estimates follow from the proof of 
\eqref{II:eq:untilted-gradient}.
\end{proof}

The following example shows that the pseudogradient need not form a
uniform positive angle with the ambient gradient.

Let $N=1$, $m_i=4$, and let $M$ be the flat square torus. Take
\[
A=
\begin{pmatrix}
0&1\\
1&0
\end{pmatrix},
\qquad
\sigma_1=\sigma_2=4,
\qquad
h_1=e^{B\phi+f},
\qquad
h_2=e^{B\phi-f},
\]
where $\phi(x)=\cos(2\pi x_1)+2\cos(2\pi x_2)$ and
$f(x)=\sin(2\pi x_1)$. Up to an additive constant, $\ff$ is a positive
multiple of $B\phi$. It is Morse, with four critical points satisfying
$\Delta\phi\ne0$.

Since $M$ is flat, $R^g$ is constant,
$\nabla_xH^g(p,p)=0$, and $\Delta_xH^g(p,p)=1$. At every critical point,
\[
T_{i,1}^{\rm raw}
=
B\Delta\phi+|\nabla f|^2+8\pi,
\qquad i=1,2.
\]
Thus, for sufficiently large $B>0$, $\mathcal L_\infty$ is nonzero at
all four critical points and
$\mathcal L_\infty(p^*)<0$ at $p^*:=(0,0)$. The matrix assumptions hold,
and balance is automatic because $N=1$.

Consider the uncorrected terminal ray $\gamma(l)$ at $p^*$. Choose a
nonzero radial function $\eta\in C_c^\infty(B_1)$ with $\eta\ge0$, and
set $g:=\nabla f(p^*)\ne0$ and $Y:=e^ly_{p^*}$. Define
$(T_l)_1:=\eta(Y)g\cdot Y$ and $(T_l)_2:=0$. Its surface mean vanishes,
and $c\le\|T_l\|_{\mathcal H}\le C$.

The free density is radial on $\operatorname{supp}T_l$, so its pairing
with $T_l$ vanishes. By
\eqref{II:eq:m4-normalized-density-error}, with $s=1$ and
$\Omega_{1,1}=1$, the first-component probability density is
\[
k_1(Y)\left[1+e^{-l}g\cdot Y+O(le^{-2l})\right]\,dY.
\]
Hence
\begin{align*}
DJ_{\rho^*}(\gamma(l))[T_l]
&=
-\rho_1^*e^{-l}
\int_{\mathbb R^2}
k_1(Y)\eta(Y)(g\cdot Y)^2\,dY
+O(le^{-2l})\le-ce^{-l}.
\end{align*}
Therefore
$\|\nabla J_{\rho^*}(\gamma(l))\|\ge ce^{-l}$.

On this ray,
$W(\gamma(l))=\gamma'(l)$ and
$c\le\|\gamma'(l)\|\le C$. Moreover,
$$-DJ_{\rho^*}(\gamma(l))[\gamma'(l)]\asymp le^{-2l}.$$ Consequently,
\[
0\le
\frac{-DJ_{\rho^*}(\gamma(l))[W(\gamma(l))]}
{\|\nabla J_{\rho^*}(\gamma(l))\|\,\|W(\gamma(l))\|}
\le Cle^{-l}\to 0.
\]

	\medskip
	\noindent\textbf{Acknowledgements.}~
	During the preparation of this manuscript, the authors used ChatGPT (OpenAI, GPT-6.1 Sol and GPT-5.6 Sol)  to assist with organization, LaTeX consistency, proof
	checking, and the formulation of  arguments. The authors subsequently reviewed, revised, and verified all content and take full responsibility for the accuracy, originality, and integrity of the manuscript.
	
	\section*{Statements and Declarations}
	\noindent\textbf{Funding.}~
	M. Ahmedou is supported by the Deutsche Forschungsgemeinschaft (DFG) under grant AH~156/2-1. Z. Hu is partially supported by the
	National Natural Science Foundation of China
	(Grant No.~12601383) and the China Postdoctoral Science Foundation
	(Grant No.~2026M793392).  L. Zhang is partially supported by Simons Foundation Grant SFI-MPS-TSM-00013752.

	\medskip
	
	\noindent\textbf{Data availability.}~	 Data sharing does not apply to this article as no datasets were generated or analysed during the current study.

	\medskip
	
	\noindent\textbf{Competing interests.} ~ The authors declare that there is no conflict of interest.

\begingroup
\small
\begin{thebibliography}{99}
\setlength{\itemsep}{2pt}
\setlength{\parsep}{0pt}

\bibitem{AhmedouBenAyed2017}
M.~Ahmedou and M.~Ben Ayed,
Theory of ``critical points at infinity'' and a resonant singular
Liouville-type equation,
\emph{Adv. Nonlinear Stud.} 17 (2017),
\href{https://doi.org/10.1515/ans-2016-6016}{139--166}.

\bibitem{ABAMorse}
M.~Ahmedou and M.~Ben Ayed,
Morse inequalities at infinity for a resonant mean field equation,
\emph{Commun. Contemp. Math.} 25 (2023),
\href{https://doi.org/10.1142/S0219199721500541}{Paper No.~2150054}.

\bibitem{MBL2017}
M.~Ahmedou, M.~Ben Ayed, and M.~Lucia,
On a resonant mean field type equation:
a ``critical point at infinity'' approach,
\emph{Discrete Contin. Dyn. Syst.} 37 (2017),
\href{https://doi.org/10.3934/dcds.2017075}{1789--1818}.

\bibitem{AFHZI}
M.~Ahmedou, J.~Fey, Z.~Hu, and L.~Zhang,
Liouville systems I: Blow-up analysis and orthogonal decomposition,
\href{https://arxiv.org/abs/2611.07857}{arXiv:2611.07857} (2026).

\bibitem{AHW}
M.~Ahmedou, Z.~Hu, and H.~Wang,
A degree-counting formula for a Keller--Segel equation on a surface
with boundary,
\emph{Commun. Contemp. Math.} (2026),
\href{https://doi.org/10.1142/S0219199726500392}
{Paper No.~2650039}.

\bibitem{AHWZ}
M.~Ahmedou, Z.~Hu, H.~Wang, and L.~Zhang,
Morse theory and topological degree for resonant mean field equations
on Riemann surfaces with boundary,
\href{https://arxiv.org/abs/2610.41760}{arXiv:2610.41760} (2026).

\bibitem{Bahri1989}
A.~Bahri,
\emph{
{Critical points at infinity in some variational problems}},
Pitman Research Notes in Mathematics Series 182,
Longman Scientific \& Technical, Harlow, 1989.

\bibitem{BahriCoron1988}
A.~Bahri and J.-M.~Coron,
On a nonlinear elliptic equation involving the critical Sobolev exponent:
the effect of the topology of the domain,
\emph{Comm. Pure Appl. Math.} 41 (1988),
\href{https://doi.org/10.1002/cpa.3160410302}{253--294}.

\bibitem{BahriLiRey1995}
A.~Bahri, Y.~Li, and O.~Rey,
On a variational problem with lack of compactness:
the topological effect of the critical points at infinity,
\emph{Calc. Var. Partial Differential Equations} 3 (1995),
\href{https://doi.org/10.1007/BF01190892}{67--93}.

\bibitem{BM1991}
H.~Br\'ezis and F.~Merle,
Uniform estimates and blow-up behavior for solutions of
$-\Delta u=V(x)e^u$ in two dimensions,
\emph{Comm. Partial Differential Equations} 16 (1991),
\href{https://doi.org/10.1080/03605309108820797}{1223--1253}.

\bibitem{CK1995}
S.~Chanillo and M.~K.-H.~Kiessling,
Conformally invariant systems of nonlinear PDE of Liouville type,
\emph{Geom. Funct. Anal.} 5 (1995),
\href{https://doi.org/10.1007/BF01902215}{924--947}.

\bibitem{CLZconstruction}
Z.~Cheng, H.~Li, and L.~Zhang,
Construction of blow-up solutions for Liouville systems,
\emph{Proc. Lond. Math. Soc.} (3) 131 (2025),
\href{https://doi.org/10.1112/plms.70093}{Paper No.~e70093}.

\bibitem{CLZ}
Z.~Cheng, H.~Li, and L.~Zhang,
Local uniqueness and non-degeneracy of blow-up solutions for regular
Liouville systems,
\href{https://arxiv.org/abs/2509.09781}{arXiv:2509.09781} (2025).

\bibitem{CSW1997}
M.~Chipot, I.~Shafrir, and G.~Wolansky,
On the solutions of Liouville systems,
\emph{J. Differential Equations} 140 (1997),
\href{https://doi.org/10.1006/jdeq.1997.3316}{59--105};
erratum, 178 (2002),
\href{https://doi.org/10.1006/jdeq.2001.4105}{630}.

\bibitem{Fontana}
L.~Fontana,
Sharp borderline Sobolev inequalities on compact Riemannian manifolds,
\emph{Comment. Math. Helv.} 68 (1993),
\href{https://doi.org/10.1007/BF02565828}{415--454}.

\bibitem{HL2009}
J.~Hor\'ak and M.~Lucia,
A minimax theorem in the presence of unbounded Palais--Smale sequences,
\emph{Israel J. Math.} 172 (2009),
\href{https://doi.org/10.1007/s11856-009-0067-0}{125--143}.

\bibitem{HuangZhang}
H.-Y.~Huang and L.~Zhang,
On Liouville systems at critical parameters, Part 2: Multiple bubbles,
\emph{Calc. Var. Partial Differential Equations} 61 (2022),
\href{https://doi.org/10.1007/s00526-021-02110-6}{Paper No.~3}.

\bibitem{LinZhang2010}
C.-S.~Lin and L.~Zhang,
Profile of bubbling solutions to a Liouville system,
\emph{Ann. Inst. H. Poincar\'e Anal. Non Lin\'eaire} 27 (2010),
\href{https://doi.org/10.1016/j.anihpc.2009.09.001}{117--143}.

\bibitem{LinZhang2011}
C.-S.~Lin and L.~Zhang,
A topological degree counting for some Liouville systems of mean field type,
\emph{Comm. Pure Appl. Math.} 64 (2011),
\href{https://doi.org/10.1002/cpa.20355}{556--590}.

\bibitem{LZone}
C.-S.~Lin and L.~Zhang,
On Liouville systems at critical parameters, Part 1: One bubble,
\emph{J. Funct. Anal.} 264 (2013),
\href{https://doi.org/10.1016/j.jfa.2013.02.022}{2584--2636}.

\bibitem{Lucia2007}
M.~Lucia,
A deformation lemma with an application to a mean field equation,
\emph{Topol. Methods Nonlinear Anal.} 30 (2007),
\href{https://www.tmna.ncu.pl/static/files/v30n1-05.pdf}{113--138}.

\end{thebibliography}

\endgroup
\end{document}